\documentclass[
  smallextended,
  final, %final, draft, referee
  envcountsame,
  natbib
  ]{svjour3}
\pdfoutput=1

\title{Shallow water equations: Split-form, entropy stable, well-balanced, and
       positivity preserving numerical methods}
\titlerunning{Entropy stable, well-balanced, positivity preserving methods for
              shallow water}
\author{Hendrik Ranocha}
\institute{
  Hendrik Ranocha \at
  TU Braunschweig, Pockelsstra\ss e 14, 38106 Braunschweig, Germany\\
  Tel.: +49-531-391-7417, Fax: +49-531-391-7409\\
  \email{h.ranocha@tu-bs.de}
}

\date{Received: date / Accepted: date}

\usepackage[utf8]{luainputenc}
\usepackage[english]{babel}
\usepackage{csquotes}

\usepackage[plainpages=false,pdfpagelabels,hidelinks]{hyperref}
\hypersetup{pdfauthor={Hendrik Ranocha}}
\hypersetup{pdftitle={Shallow water equations: Split-form, entropy stable,
                      well-balanced, and positivity preserving numerical methods}}

\usepackage{amsmath}
\allowdisplaybreaks
\usepackage{amssymb}
\usepackage{commath}
\usepackage{mathtools}

\usepackage{siunitx}

\usepackage{color}
\usepackage{graphicx}
\usepackage[small]{caption}
\usepackage{subcaption}

\begingroup\expandafter\expandafter\expandafter\endgroup
\expandafter\ifx\csname pdfsuppresswarningpagegroup\endcsname\relax
\else
\fi

\usepackage{booktabs}
\usepackage{rotating}

\usepackage{multicol}
\usepackage{enumitem}

\usepackage{lineno}
\newcommand*\patchAmsMathEnvironmentForLineno[1]{%
  \expandafter\let\csname old#1\expandafter\endcsname\csname #1\endcsname
  \expandafter\let\csname oldend#1\expandafter\endcsname\csname end#1\endcsname
  \renewenvironment{#1}%
     {\linenomath\csname old#1\endcsname}%
     {\csname oldend#1\endcsname\endlinenomath}}% 
\newcommand*\patchBothAmsMathEnvironmentsForLineno[1]{%
  \patchAmsMathEnvironmentForLineno{#1}%
  \patchAmsMathEnvironmentForLineno{#1*}}%
\AtBeginDocument{%
\patchBothAmsMathEnvironmentsForLineno{equation}%
\patchBothAmsMathEnvironmentsForLineno{align}%
\patchBothAmsMathEnvironmentsForLineno{flalign}%
\patchBothAmsMathEnvironmentsForLineno{alignat}%
\patchBothAmsMathEnvironmentsForLineno{gather}%
\patchBothAmsMathEnvironmentsForLineno{multline}%
}
\usepackage{calc}
\usepackage{xparse}

\newcommand{\widesmall}[2]{#2}

\renewcommand{\vec}[1]{\underline{#1}}
\NewDocumentCommand{\mat}{mo}{%
  \IfValueTF{#2}{%
    \underline{\underline{#1}}{#2}
  }{%
    \underline{\underline{#1}}\,
  }%
}
\newcommand{\diag}[1]{\operatorname{diag}\left(#1\right)}

\newcommand{\I}{\operatorname{I}}
\newcommand{\fnum}{f^{\mathrm{num}}}
\newcommand{\vecfnum}{\vec{f}^{\mathrm{num}}}
\newcommand{\fvol}{f^{\mathrm{vol}}}

\newcommand{\VOL}{\vec{\mathrm{VOL}}}
\newcommand{\SURF}{\vec{\mathrm{SURF}}}

\renewcommand{\epsilon}{\varepsilon}
\renewcommand{\phi}{\varphi}

\newcommand{\R}{\mathbb{R}}

\newsavebox{\DelimiterBox}
\newlength{\DelimiterHeight}
\newlength{\DelimiterDepth}
\newsavebox{\ArgumentBox}
\newlength{\ArgumentHeight}
\newlength{\ArgumentDepth}
\newlength{\ResizedDelimiterHeight}

\newcommand{\mean}[1]{%
  \savebox{\ArgumentBox}{$ \displaystyle #1 $}%
  \settoheight{\ArgumentHeight}{\usebox{\ArgumentBox}}%
  \settodepth{\ArgumentDepth}{\usebox{\ArgumentBox}}%
  \savebox{\DelimiterBox}{$ \{\!\!\{ $}%
  \settoheight{\DelimiterHeight}{\usebox{\DelimiterBox}}%
  \settodepth{\DelimiterDepth}{\usebox{\DelimiterBox}}%
  \setlength{\ResizedDelimiterHeight}{\maxof{1.2\ArgumentHeight}{\DelimiterHeight}}
  \!\!
  \resizebox{\width}{\ResizedDelimiterHeight}{ $\{\!\!\{$ }
  \mkern-6.5mu
  #1
  \mkern-6.5mu
  \resizebox{\width}{\ResizedDelimiterHeight}{ $\}\!\!\}$ }
  \!\!
}

\newcommand{\jump}[1]{%
  \savebox{\ArgumentBox}{$ \displaystyle #1 $}%
  \settoheight{\ArgumentHeight}{\usebox{\ArgumentBox}}%
  \settodepth{\ArgumentDepth}{\usebox{\ArgumentBox}}%
  \savebox{\DelimiterBox}{$[\![$}%
  \settoheight{\DelimiterHeight}{\usebox{\DelimiterBox}}%
  \settodepth{\DelimiterDepth}{\usebox{\DelimiterBox}}%
  \setlength{\ResizedDelimiterHeight}{\maxof{1.2\ArgumentHeight}{\DelimiterHeight}}
  \!\!
  \resizebox{\width}{\ResizedDelimiterHeight}{ $[\![$ }
  \mkern-6.5mu
  #1
  \mkern-6.5mu
  \resizebox{\width}{\ResizedDelimiterHeight}{ $]\!]$ }
  \!\!
}

\begin{document}

\maketitle

\begin{abstract}
  
For the first time, a general two-parameter family of entropy conservative
numerical fluxes for the shallow water equations is developed and investigated.
These are adapted to a varying bottom topography in a well-balanced way, i.e.
preserving the lake-at-rest steady state.
Furthermore, these fluxes are used to create entropy stable and well-balanced
split-form semidiscretisations based on general summation-by-parts (SBP) operators,
including Gauß nodes.
Moreover, positivity preservation is ensured using the framework of Zhang and Shu
(\emph{Maximum-principle-satisfying and positivity-preserving high-order schemes for
conservation laws: survey and recent developments}, 2011. In: Proceedings of the
Royal Society of London A: Mathematical, Physical and Engineering Sciences, The
Royal Society, vol 467, pp. 2752--2766).
Therefore, the new two-parameter family of entropy conservative fluxes is
enhanced by dissipation operators and investigated with respect to positivity
preservation. Additionally, some known entropy stable and positive numerical
fluxes are compared.
Furthermore, finite volume subcells adapted to nodal SBP bases with diagonal
mass matrix are used.
Finally, numerical tests of the proposed schemes are performed and some conclusions
are presented.

\end{abstract}

\section{Introduction}
\label{sec:introduction}

For the first time, a two-parameter family of entropy conservative and well-balanced
numerical fluxes for the shallow water equations is investigated, resulting in
genuinely high-order semidiscretisations that are both entropy stable and
positivity preserving.
These semidiscretisations of the shallow water equations in one space dimension
are based on \emph{summation-by-parts} (SBP) operators, see
inter alia the review articles by \citet{svard2014review, fernandez2014review}
and references cited therein.
This setting of SBP operators originates in the \emph{finite difference} (FD)
setting, but can also be used in polynomial methods as nodal \emph{discontinuous
Galerkin} (DG) \citep{gassner2013skew} or flux reconstruction / correction procedure
via reconstruction \citep{ranocha2016summation}.

Entropy stability has long been known as a desirable stability property for
conservation laws. Here, the semidiscrete setting of \citet{tadmor1987numerical,
tadmor2003entropy} will be used.
Other desirable stability properties for the shallow water equations are the
preservation of non-negativity of the water height and the correct handling
of steady states, especially the lake-at-rest initial condition, resulting in
well-balanced methods.

References for numerical methods for the shallow water equations can be found
in the review article of \citet{xing2014survey} and references cited therein.

This article extends the entropy conservative split form of \citet{gassner2016well,
wintermeyer2016entropy} to a new two-parameter family of well-balanced and entropy
conservative splittings. Moreover, SBP bases not including boundary nodes are
considered and corresponding semidiscretisations are developed in one spatial
dimension. To the author's knowledge, this general two-parameter family of
numerical fluxes has not appeared in the literature before.

Furthermore, positivity preservation in the framework of \citet{zhang2011maximum}
is considered. Therefore, the new fluxes are investigated regarding positivity
and \emph{finite volume} (FV) subcells are introduced. Additionally, some known
entropy stable and positivity preserving numerical fluxes are compared.

The extension to general SBP bases can be seen as some kind of 'positive negative
result': It is possible to get the desired properties of the schemes using Gauß
nodes, including a higher accuracy compared to Lobatto nodes, but this comes at
the cost of complicated correction terms. Therefore, although a similar extension
to two-dimensional unstructured and curvilinear grids can be conjectured to exist,
it is expected to be even more complex and thus not suited for high-performance
production codes.

However, the combination of entropy stability and positivity preservation for the
shallow water equations using high-order SBP methods has not been considered
before. The results can be expected to be extendable to high-performance codes
using Lobatto nodes and the flux-differencing form of \citet{fisher2013high}
on two-dimensional unstructured and curvilinear grids if tensor product bases
on quadrilaterals are used.

At first, some analytical properties of the shallow water equations are reviewed
in section \ref{sec:SWE} and the existing split form of \citet{gassner2016well,
wintermeyer2016entropy} is described in section \ref{sec:existing-split-form}.
Afterwards, a new two-parameter family of entropy conservative numerical fluxes for
the shallow water equations with constant bottom topography is developed in
section \ref{sec:fnum} and extended to a varying bottom in section \ref{sec:source}.
The corresponding semidiscretisation using general SBP bases is designed in
section \ref{sec:gauss}. The positivity preserving framework of \citet{zhang2011maximum}
is introduced to this setting and numerical
fluxes based on the entropy conserving schemes are investigated with respect to
positivity preservation in section \ref{sec:fluxes}. Additionally, some known fluxes
are presented. An extension of the idea to use FV subcells to the setting of
nodal SBP bases with diagonal mass matrix is proposed in section \ref{sec:subcell-FV}
and numerical experiments are presented in section \ref{sec:numerical-tests}.
Finally, the results are summed up in section \ref{sec:summary} and some
conclusions and directions of further research are presented.

\section{Review of some properties of the shallow water equations}
\label{sec:SWE}

The shallow water equations in one space dimension are
\begin{equation}
\label{eq:SWE}
\begin{aligned}
  \partial_t
  \underbrace{
  \begin{pmatrix}
    h
    \\
    h v
  \end{pmatrix}
  }_{= u}
  + \,\partial_x
  \underbrace{
  \begin{pmatrix}
    h v
    \\
    h v^2 + \frac{1}{2} g h^2
  \end{pmatrix}
  }_{= f(u)}
  =
  \underbrace{
  \begin{pmatrix}
    0
    \\
    - g h \partial_x b
  \end{pmatrix}
  }_{=s(u,x)},
\end{aligned}
\end{equation}
where $h$ is the water height, $v$ its speed, $hv$ the discharge, $b$ describes the bottom
topography, and $g$ is the gravitational constant. In the following, some well
known results that will be used in the remainder of this work are presented.

As described inter alia by \citet[Section 3.2]{bouchut2004nonlinear},
\citet[Section 3.3]{dafermos2010hyperbolic},
\citet{fjordholm2011well, wintermeyer2016entropy}, the entropy / total energy
$U = \frac{1}{2} h v^2 + \frac{1}{2} g h^2 + g h b
= \frac{1}{2} \frac{u_2^2}{u_1} + \frac{1}{2} g u_1^2 + g u_1 b$,
is strictly convex for positive water heights $h > 0$. Thus, in this case, the
associated entropy variables
\begin{equation}
\label{eq:entropy-variables}
  w
  =
  U'(u)
  =
  \begin{pmatrix}
    g (h + b) - \frac{1}{2} v^2
    \\
    v
  \end{pmatrix}
\end{equation}
and conserved variables $u$ can be used interchangeably. With the entropy flux
\begin{equation}
\label{eq:entropy-flux}
  F
  =
  \frac{1}{2} h v^3 + g h^2 v
  + g b h v,
\end{equation}
smooth solutions satisfy $\partial_t U + \partial_x F = 0$, and the entropy
inequality $\partial_t U + \partial_x F \leq 0$ will be used as an additional
admissibility criterion for weak solutions.

If the bottom topography $b \equiv 0$ is constant, the entropy variables are
\begin{equation}
\label{eq:entropy-variables-b-constant}
  w
  =
  \begin{pmatrix}
    g h - \frac{1}{2} v^2
    \\
    v
  \end{pmatrix}.
\end{equation}
In this case, the flux expressed in terms of the entropy variables
($b \equiv 0$) is
\begin{equation}
\label{eq:g}
\begin{aligned}
  f\left( u(w) \right)
  =&
  \begin{pmatrix}
    \frac{w_1 + \frac{1}{2} w_2^2}{g} w_2
    \\
    \frac{w_1 + \frac{1}{2} w_2^2}{g} w_2^2
    + \frac{1}{2} g \left( \frac{w_1 + \frac{1}{2} w_2^2}{g} \right)^2
  \end{pmatrix}
  =
  \frac{1}{g}
  \begin{pmatrix}
    w_1 w_2 + \frac{1}{2} w_2^3
    \\
    \frac{1}{2} w_1^2 + \frac{3}{2} w_1 w_2^2 + \frac{5}{8} w_2^4
  \end{pmatrix},
\end{aligned}
\end{equation}
and the flux potential ($b \equiv 0$) is given by
\begin{equation}
\label{eq:flux-potential}
\begin{aligned}
  \psi
  =&
  \frac{1}{2 g} w_1^2 w_2 + \frac{1}{2 g} w_1 w_2^3  + \frac{1}{8 g} w_2^5
  =
  \frac{1}{2} g h^2 v,
\end{aligned}
\end{equation}
fulfilling $\psi'(w) = f\left( u(w) \right)$. Finally, the entropy Jacobian
\begin{equation}
\label{eq:entropy-jacobian}
  \partial_w u
  =
  \left( \partial_u w \right)^{-1}
  =
  \begin{pmatrix}
    \frac{1}{g} & \frac{v}{g}
    \\
    \frac{v}{g} & h + \frac{v^2}{g}
  \end{pmatrix},
  \quad
  \partial_u w
  =
  \begin{pmatrix}
    g + \frac{u_2^2}{u_1^3} & - \frac{u_2}{u_1^2}
    \\
    - \frac{u_2}{u_1^2}     & \frac{1}{u_1}
  \end{pmatrix}
  =
  \begin{pmatrix}
    g + \frac{v^2}{h} & - \frac{v}{h}
    \\
    - \frac{v}{h}     & \frac{1}{h}
  \end{pmatrix},
\end{equation}
can be expressed by using a scaling of the eigenvectors in the form proposed
by \citet[Theorem 4]{barth1999numerical} as
\begin{equation}
\label{eq:scaling-barth}
  \partial_w u
  =
  R R^T,
  \qquad
  R = 
  \frac{1}{\sqrt{2 g}}
  \begin{pmatrix}
    1 & 1
    \\
    v - \sqrt{g h} & v + \sqrt{g h}
  \end{pmatrix},
\end{equation}
where the columns of $R$ are eigenvectors of the flux Jacobian $f'(u)$. This scaling
has also been used inter alia by \citet[Section 2.3]{fjordholm2011well}.

\section{Review of an existing split-form SBP method}
\label{sec:existing-split-form}

In order to fix some notation, present the general setting and motivate the
extensions developed in this work, some existing results will be reviewed.

A general SBP SAT semidiscretisation is obtained by a partition of the domain into
disjoint elements. On each element, the solution is represented in some basis,
mostly nodal bases. These cells are mapped to a standard element for the following
computations. There, the symmetric and positive definite mass matrix $\mat{M}$
induces a scalar product, approximating the $L_2$ scalar product. The derivative
is represented by the matrix $\mat{D}$. Interpolation to the (two point) boundary
of the cell (interval) is performed via the restriction operator $\mat{R}$ and
evaluation of the values at the right boundary minus values at the left boundary
is conducted by the boundary matrix $\mat{B} = \diag{-1,1}$. Together, these
operators fulfil the summation-by-parts (SBP) property
\begin{equation}
\label{eq:SBP}
  \mat{M} \mat{D} + \mat{D}[^T] \mat{M} = \mat{R}[^T] \mat{B} \mat{R},
\end{equation}
mimicking integration by parts on a discrete level
$\int_{\Omega} u (\partial_x v) + \int_{\Omega} (\partial_x u) v
= u v \big|_{\partial\Omega}$.
Here, the notation of \citet{ranocha2016summation, ranocha2015extended} has been
used. Then, similar to strong form discontinuous Galerkin methods, the
semidiscretisation can be written as the sum of volume terms, surface terms,
and numerical fluxes at the boundaries.

\citet{gassner2016well} proposed as semidiscretisation of the shallow water
equations \eqref{eq:SWE} with continuous bottom topography $b$ in the setting of
a discontinuous Galerkin spectral element method (DGSEM) using Lobatto-Legendre
nodes in each element, that can be generalised to diagonal norm SBP operators
with nodal bases including boundary nodes. \citet{wintermeyer2016entropy}
extended this setting to two space dimensions, curvilinear grids and discontinuous
bottom topographies. In one space dimension on a linear grid, this semidiscretisation
can be written as
\begin{equation}
\label{eq:semidisc-gassner}
\begin{aligned}
  \partial_t \vec{h}
  =&
  - \mat{D} \vec{hv}
  - \mat{M}[^{-1}] \mat{R}[^T] \mat{B} \left(
      \vecfnum_{h} - \mat{R} \vec{hv}
    \right),
  \\
  \partial_t \vec{hv}
  =&
  - \frac{1}{2} \left( 
      \mat{D} \vec{h v^2} + \mat{hv} \mat{D} \vec{v} + \mat{v} \mat{D} \vec{hv}
    \right)
  - g \mat{h} \mat{D} \vec{h}
  - g \mat{h} \mat{D} \vec{b}
  \\&
  \widesmall{
  - \mat{M}[^{-1}] \mat{R}[^T] \mat{B} \left(
      \vecfnum_{hv} - \mat{R} \vec{h v^2} - \frac{1}{2} g \mat{R} \vec{h^2}
      + \frac{1}{2} g \mean{h}_- \jump{b}_- \vec{e}_0
      + \frac{1}{2} g \mean{h}_+ \jump{b}_+ \vec{e}_1
    \right)
  }{
  - \mat{M}[^{-1}] \mat{R}[^T] \mat{B} \Big(
      \vecfnum_{hv} - \mat{R} \vec{h v^2} - \frac{1}{2} g \mat{R} \vec{h^2}
      \\&\qquad\qquad\qquad
      + \frac{1}{2} g \mean{h}_- \jump{b}_- \vec{e}_0
      + \frac{1}{2} g \mean{h}_+ \jump{b}_+ \vec{e}_1
    \Big)
  },
\end{aligned}
\end{equation}
where $\vec{e}_k$ is the $k$-th unit vector and for cell $i$
\begin{equation}
\begin{aligned}
  \mean{h}_- &= \frac{h_{i,0} + h_{i-1,p}}{2},
  \quad&
  \mean{h}_+ &= \frac{h_{i,p} + h_{i+1,0}}{2},
  \\
  \jump{b}_- &= b_{i,0} - b_{i-1,p},
  \quad&
  \jump{b}_+ &= b_{i+1,0} - b_{i,p}.
\end{aligned}
\end{equation}
Here, $h_{i,0}, h_{i,p}$ are the values of $h$ at the first and last node $0, p$
in cell $i$, respectively.

Using
\begin{equation}
\label{eq:fnum-gassner}
  \fnum_h
  =
  \mean{h} \mean{v},
  \qquad
  \fnum_{hv}
  =
  \mean{h} \mean{v}^2 + \frac{1}{2} g \mean{h^2},
\end{equation}
as numerical (surface) flux, where $\mean{a} = \frac{a_- + a_+}{2}$, the
resulting scheme
\begin{enumerate}
  \item
  conserves the mass in general and the discharge for a constant bottom topography,
  
  \item
  conserves the total energy which is used as entropy,
  
  \item
  handles the lake-at-rest stationary state correctly,
\end{enumerate}
i.e. it is \emph{conservative}, \emph{stable} and \emph{well-balanced}, as proved
by \citet[Theorem 1]{wintermeyer2016entropy}. The split form discretisation
has been recast into the flux differencing framework of \citet{fisher2013discretely,
fisher2013high} using the "translations" provided by \citet[Lemma 1]{gassner2016split}.
The resulting volume fluxes are
\begin{equation}
\label{eq:fvol-gassner}
  \fvol_h
  =
  \mean{h v},
  \qquad
  \fvol_{hv}
  =
  \mean{h v} \mean{v} + g \mean{h}^2 - \frac{1}{2} g \mean{h^2},
\end{equation}
i.e. the volume terms in \eqref{eq:semidisc-gassner} can be rewritten using the
differentiation matrix $\mat{D}$ as $\sum_{k=0}^p 2 D_{i k} \fvol_{ik}$,
where $\fvol_{ik} = \fvol(u_i,u_k)$. 

The split form in \eqref{eq:semidisc-gassner} corresponds to the entropy
conservative fluxes $\fvol$ \eqref{eq:fvol-gassner}.
$\fnum$ \eqref{eq:fnum-gassner} corresponds to a splitting, too, but
exchanging $\fnum$ for $\fvol$ does not yield a well-balanced method. Similarly,
exchanging $\fvol$ for $\fnum$ as surface flux does not work properly, since
the numerical flux and source discretisation have to be coupled properly in order
to result in a well-balanced discretisation, as described in section \ref{sec:source}.

The remaining parts of this paper are dedicated to the investigation of the
following questions
\begin{enumerate}
  \item 
  Are there other entropy conservative fluxes than $\fvol$, $\fnum$ and
  corresponding split forms? (Section \ref{sec:fnum})
  
  \item
  Are there other discretisations of $g h \partial_x b$ that can be used to get
  a well-balanced scheme, respecting the lake-at-rest stationary state?
  (Section \ref{sec:source})
  
  \item
  Can the split forms be used for a nodal SBP method without boundary nodes,
  e.g. for Gauß nodes? (Section \ref{sec:gauss})
  
  \item
  Are there entropy conservative / stable and positivity preserving numerical
  fluxes that can be used to apply the bound preserving framework of
  \citet{zhang2011maximum}? If so, is the resulting method still entropy stable
  and well-balanced? (Section \ref{sec:fluxes})
\end{enumerate}

\section{Entropy conservative fluxes and split forms for vanishing bottom
         topography \texorpdfstring{$b \equiv 0$}{b=0}}
\label{sec:fnum}

In this section, several numerical fluxes and associated split forms of the
shallow water equations \eqref{eq:SWE} with constant bottom topography $b \equiv 0$
will be considered.

The numerical flux $\fnum(u_-, u_+)$ has to be \emph{consistent}, i.e.
$\fnum(u, u) = f(u)$.
In order to be \emph{entropy conservative}, the condition
\begin{equation}
\label{eq:fnum-entropy-conservative}
  \jump{w} \cdot \fnum(u_-, u_+) = \jump{\psi}
\end{equation}
of \citet{tadmor1987numerical, tadmor2003entropy} has to be fulfilled in a semidiscrete
setting. Similarly, if $\jump{w} \cdot \fnum(u_-, u_+) \leq \jump{\psi}$,
the numerical flux is \emph{entropy stable}, since it contains more dissipation
than an entropy conservative flux. Here, $\jump{a} = a_+ - a_-$.

\subsection{A two-parameter family of entropy conservative numerical fluxes}
\label{sec:fnum-w}

In order to investigate the entropy conservation condition
\eqref{eq:fnum-entropy-conservative}, the jump of the entropy variables $w$
and the flux potential $\psi$ have to be expressed in a common set of variables.
Here, the entropy variables $w$ will be used. Thus, the jump of the flux potential
$\psi = \frac{1}{2 g} w_1^2 w_2 + \frac{1}{2 g} w_1 w_2^3  + \frac{1}{8 g} w_2^5$
\eqref{eq:flux-potential} can be written as
\begin{equation}
\begin{aligned}
  \jump{\psi}
  =&
    \frac{1}{2 g} \jump{w_1^2 w_2}
  + \frac{1}{2 g} \jump{w_1 w_2^3}
  + \frac{1}{8 g} \jump{w_2^5}.
\end{aligned}
\end{equation}
Using a discrete analogue of the product rule
\begin{equation}
\label{eq:product-2}
\begin{aligned}
  \jump{a b}
  &=
  a_+ b_+ - a_- b_-
  =
  \frac{a_+ + a_-}{2} (b_+ - b_-) + (a_+ - a_-) \frac{b_+ + b_-}{2}
  \\&=
  \mean{a} \jump{b} + \jump{a} \mean{b},
\end{aligned}
\end{equation}
the first jump term $\jump{w_1^2 w_2}$ can be written in two different ways.
Weighting both variants with weights $a_1 \in \R$ (first line) and $1-a_1$
(second line), respectively, results in
\begin{align*}
\stepcounter{equation}\tag{\theequation}
  \jump{w_1^2 w_2}
  \stackrel{a_1}{=}&
  \mean{w_1^2} \jump{w_2} + \mean{w_2} \jump{w_1^2}
  =
  \mean{w_1^2} \jump{w_2} + 2 \mean{w_1} \mean{w_2} \jump{w_1}
  \\
  \stackrel{1-a_1}{=}&
  \mean{w_1 w_2} \jump{w_1} + \mean{w_1} \jump{w_1 w_2}
  \\&=
  \mean{w_1 w_2} \jump{w_1} + \mean{w_1}^2 \jump{w_2} + \mean{w_1} \mean{w_2} \jump{w_1}
  \\
  =&
  \left( (1 + a_1) \mean{w_1} \mean{w_2} + (1 - a_1) \mean{w_1 w_2} \right) \jump{w_1}
  \\&
  + \left( a_1 \mean{w_1^2} + (1 - a_1) \mean{w_1}^2 \right) \jump{w_2},
\end{align*}
for $a_1 \in \R$. Similarly, the second jump term $\jump{w_1 w_2^3}$ can be
expressed in four different ways as
\begin{align*}
\stepcounter{equation}\tag{\theequation}
  &
  \jump{w_1 w_2^3}
  \\
  \stackrel{a_2}{=}&
  \mean{w_2^3} \jump{w_1} + \mean{w_1} \jump{w_2^3}
  \widesmall{=}{\\&=}
  \mean{w_2^3} \jump{w_1} + \mean{w_1} \mean{w_2^2} \jump{w_2}
  + 2 \mean{w_1} \mean{w_2}^2 \jump{w_2}
  \\
  \stackrel{a_3}{=}&
  \mean{w_2^2} \jump{w_1 w_2} + \mean{w_1 w_2} \jump{w_2^2}
  \\&=
  \mean{w_2^2} \mean{w_2} \jump{w_1} + \mean{w_1} \mean{w_2^2} \jump{w_2}
  + 2 \mean{w_1 w_2} \mean{w_2} \jump{w_2}
  \\
  \stackrel{a_4}{=}&
  \mean{w_2} \jump{w_1 w_2^2} + \mean{w_1 w_2^2} \jump{w_2}
  \\&=
  \mean{w_2^2} \mean{w_2} \jump{w_1} + 2 \mean{w_1} \mean{w_2}^2 \jump{w_2}
  + \mean{w_1 w_2^2} \jump{w_2}
  \\
  =&
  \mean{w_2} \jump{w_1 w_2^2} + \mean{w_1 w_2^2} \jump{w_2}
  \\&=
  \mean{w_2}^2 \jump{w_1 w_2} + \mean{w_1 w_2} \mean{w_2} \jump{w_2}
  + \mean{w_1 w_2^2} \jump{w_2}
  \\&=
  \mean{w_2}^3 \jump{w_1} + \mean{w_1} \mean{w_2}^2 \jump{w_2}
  + \mean{w_1 w_2} \mean{w_2} \jump{w_2} + \mean{w_1 w_2^2} \jump{w_2}
  \\
  =&
  \left(
    a_2 \mean{w_2^3} + (a_3+a_4) \mean{w_2} \mean{w_2^2}
    + (1-a_2-a_3-a_4) \mean{w_2}^3
  \right) \jump{w_1}
  \\&
  + \Big(
    (a_2+a_3) \mean{w_1} \mean{w_2^2}
    + (1+a_2-a_3+a_4) \mean{w_1} \mean{w_2}^2
    \\&\qquad
    + (1-a_2+a_3-a_4) \mean{w_1 w_2} \mean{w_2}
    + (1-a_2-a_3) \mean{w_1 w_2^2}
  \Big) \jump{w_2},
\end{align*}
where $a_2, a_3, a_4 \in \R$. However, expanding the terms in the last expression
results in
\begin{align*}
\stepcounter{equation}\tag{\theequation}
  &
  a_2 \mean{w_2^3}
  + (a_3+a_4) \mean{w_2} \mean{w_2^2}
  - (a_2+a_3+a_4) \mean{w_2}^3
  \\
  =&
  \frac{3 a_{2} + a_{3} + a_{4}}{8} \left(
    w_{2+}^{3}  - w_{2+}^{2} w_{2-} - w_{2+} w_{2-}^{2} + w_{2-}^{3}
  \right),
\end{align*}
and
\begin{align*}
\stepcounter{equation}\tag{\theequation}
  &
    (a_2+a_3) \mean{w_1} \mean{w_2^2}
  + (a_2-a_3+a_4) \mean{w_1} \mean{w_2}^2
  \\&
  - (a_2-+a_3+a_4) \mean{w_1 w_2} \mean{w_2}
  - (a_2+a_3) \mean{w_1 w_2^2}
  \\
  =&
  \frac{3 a_{2} + a_{3} + a_{4}}{8} \left(
    - w_{1+} w_{2+}^{2} + w_{1+} w_{2-}^{2} + w_{1-} w_{2+}^{2} - w_{1-} w_{2-}^{2} 
  \right).
\end{align*}
Thus, the expression depends only on $3 a_2 + a_3 + a_4$ and can be simplified by
setting $a_3 = a_4 = 0$ to
\begin{equation}
\begin{aligned}
  \jump{w_1 w_2^3}
  =&
  \left(
    a_2 \mean{w_2^3}
    + (1-a_2) \mean{w_2}^3
  \right) \jump{w_1}
  \\&
  + \Big(
      a_2 \mean{w_1} \mean{w_2^2}
    + (1+a_2) \mean{w_1} \mean{w_2}^2
    \\&\qquad
    + (1-a_2) \mean{w_1 w_2} \mean{w_2}
    + (1-a_2) \mean{w_1 w_2^2}
  \Big) \jump{w_2}
\end{aligned}
\end{equation}
Finally, the last jump term can be expressed as
\begin{align*}
\stepcounter{equation}\tag{\theequation}
  &
  \jump{w_2^5}
  \\
  \stackrel{a_5}{=}&
  \mean{w_2^4} \jump{w_2} + \mean{w_2} \jump{w_2^4}
  =
  \mean{w_2^4} \jump{w_2} + 4 \mean{w_2^2} \mean{w_2}^2 \jump{w_2}
  \\
  \stackrel{a_6}{=}&
  \mean{w_2^4} \jump{w_2} + \mean{w_2} \jump{w_2^4}
  \widesmall{=}{\\&=}
  \mean{w_2^4} \jump{w_2} + \mean{w_2}^2 \jump{w_2^3} + \mean{w_2} \mean{w_2^3} \jump{w_2}
  \\&=
  \mean{w_2^4} \jump{w_2} + 2 \mean{w_2}^4 \jump{w_2} + \mean{w_2}^2 \mean{w_2^2} \jump{w_2}
  + \mean{w_2} \mean{w_2^3} \jump{w_2}
  \\
  =&
  \mean{w_2^3} \jump{w_2^2} + \mean{w_2^2} \jump{w_2^3}
  \\&=
  2 \mean{w_2} \mean{w_2^3} \jump{w_2} + 2 \mean{w_2}^2 \mean{w_2^2} \jump{w_2}
  + \mean{w_2^2}^2 \jump{w_2}
  \\
  =&
  \Bigg(
    (a_5 + a_6) \mean{w_2^4}
    + (2 + 2 a_5 - a_6) \mean{w_2}^2 \mean{w_2^2}
    + 2 a_6 \mean{w_2}^4
    \\&\quad
    + (2 - 2 a_5 - a_6) \mean{w_2} \mean{w_2^3}
    + (1 - a_5 - a_6) \mean{w_2^2}^2
  \Bigg) \jump{w_2},
\end{align*}
where $a_5, a_6 \in \R$. However, these parameters $a_5, a_6$ are also redundant,
since the last expression can be simplified as
\begin{equation}
\begin{aligned}
  \jump{w_2^5}
  =&
  \left(
    w_{2+}^{4} + w_{2+}^{3} w_{2-} + w_{2+}^{2} w_{2-}^{2} + w_{2+} w_{2-}^{3} + w_{2-}^{4}
  \right) \jump{w_2}
  \\
  =&
  \left( \mean{w_2^4} + 4 \mean{w_2}^2 \mean{w_2^2} \right) \jump{w_2}.
\end{aligned}
\end{equation}
Inserting these forms in the condition $\jump{w} \cdot \fnum(u_-, u_+) = \jump{\psi}$
\eqref{eq:fnum-entropy-conservative} for an entropy conservative flux,
\begin{align*}
\stepcounter{equation}\tag{\theequation}
\label{eq:gnum}
  f^{a_1,a_2}_{h}(w_-, w_+)
  =&
    \frac{1 + a_1}{2 g} \mean{w_1} \mean{w_2} 
  + \frac{1 - a_1}{2 g} \mean{w_1 w_2}
  + \frac{a_2}{2 g} \mean{w_2^3}
  \widesmall{}{\\&}
  + \frac{1-a_2}{2 g} \mean{w_2}^3,
  \\
  f^{a_1,a_2}_{hv}(w_-, w_+)
  =&
    \frac{a_1}{2 g} \mean{w_1^2}
  + \frac{1 - a_1}{2 g} \mean{w_1}^2
  + \frac{a_2}{2 g} \mean{w_1} \mean{w_2^2}
  \\&
  + \frac{1+a_2}{2 g} \mean{w_1} \mean{w_2}^2
  + \frac{1-a_2}{2 g} \mean{w_1 w_2} \mean{w_2}
  \widesmall{}{\\&}
  + \frac{1-a_2}{2 g} \mean{w_1 w_2^2}
  \widesmall{\\&}{}
  + \frac{1}{8 g} \mean{w_2^4}
  + \frac{1}{2 g} \mean{w_2}^2 \mean{w_2^2},
\end{align*}
turns out to be a two-parameter family of entropy conservative numerical fluxes
expressed in terms of the entropy variables
$w = \begin{pmatrix} g h - \frac{1}{2} v^2, v \end{pmatrix}^T$
\eqref{eq:entropy-variables-b-constant} for vanishing bottom topography
$b \equiv 0$.
Expanding the terms in entropy and primitive variables, respectively, results
after direct but tedious calculations in
\begin{align*}
\stepcounter{equation}\tag{\theequation}
\label{eq:gnum-w}
  f^{a_1,a_2}_{h}
  =&
    \frac{ 3 - a_{1} }{8 g} \left( w_{1+} w_{2+} + w_{1-} w_{2-} \right)
  + \frac{ 1 + a_{1} }{8 g} \left( w_{1+} w_{2-} + w_{1-} w_{2+} \right)
  \\&
  + \frac{ 1 + 3 a_{2} }{16 g} \left( w_{2+}^{3} + w_{2-}^{3} \right)
  + \frac{ 3 - 3 a_{2} }{16 g} \left( w_{2+}^{2} w_{2-} + w_{2+} w_{2-}^{2} \right),
  \\
  f^{a_1,a_2}_{hv}
  =&
    \frac{ 1 + a_{1} }{8 g} \left( w_{1+}^{2} + w_{1-}^{2} \right)
  + \frac{ 1 - a_{1} }{4 g} w_{1+} w_{1-}
  \widesmall{}{\\&}
  + \frac{ 7 - 3 a_{2} }{16 g} \left( w_{1+} w_{2+}^{2} + w_{1-} w_{2-}^{2} \right)
  \widesmall{\\&}{}
  + \frac{ 1 + 3 a_{2} }{16 g} \left( w_{1+} w_{2-}^{2} + w_{1-} w_{2+}^{2} \right)
  \widesmall{}{\\&}
  + \frac{w_{1+} w_{2+} w_{2-}}{4 g}
  + \frac{w_{1-} w_{2+} w_{2-}}{4 g}
  \\&
  + \frac{1}{8 g} \left(
      w_{2+}^{4}
    + w_{2+}^{3} w_{2-}
    + w_{2+}^{2} w_{2-}^{2}
    + w_{2+} w_{2-}^{3}
    + w_{2-}^{4}
  \right),
\end{align*}
and
\begin{align*}
\stepcounter{equation}\tag{\theequation}
\label{eq:gnum-p}
  f^{a_1,a_2}_{h}
  =&
    \frac{ 3 - a_{1} }{8} \left( h_{+} v_{+} + h_{-} v_{-} \right)
  + \frac{ 1 + a_{1} }{8} \left( h_{+} v_{-} + h_{-} v_{+} \right)
  \\&
  + \frac{ a_{1} + 3 a_{2} - 2 }{16 g} \left(
    v_{+}^{3} - v_{+}^{2} v_{-} - v_{+} v_{-}^{2} + v_{-}^{3}
  \right),
  \\
  f^{a_1,a_2}_{hv}
  =&
    \frac{1 + a_{1} }{8} g \left( h_{+}^{2} + h_{-}^{2} \right)
  + \frac{1 - a_{1} }{4} g h_{+} h_{-}
  - \frac{2 a_{1} + 3 a_{2} - 5 }{16} \left( h_{+} v_{+}^{2} + h_{-} v_{-}^{2} \right)
  \\&
  + \frac{2 a_{1} + 3 a_{2} - 1 }{16} \left( h_{+} v_{-}^{2} + h_{-} v_{+}^{2} \right)
  + \frac{h_{+} v_{+} v_{-}}{4}
  + \frac{h_{-} v_{+} v_{-}}{4}
  \\&
  + \frac{ a_{1} + 3 a_{2} - 2 }{32 g} \left(
    v_{+}^{4} - 2 v_{+}^{2} v_{-}^{2} + v_{-}^{4}
  \right).
\end{align*}
This proves
\begin{lemma}
\label{lem:gnum}
  The two-parameter family \eqref{eq:gnum} of numerical fluxes $f^{a_1,a_2}$,
  expressed also as \eqref{eq:gnum-w} and \eqref{eq:gnum-p}, is a family of
  consistent and entropy conservative numerical fluxes for the shallow water
  equations \eqref{eq:SWE} with vanishing bottom topography $b \equiv 0$.
\end{lemma}

To the author's knowledge, this family of entropy conservative fluxes has not
been considered en bloc in the literature before. Instead, only the fluxes corresponding
to three distinct choices of the parameters $a_1, a_2$ have been proposed,
as described in the Remarks \ref{rem:tadmor} and \ref{rem:gassner} in the
following section.

\subsection{Relations with known fluxes and methods}

The crucial ingredient to obtain the entropy conservative fluxes in Lemma
\ref{lem:gnum} has been the expression of both the entropy
variables $w$ and the flux potential $\psi$ as polynomials in the same variables.
Therefore, it may be conjectured, that if such a condition is complied with,
there are entropy conservative fluxes expressed in terms of these variables,
corresponding to a split form as described inter alia by \citet{gassner2016split}
and in section \ref{sec:gauss}.

\begin{remark}
\label{rem:tadmor}
  The general entropy conservative flux of
  \citet[Equation (4.6a)]{tadmor1987numerical}, obtained by integration in phase
  space, can be recovered by the family \eqref{eq:gnum} of entropy conservative
  fluxes. Indeed, for $f\left( u(w) \right)$ as in \eqref{eq:g},
  \begin{align*}
  \stepcounter{equation}\tag{\theequation}
    &
    \int_0^1 f_{h} \circ u \left( (1-s) w_- + s w_+ \right) \dif s
    \\
    =&
%     \frac{1}{g} \int_0^1 \Bigg(
%       \left( (1-s) w_{1-} + s w_{1+} \right) \left( (1-s) w_{2-} + s w_{2+} \right)
%       \widesmall{}{\\&\qquad\quad}
%       + \frac{1}{2} \left( (1-s) w_{2-} + s w_{2+} \right)^3
%     \Bigg) \dif s
%     \\
%     =&
    \frac{1}{g} \Bigg(
      \frac{w_{1+} w_{2+}}{3}
      + \frac{w_{1+} w_{2-}}{6}
      + \frac{w_{1-} w_{2+}}{6}
      + \frac{w_{1-} w_{2-}}{3}
      \widesmall{}{\\&\qquad}
      + \frac{w_{2+}^{3}}{8}
      + \frac{w_{2+}^{2} w_{2-}}{8}
      + \frac{w_{2+} w_{2-}^{2}}{8}
      + \frac{w_{2-}^{3}}{8}
    \Bigg),
  \end{align*}
  and
  \begin{align*}
  \stepcounter{equation}\tag{\theequation}
    &
    \int_0^1 f_{hv} \circ u\left( (1-s) w_- + s w_+ \right) \dif s
    \\
    =&
%     \frac{1}{g} \int_0^1 \Bigg(
%       \frac{1}{2} \left( (1-s) w_{1-} + s w_{1+} \right)^2
%       \widesmall{}{\\&\qquad\quad}
%       + \frac{3}{2} \left( (1-s) w_{1-} + s w_{1+} \right) \left( (1-s) w_{2-} + s w_{2+} \right)^2
%       \\&\qquad\quad
%       + \frac{5}{8} \left( (1-s) w_{2-} + s w_{2+} \right)^4
%     \Bigg) \dif s
%     \\
%     =&
    \frac{1}{24 g} \Big(
      4 w_{1+}^{2} + 4 w_{1+} w_{1-}
      + 9 w_{1+} w_{2+}^{2}
      + 6 w_{1+} w_{2+} w_{2-}
      + 3 w_{1+} w_{2-}^{2}
      \widesmall{}{\\&\qquad}
      + 4 w_{1-}^{2}
      + 3 w_{1-} w_{2+}^{2}
      \widesmall{\\&\qquad}{}
      + 6 w_{1-} w_{2+} w_{2-}
      + 9 w_{1-} w_{2-}^{2}
      \widesmall{}{\\&\qquad}
      + 3 w_{2+}^{4}
      + 3 w_{2+}^{3} w_{2-}
      + 3 w_{2+}^{2} w_{2-}^{2}
      + 3 w_{2+} w_{2-}^{3}
      + 3 w_{2-}^{4}
    \Big).
  \end{align*}
  Comparing this with the numerical fluxes $f^{a_1,a_2}$ \eqref{eq:gnum-w},
  it can be seen that Tadmor's flux as above is recovered by setting
  $a_1 = a_2 = \frac{1}{3}$.
\end{remark}

\begin{remark}
  Using entropy variables, there does not seem to be a special choice of the
  parameters $a_1, a_2$ for the numerical flux \eqref{eq:gnum-w}. However,
  the expression in primitive variables $h, v$ reveals that the choice
  $a_2 = \frac{2 - a_1}{3}$ is special, since higher order terms in the velocity
  that are consistent with zero are removed. In this case, the one-parameter family
  of entropy conservative fluxes can be written as
  \begin{equation}
  \label{eq:fnum}
  \begin{aligned}
    f^{a_1}_{h}
    =&
    \frac{1-a_1}{2} \mean{h v} + \frac{1+a_1}{2} \mean{h} \mean{v},
    \\
    f^{a_1}_{hv}
    =&
    \frac{1-a_1}{2} \mean{h v} \mean{v} + \frac{1+a_1}{2} \mean{h} \mean{v}^2
    + g \frac{a_1}{2} \mean{h^2} + \frac{1-a_1}{2} g \mean{h}^2,
  \end{aligned}
  \end{equation}
  or
  \begin{align*}
  \stepcounter{equation}\tag{\theequation}
  \label{eq:fnum-p}
    f^{a_1}_{h}
    =&
      \frac{ 3 - a_1 }{8} \left( h_{+} v_{+} + h_{-} v_{-} \right)
    + \frac{ 1 + a_1 }{8} \left( h_{+} v_{-} + h_{-} v_{+} \right),
    \\
    f^{a_1}_{hv}
    =&
      \frac{1 + a_1 }{8} g \left( h_{+}^{2} + h_{-}^{2} \right)
    + \frac{1 - a_1 }{4} g h_{+} h_{-}
    + \frac{ 3 - a_1 }{16} \left( h_{+} v_{+}^{2} + h_{-} v_{-}^{2} \right)
    \\&
    + \frac{ 1 + a_1 }{16} \left( h_{+} v_{-}^{2} + h_{-} v_{+}^{2} \right)
    + \frac{h_{+} v_{+} v_{-}}{4}
    + \frac{h_{-} v_{+} v_{-}}{4}.
  \end{align*}
\end{remark}

\begin{remark}
\label{rem:gassner}
  The numerical surface and volume fluxes $\fvol$ \eqref{eq:fvol-gassner},
  $\fnum$ \eqref{eq:fnum-gassner} of \citet{gassner2016well} and
  \citet{wintermeyer2016entropy} are members of this family with parameters
  $a_1 = -1$ and $a_1 = 1$, respectively. Therefore, this one-parameter
  family of entropy conservative fluxes \eqref{eq:fnum}
  can also be seen as linear combinations of the two fluxes
  \eqref{eq:fnum-gassner} and \eqref{eq:fvol-gassner}
  with coefficients summing up to one.
\end{remark}

\begin{remark}
  The derivations of section \ref{sec:fnum-w} can also be conducted using
  primitive variables $h, v$ instead of the entropy variables $w$. In this case,
  only the one-parameter family of entropy conservative numerical fluxes
  \eqref{eq:fnum}, \eqref{eq:fnum-p} can be obtained.
\end{remark}

\begin{remark}
  As proved by \citet[Theorem 3.2]{fisher2013high}, a two-point entropy
  conservative flux as $f^{a_1,a_2}$ \eqref{eq:gnum} can be used to construct
  a high-order spatial discretisation for diagonal-norm SBP operators with nodal
  basis including boundary nodes.
  \citet[Lemma 1]{gassner2016split} provided some examples for analogous split 
  forms and numerical fluxes given by simple products of mean values.
  Analogously, using a diagonal-norm SBP derivative operator $\mat{D}$ (i.e. an
  SBP derivative operator $\mat{D}$, where the corresponding norm / mass matrix
  $\mat{M}$ is diagonal) with nodal basis including boundary nodes, the split
  form corresponding to the flux $f^{a_1,a_2}$ expressed via primitive
  variables \eqref{eq:gnum-p} is for the $h$ component given by
  \begin{align*}
  \stepcounter{equation}\tag{\theequation}
  \label{eq:gnum-volume-terms-h}
    &
    \left[ \VOL^{a_1,a_2}_{h} \right]_i
    =
    \sum_{k=0}^p 2 D_{i,k} f^{a_1,a_2}_{h}(u_i,u_k)
    \\
    =&
    \sum_{k=0}^p 2 D_{i,k} \Bigg(
        \frac{ 3 - a_{1} }{8} \left( h_i v_i + h_k v_k \right)
      + \frac{ 1 + a_{1} }{8} \left( h_i v_k + h_k v_i \right)
      \\&\qquad\qquad\quad
      + \frac{ a_{1} + 3 a_{2} - 2 }{16 g} \left(
        v_i^3 - v_i^2 v_k - v_i v_k^2 + v_k^3
      \right)
    \Bigg)
    \\
    =&
    \sum_{k=0}^p D_{i,k} \Bigg(
        \frac{ 3 - a_{1} }{4}  h_k v_k
      + \frac{ 1 + a_{1} }{4} \left( h_i v_k + h_k v_i \right)
      \widesmall{}{\\&\qquad\qquad}
      + \frac{ a_{1} + 3 a_{2} - 2 }{8 g} \left( v_k^3 - v_i v_k^2 - v_i^2 v_k \right)
    \Bigg)
    \\
    =&
    \Bigg[
        \frac{ 3 - a_{1} }{4} \mat{D} \vec{h v}
      + \frac{ 1 + a_{1} }{4} \left( \mat{h} \mat{D} \vec{v} + \mat{v} \mat{D} \vec{h} \right)
      \widesmall{}{\\&\quad}
      + \frac{ a_{1} + 3 a_{2} - 2 }{8 g} \left(
        \mat{D} \vec{v^3} - \mat{v} \mat{D} \vec{v^2} - \mat{v^2} \mat{D} \vec{v} \right)
    \Bigg]_i.
  \end{align*}
  Here, some summands have been dropped, because the derivative is exact for
  constants, i.e. $\mat{D} \vec{1} = 0$, resulting in $\sum_{k=0}^p D_{i,k} = 0$.
  The first two terms form a consistent discretisation of $\partial_x (h v)
  = (\partial_x h) v + h (\partial_x v)$ for smooth solutions. The third term
  is consistently zero, since $\partial_x v^3 = (\partial_x v^2) v
  + v^2 (\partial_x v)$ by the product rule for smooth solutions.
  Similarly, the $hv$ component can be computed via
  \begin{align*}
  \stepcounter{equation}\tag{\theequation}
  \label{eq:gnum-volume-terms-hv}
    &
    \left[ \VOL^{a_1,a_2}_{hv} \right]_i
    =
    \sum_{k=0}^p 2 D_{i,k} f^{a_1,a_2}_{hv}(u_i,u_k)
    \\
    =&
    \sum_{k=0}^p 2 D_{i,k} \Bigg(
        \frac{1 + a_{1} }{8} g \left( h_i^2 + h_k^2 \right)
      + \frac{1 - a_{1} }{4} g h_i h_k
      \widesmall{}{\\&\qquad\qquad\quad}
      - \frac{2 a_{1} + 3 a_{2} - 5 }{16} \left( h_i v_i^2 + h_k v_k^2 \right)
      \widesmall{\\&\qquad\qquad\quad}{}
      + \frac{2 a_{1} + 3 a_{2} - 1 }{16} \left( h_i v_k^2 + h_k v_i^2 \right)
      \widesmall{}{\\&\qquad\qquad\quad}
      + \frac{h_i v_i v_k}{4}
      + \frac{h_k v_i v_k}{4}
      \widesmall{\\&\qquad\qquad\quad}{}
      + \frac{ a_{1} + 3 a_{2} - 2 }{32 g} \left( v_i^4 - 2 v_i^2 v_k^2 + v_k^4 \right)
    \Bigg)
    \\
    =&
    \sum_{k=0}^p D_{i,k} \Bigg(
        \frac{1 + a_{1} }{4} g h_k^2
      + \frac{1 - a_{1} }{2} g h_i h_k
      - \frac{2 a_{1} + 3 a_{2} - 5 }{8} h_k v_k^2
      \widesmall{}{\\&\qquad\qquad\quad}
      + \frac{2 a_{1} + 3 a_{2} - 1 }{8} \left( h_i v_k^2 + h_k v_i^2 \right)
      \widesmall{\\&\qquad\qquad\quad}{}
      + \frac{h_i v_i v_k}{2}
      + \frac{h_k v_i v_k}{2}
      \widesmall{}{\\&\qquad\qquad\quad}
      + \frac{ a_{1} + 3 a_{2} - 2 }{16 g} \left( v_k^4 - 2 v_i^2 v_k^2 \right)
    \Bigg)
    \\
    =&
    \Bigg[
        \frac{1 + a_{1} }{4} g \mat{D} \vec{h^2}
      + \frac{1 - a_{1} }{2} g \mat{h} \mat{D} \vec{h}
      - \frac{2 a_{1} + 3 a_{2} - 5 }{8} \mat{D} \vec{h v^2}
      \widesmall{}{\\&\quad}
      + \frac{2 a_{1} + 3 a_{2} - 1 }{8} \left( \mat{h} \mat{D} \vec{v^2}
                                                + \mat{v^2} \mat{D} \vec{h} \right)
      \widesmall{\\&\quad}{}
      + \frac{1}{2} \left( \mat{h v} \mat{D} \vec{v} + \mat{v} \mat{D} \vec{h v} \right)
      \widesmall{}{\\&\quad}
      + \frac{ a_{1} + 3 a_{2} - 2 }{16 g} \left( \mat{D} \vec{v^4}
        - 2 \mat{v^2} \mat{D} \vec{v^2} \right)
    \Bigg]_i,
  \end{align*}
  where again $\sum_{k=0}^p D_{i,k} = 0$ has been used. The first two terms form
  a consistent discretisation of $\frac{1}{2} g \partial_x h^2$, the three following
  terms are consistent with $\partial_x h v^2$ and the last two terms are
  consistently zero.
  
  The entropy conservation follows from the general result of
  \citet[Theorem 3.2]{fisher2013high} and will also be investigated in more detail
  for Gauß nodes and other SBP bases in section \ref{sec:gauss}.
\end{remark}

\begin{remark}
  It is also possible to start with a general ansatz of the split form and use
  conditions for consistency, conservation, and entropy stability (similar to
  section \ref{sec:gauss}), in order to determine the coefficients.
  This yields the same two-parameter family of fluxes and corresponding split
  forms, but is much more tedious.
\end{remark}

\section{Adding well-balanced source discretisations}
\label{sec:source}

In this section, the discretisation of the source term $g h \partial_x b$ in
the shallow water equations \eqref{eq:SWE} will be investigated. It should be
\emph{consistent}, \emph{stable}, and \emph{well-balanced}, if combined with
the remaining semidiscretisation derived in the previous section.

\subsection{Connections between finite volume and SBP SAT schemes}

A general semidiscretisation of a conservation law $\partial_t u + \partial_x f(u) = 0$
with a polynomial SBP method using the notation of \citet{ranocha2016summation}
can be written as
\begin{equation}
  \partial_t \vec{u}
  =
  - \VOL
  + \SURF
  - \mat{M}[^{-1}] \mat{R}[^T] \mat{B} \vecfnum,
\end{equation}
where $\VOL$ contains the volume terms consistent with $\partial_x f(u)$,
possibly using some split form, $\fnum$ is the numerical (surface) flux, and
$\SURF$ contains additional surface terms, consistent with the difference of the
flux values $f(u)$ at the boundaries, i.e. $\mat{M}[^{-1}] \mat{R}[^T] \mat{B}
\mat{R} \vec{f}$ in the simplest case, but additional terms may also appear,
especially if a nodal basis without boundary nodes is used, see also section
\ref{sec:gauss}.

If the polynomial degree $p$ is set to zero, the volume terms vanish, since the
derivative is exact for constants, i.e. $\mat{D} \vec{1} = 0$. Additionally,
since the extra surface terms $\SURF$ are a consistent evaluation of the difference
of boundary values, they vanish, too, because this difference is zero for constants.
Therefore, this method reduces to a simple finite volume method. If the cell $i$
is of size $\Delta x$ and the flux $\fnum$ between the cells $i$ and $k$ is
denoted as $\fnum_{i,k}$, the FV method can be written as
\begin{equation}
\label{eq:FV}
  \partial_t u_i
  =
  - \frac{1}{\Delta x} \left( \fnum_{i,i+1} - \fnum_{i,i-1} \right),
\end{equation}
and is determined solely by the numerical flux $\fnum$ used at the boundaries.

On the other hand, using the theory of \citet{fisher2013high}, a finite volume
method with entropy conservative flux $\fnum$ can be used to construct the volume
terms $\VOL$, if a nodal SBP basis including boundary nodes is used. In this case,
since the evaluation at the boundary is exact and commutes with nonlinear
operations, the surface terms are simply $\mat{M}[^{-1}] \mat{R}[^T] \mat{B}
\mat{R} \vec{f}$.

This strong correspondence between SBP schemes and FV methods will be used in
the following sections to extend results from one area to the other and vice
versa.

\subsection{Extended numerical fluxes and entropy conservation}

In order to incorporate source terms in a finite volume method, the additional
contributions can be incorporated into the numerical fluxes, resulting in
extended numerical fluxes as described inter alia in the monograph by
\citet[Section 4]{bouchut2004nonlinear} and references cited therein. Here,
the extended flux for the discharge $hv$ will have the form
\begin{equation}
\label{eq:extended-flux}
  f^\mathrm{num,ext}_{i,k} = \fnum_{i,k} + S_{i,k},
\end{equation}
where $\fnum_{i,k}$ is a usual (entropy conservative and symmetric) numerical flux
of the problem without source terms and $S_{i,k}$ describes the source terms and
is in general not symmetric.

\begin{remark}
  Setting the polynomial degree $p$ to zero, the surface flux $\fnum_{hv}$
  and the surface source terms of the SBP SAT semidiscretisation of
  \citet{wintermeyer2016entropy} can be written as the FV method of
  \citet{fjordholm2011well}, which can be written using the extended numerical flux
  \begin{equation}
  \label{eq:extended-flux-fnum}
    \fnum_{{hv}_{i,k}} + \frac{1}{2} g \mean{h}_{i,k} \jump{b}_{i,k},
    \qquad
    \fnum_{{hv}_{i,k}}
    =
    \mean{h}_{i,k} \mean{v}_{i,k}^2 + \frac{1}{2} g \mean{h^2}_{i,k},
  \end{equation}
  where $\mean{a}_{i,k} = \frac{a_i + a_k}{2}$, $\jump{a}_{i,k} = a_k - a_i$.

\end{remark}

Rewriting the FV evolution equation \eqref{eq:FV} by adding $f_i - f_i = 0$ (motivated
by the form of SBP SAT methods) and using extended numerical fluxes yields
\begin{equation}
  \partial_t u_i
  =
  -\frac{1}{\Delta x} \left(
      \left( f^\mathrm{num,ext}_{i,i+1} - f_i \right)
    - \left( f^\mathrm{num,ext}_{i,i-1} - f_i \right)
  \right).
\end{equation}
Therefore, the rate of change of the entropy $U$ can be calculated as
\begin{equation}
\begin{aligned}
  \partial_t U_i
  =&
  w_i \cdot \partial_t u_i
  =
  -\frac{1}{\Delta x} \left(
      w_i \cdot \left( f^\mathrm{num,ext}_{i,i+1} - f_i \right)
    - w_i \cdot \left( f^\mathrm{num,ext}_{i,i-1} - f_i \right)
  \right)
  \\
  =&
  -\frac{1}{\Delta x} \left(
      \left[ w_i \cdot \left( f^\mathrm{num,ext}_{i,i+1} - f_i \right) + F_i \right]
    - \left[ w_i \cdot \left( f^\mathrm{num,ext}_{i,i-1} - f_i \right) + F_i \right]
  \right).
\end{aligned}
\end{equation}
Thus, adding the contributions of the right hand side of cell $i$ and the left
hand side of cell $i+1$ yields after multiplication with $\Delta x$
\begin{align*}
\stepcounter{equation}\tag{\theequation}
\label{eq:FV-entropy-contribution-one-boundary-1}
  &
    \left[ w_{i+1} \cdot \left( f^\mathrm{num,ext}_{i+1,i} - f_{i+1} \right) + F_{i+1} \right]
  - \left[ w_{i} \cdot \left( f^\mathrm{num,ext}_{i,i+1} - f_{i} \right) + F_{i} \right]
  \\
  =&
  (w_{i+1} - w_{i}) \cdot \fnum_{i,i+1}
  - \Big( \underbrace{ \left[ w_{i+1} \cdot f_{i+1} + F_{i+1} \right] }_{=\psi_{i+1}}
            - \underbrace{ \left[ w_{i} \cdot f_{i} + F_{i} \right] }_{=\psi_{i}}
    \Big)
  \widesmall{}{\\&}
  + w_{i+1} \cdot S_{i+1,i}
  - w_{i} \cdot S_{i,i+1},
\end{align*}
where the extended flux $f^\mathrm{num,ext}$ \eqref{eq:extended-flux} has been
inserted and the symmetry of the numerical flux $\fnum$ has been used.

Assume now that the numerical flux $\fnum_{i,i+1}$ is chosen as an entropy
conservative one, fulfilling $\jump{w}_{i,i+1} \cdot \fnum_{i,i+1} = \jump{\psi}_{i,i+1}$
\eqref{eq:fnum-entropy-conservative}, for vanishing bottom topography $b \equiv 0$.
Here, the entropy variables are  $w = \begin{pmatrix} g h - \frac{1}{2} v^2,
v \end{pmatrix}^T$, since $b \equiv 0$. In the general case, the entropy variables
are  $w = \begin{pmatrix} g (h + b) - \frac{1}{2} v^2, v \end{pmatrix}^T$,
resulting in $\jump{w}_{i,i+1} \cdot \fnum_{i,i+1} = \jump{\psi}_{i,i+1}
+ g \fnum_{h_{i,i+1}} \jump{b}_{i,i+1}$.
Thus, the contribution of one boundary to the rate of change of the entropy
\eqref{eq:FV-entropy-contribution-one-boundary-1} is
\begin{equation}
\label{eq:FV-entropy-contribution-one-boundary-2}
  g \fnum_{h_{i,i+1}} \jump{b}_{i,i+1}
  + w_{i+1} \cdot S_{i+1,i}
  - w_{i} \cdot S_{i,i+1}
  \stackrel{!}{=} 0.
\end{equation}
This proves
\begin{lemma}
\label{lem:entropy-conservation-condition-FV}
  If the source discretisation $S_{i,k}$ in the extended numerical flux
  \eqref{eq:extended-flux} is chosen such that the expression
  \eqref{eq:FV-entropy-contribution-one-boundary-2} is zero for an entropy
  conservative numerical flux $\fnum$ \eqref{eq:fnum-entropy-conservative}
  for the shallow water equations with vanishing bottom topography $b \equiv 0$,
  then the resulting scheme is entropy conservative for general bottom topography.
\end{lemma}

To the author's knowledge, this general result has not been formulated before.
Instead, special source discretisations have been chosen and similar calculations
adapted to the specific discretisation have been performed.

\begin{remark}
  The source discretisation of \citet[Lemma 2.1]{fjordholm2011well} and
  \citet{wintermeyer2016entropy} results in the extended numerical flux
  \eqref{eq:extended-flux-fnum} with source terms $S_{i,k} = \frac{1}{2}
  g \mean{h}_{i,k} \jump{b}_{i,k}$. Thus, inserting this and their numerical flux
  $\fnum_{h_{i,k}} = \mean{h}_{i,k} \mean{v}_{i,k}$ into
  \eqref{eq:FV-entropy-contribution-one-boundary-2} results in
  \begin{align*}
  \stepcounter{equation}\tag{\theequation}
    &
      g \mean{h}_{i,i+1} \mean{v}_{i,i+1} \jump{b}_{i,i+1}
    + v_{i+1} \, \frac{g}{2} \mean{h}_{i,i+1} \jump{b}_{i+1,i}
    - v_{i} \, \frac{g}{2} \mean{h}_{i+1,i} \jump{b}_{i,i+1}
    \\
    =&
    g \mean{h}_{i,i+1} \mean{v}_{i,i+1} \jump{b}_{i,i+1}
    - g \mean{h}_{i,i+1} \frac{v_i + v_{i+1}}{2} \jump{b}_{i,i+1}
    = 0,
  \end{align*}
  fulfilling the condition of Lemma \ref{lem:entropy-conservation-condition-FV}.
\end{remark}

\begin{remark}
  The discretisation of the volume terms for the discharge of \citet{gassner2016well}
  and \citet{wintermeyer2016entropy} can be formulated  using the flux difference
  form and an extended numerical flux
  \begin{equation}
  \label{eq:extended-flux-fvol}
    \fvol_{hv_{i,k}} + \frac{1}{2} g \, h_i \jump{b}_{i,k},
    \qquad
    \fvol_{hv_{i,k}}
    =
    \mean{h v}_{i,k} \mean{v}_{i,k} + g \mean{h}_{i,k}^2 - \frac{1}{2} g \mean{h^2}_{i,k}.
  \end{equation}
  Indeed, since $\mat{D} \vec{1} = 0$, the discretisation of 
  $\frac{1}{2} g \partial_x h^2 + g h \partial_x b$ is given by
  \begin{align*}
  \stepcounter{equation}\tag{\theequation}
    &
    \sum_{k=0}^p 2 D_{i,k} \left(
      g \mean{h}_{i,k}^2 - \frac{1}{2} g \mean{h^2}_{i,k}
      + \frac{1}{2} g \, h_i \jump{b}_{i,k}
    \right)
    \\
    =&
    g \sum_{k=0}^p 2 D_{i,k} \left(
      \left( \frac{h_i + h_k}{2} \right)^2 - \frac{1}{2} \frac{h_i^2 + h_k^2}{2}
      + \frac{1}{2} h_i (b_k - b_i)
    \right)
    \\
    =&
    g \sum_{k=0}^p D_{i,k} \left(
      h_i h_k + h_i (b_k - b_i)
    \right)
    =
    g h_i \underbrace{\sum_{k=0}^p D_{i,k} (h_k + b_k)}_{= [\mat{D} (\vec{h} + \vec{b})]_i}
    - g h_i b_i \underbrace{\sum_{k=0}^p D_{i,k}}_{= [\mat{D} \vec{1}]_i=0}.
  \end{align*}
\end{remark}

In the same way the entropy conservative numerical fluxes \eqref{eq:fnum-gassner}
and \eqref{eq:fvol-gassner} can be combined to get the one-parameter family of
entropy conservative fluxes $f^{a_1}$ \eqref{eq:fnum} for the shallow water
equations with vanishing bottom topography $b \equiv 0$, the extended numerical
fluxes \eqref{eq:extended-flux-fnum} and \eqref{eq:extended-flux-fvol}
can be combined to get entropy conservative extended numerical fluxes for
the shallow water equations with general bottom topography.

For the two-parameter family \eqref{eq:gnum}, the source discretisation $S_{i,k}$
in the extended numerical flux \eqref{eq:extended-flux} has to be adapted to the
additional terms with $a_2$ in order to fulfil the condition
\eqref{eq:FV-entropy-contribution-one-boundary-2} of Lemma
\ref{lem:entropy-conservation-condition-FV}. Since the two-parameter flux
\eqref{eq:gnum-p} for $h$ contains an additional term
$\frac{a_1 + 3 a_2 - 2}{16 g} \left(v_+^3 - v_+^2 v_- - v_+ v_-^2 + v_-^3 \right)$
compared to the one-parameter flux \eqref{eq:fnum-p}, the new source term
$S_{i,k}$ can be written as the sum of the source term
$\frac{1}{4} g \left( \frac{3-a_1}{2} h_i + \frac{1+a_1}{2} h_k \right)(b_k - b_i)$
for the one-parameter flux \eqref{eq:fnum-p} and an additional source term
$\tilde S_{i,k}$, obeying
\begin{align*}
\stepcounter{equation}\tag{\theequation}
  &
  g \fnum_{h_{i,k}} \jump{b}_{i,k}
  + w_{k} \cdot S_{k,i}
  - w_{i} \cdot S_{i,k}
  \\
  =&
  \frac{a_1 + 3 a_2 - 2}{16} \left(v_i^3 - v_i^2 v_k - v_i v_k^2 + v_k^3 \right) (b_k - b_i)
  + v_k \tilde S_{k,i} - v_i \tilde S_{i,k}
  \stackrel{!}{=}
  0.
\end{align*}
This can be rewritten using
$v_i^3 - v_i^2 v_k - v_i v_k^2 + v_k^3 = (v_i + v_k) (v_i - v_k)^2$.
Thus, choosing $\tilde S_{i,k} = \frac{a_1 + 3 a_2 - 2}{16} (v_k - v_i)^2 (b_k - b_i)$
results in the desired equality
\begin{align*}
\stepcounter{equation}\tag{\theequation}
  &
  \frac{a_1 + 3 a_2 - 2}{16} \left(v_i^3 - v_i^2 v_k - v_i v_k^2 + v_k^3 \right) (b_k - b_i)
  + v_k \tilde S_{k,i} - v_i \tilde S_{i,k}
  \\
  =&
  \frac{a_1 + 3 a_2 - 2}{16} \Big(
    (v_i + v_k) (v_i - v_k)^2  (b_k - b_i)
    + v_k (v_i - v_k)^2 (b_i - b_k)
    \widesmall{}{\\&\qquad\qquad\qquad\quad}
    - v_i (v_k - v_i)^2 (b_k - b_i)
  \Big)
  \\
  =&
  0.
\end{align*}

\subsection{Well-balancedness (preserving the lake-at-rest steady state)}

The extended numerical flux \eqref{eq:extended-flux-fnum} of \citet{fjordholm2011well}
and \citet{wintermeyer2016entropy} preserves the lake-at-rest steady state
$h+b \equiv \mathrm{const}, hv = 0$ in a finite volume method \eqref{eq:FV}.

Similarly, the extended numerical flux \eqref{eq:extended-flux-fvol} is well-balanced,
since for $h+b \equiv \mathrm{const}$ and $hv = 0$
\begin{align*}
\stepcounter{equation}\tag{\theequation}
  &
  \frac{1}{g} \left( \fvol_{hv_{i,,i+1}} + \frac{1}{2} g \, h_i \jump{b}_{i,i+1} \right)
  - \frac{1}{g} \left( \fvol_{hv_{i,,i-1}} + \frac{1}{2} g \, h_i \jump{b}_{i,i-1} \right)
  \\
  =&
  \mean{h}_{i,i+1}^2 - \frac{1}{2} \mean{h^2}_{i,i+1}
  + \frac{1}{2} h_i \jump{b}_{i,i+1}
  - \mean{h}_{i,i-1}^2 + \frac{1}{2} \mean{h^2}_{i,i-1}
  - \frac{1}{2} h_i \jump{b}_{i,i-1}
  \\
  =&
  \frac{1}{2} h_{i} h_{i+1} + \frac{1}{2} h_i b_{i+1}
  - \frac{1}{2} h_{i} h_{i-1} - \frac{1}{2} h_i b_{i-1}
  = 0.
\end{align*}

Comparing the two-parameter family of numerical fluxes $f^{a_1,a_2}$ \eqref{eq:gnum-p}
with the one-parameter family $f^{a_1}$ \eqref{eq:fnum-p}, the additional parameter
$a_2$ contributes only to terms containing the velocity $v$. Thus, it is irrelevant
for well-balancing, since these terms vanish for the lake-at-rest initial condition.

Together with the results of the previous subsection, this proves
\begin{lemma}
\label{lem:gnum-extended}
  The two-parameter family
  \begin{equation}
  \label{eq:gnum-extended}
  \begin{aligned}
    f^{a_1,a_2}_{h_{i,k}}
    =&
      \frac{ 3 - a_{1} }{8} \left( h_i v_i + h_k v_k \right)
    + \frac{ 1 + a_{1} }{8} \left( h_i v_k + h_k v_i \right)
    \widesmall{}{\\&}
    + \frac{ a_{1} + 3 a_{2} - 2 }{16 g} \left(
      v_i^{3} - v_i^{2} v_k - v_i v_k^{2} + v_k^{3}
    \right),
    \\
    f^{a_1,a_2}_{hv_{i,k}}
    =&
      \frac{1 + a_{1} }{8} g \left( h_i^{2} + h_k^{2} \right)
    + \frac{1 - a_{1} }{4} g h_i h_k
    - \frac{2 a_{1} + 3 a_{2} - 5 }{16} \left( h_i v_i^{2} + h_k v_k^{2} \right)
    \\&
    + \frac{2 a_{1} + 3 a_{2} - 1 }{16} \left( h_i v_k^{2} + h_k v_i^{2} \right)
    + \frac{h_i v_i v_k}{4}
    + \frac{h_k v_i v_k}{4}
    \widesmall{}{\\&}
    + \frac{ a_{1} + 3 a_{2} - 2 }{32 g} \left(
      v_i^{4} - 2 v_i^{2} v_k^{2} + v_k^{4}
    \right)
    + \frac{ a_{1} + 3 a_{2} - 2 }{16} (v_k - v_i)^2 (b_k - b_i)
    \\&
    + \frac{1}{4} g \left( \frac{3-a_1}{2} h_i + \frac{1+a_1}{2} h_k \right)
      (b_k - b_i),
  \end{aligned}
  \end{equation}
  is a family of entropy conservative and well-balanced extended numerical fluxes
  (including contributions of the source term $g h \partial_x b$) for the shallow
  water equations \eqref{eq:SWE} with general bottom topography $b$.
\end{lemma}

To the author's knowledge, this general family of entropy conservative and
well-balanced extended numerical fluxes has only been considered for the special
choices $a_2 = \frac{2 - a_1}{3}, a_1 \in \set{-1,1}$ before.

\begin{remark}
  Again, a one-parameter family \eqref{eq:fnum-extended} is given by the special
  choice $a_2 = \frac{2 - a_1}{3}$ as
  \begin{align*}
  \stepcounter{equation}\tag{\theequation}
  \label{eq:fnum-extended}
    f^{a_1}_{h_{i,k}}
    =&
      \frac{ 3 - a_1 }{8} \left( h_i v_i + h_k v_k \right)
    + \frac{ 1 + a_1 }{8} \left( h_i v_k + h_k v_i \right),
    \\
    f^{a_1}_{hv_{i,k}}
    =&
      \frac{1 + a_1 }{8} g \left( h_i^{2} + h_k^{2} \right)
    + \frac{1 - a_1 }{4} g h_i h_k
    + \frac{ 3 - a_1 }{16} \left( h_i v_i^{2} + h_k v_k^{2} \right)
    \widesmall{}{\\&}
    + \frac{ 1 + a_1 }{16} \left( h_i v_k^{2} + h_k v_i^{2} \right)
    \widesmall{\\&}{}
    + \frac{h_i v_i v_k}{4}
    + \frac{h_k v_i v_k}{4}
    \widesmall{}{\\&}
    + \frac{1}{4} g \left( \frac{3-a_1}{2} h_i + \frac{1+a_1}{2} h_k \right)
      (b_k - b_i).
  \end{align*}
\end{remark}

\begin{remark}
  The volume terms corresponding to the two-parameter family of extended fluxes
  \eqref{eq:gnum-extended} are
  \begin{align*}
  \stepcounter{equation}\tag{\theequation}
  \label{eq:gnum-extended-volume-terms}
    \VOL^{a_1,a_2}_{h}
    \stackrel{\eqref{eq:gnum-volume-terms-h}}{=}&
      \frac{ 3 - a_{1} }{4} \mat{D} \vec{h v}
    + \frac{ 1 + a_{1} }{4} \left( \mat{h} \mat{D} \vec{v} + \mat{v} \mat{D} \vec{h} \right)
    \widesmall{}{\\&}
    + \frac{ a_{1} + 3 a_{2} - 2 }{8 g} \left(
      \mat{D} \vec{v^3} - \mat{v} \mat{D} \vec{v^2} - \mat{v^2} \mat{D} \vec{v} \right),
    \\
    \VOL^{a_1,a_2}_{hv}
    \stackrel{\eqref{eq:gnum-volume-terms-hv}}{=}&
      \frac{1 + a_{1} }{4} g \mat{D} \vec{h^2}
    + \frac{1 - a_{1} }{2} g \mat{h} \mat{D} \vec{h}
    - \frac{2 a_{1} + 3 a_{2} - 5 }{8} \mat{D} \vec{h v^2}
    \\&
    + \frac{2 a_{1} + 3 a_{2} - 1 }{8} \left( \mat{h} \mat{D} \vec{v^2}
      + \mat{v^2} \mat{D} \vec{h} \right)
    + \frac{1}{2} \left( \mat{h v} \mat{D} \vec{v} + \mat{v} \mat{D} \vec{h v} \right)
    \\&
    + \frac{ a_{1} + 3 a_{2} - 2 }{16 g} \left( \mat{D} \vec{v^4}
      - 2 \mat{v^2} \mat{D} \vec{v^2} \right)
    + \frac{3-a_1}{4} g \mat{h} \mat{D} \vec{b}
    \widesmall{}{\\&}
    + \frac{1+a_1}{4} g \left( \mat{D} \vec{h b} - \mat{b} \mat{D} \vec{h} \right)
    \\&
    + \frac{ a_{1} + 3 a_{2} - 2 }{8} \left(
          \mat{D} \vec{b v^2}
        - \mat{b} \mat{D} \vec{v^2}
        - 2 \mat{v} \mat{D} \vec{b v}
        +   \mat{v^2} \mat{D} \vec{b}
        + 2 \mat{b v} \mat{D} \vec{v}
      \right),
  \end{align*}
  where
  \begin{align*}
  \stepcounter{equation}\tag{\theequation}
    &
    \sum_{k=0}^p 2 D_{i,k} \Bigg(
      \frac{1}{4} g \left( \frac{3-a_1}{2} h_i + \frac{1+a_1}{2} h_k \right)(b_k - b_i)
      \widesmall{}{\\&\qquad\qquad\qquad}
      + \frac{ a_{1} + 3 a_{2} - 2 }{16} (v_k - v_i)^2 (b_k - b_i)
    \Bigg)
    \\
    =&
    \sum_{k=0}^p D_{i,k} \Bigg(
        \frac{3-a_1}{4} g h_i b_k
      + \frac{1+a_1}{4} g h_k (b_k - b_i)
      \\&\qquad\qquad
      + \frac{ a_{1} + 3 a_{2} - 2 }{8} \left(
          v_k^2 b_k - 2 v_i v_k b_k + v_i^2 b_k - v_k^2 b_i + 2 v_i b_i v_k \right)
    \Bigg)
    \\
    =&
    \Big[
        \frac{3-a_1}{4} g \mat{h} \mat{D} \vec{b}
      + \frac{1+a_1}{4} g \left( \mat{D} \vec{h b} - \mat{b} \mat{D} \vec{h} \right)
      \\&\quad
      + \frac{ a_{1} + 3 a_{2} - 2 }{8} \left(
        \mat{D} \vec{b v^2} - 2 \mat{v} \mat{D} \vec{b v} + \mat{v^2} \mat{D} \vec{b}
        - \mat{b} \mat{D} \vec{v^2} + 2 \mat{b v} \mat{D} \vec{v}
        \right)
    \Big]_i
  \end{align*}
  has been used. The corresponding surface terms using nodal bases including boundary
  nodes are simply given by
  \begin{equation*}
  \stepcounter{equation}\tag{\theequation}
  \label{eq:gnum-extended-surface-terms}
    \SURF_{h}
    =
    \mat{M}[^{-1}] \mat{R}[^T] \mat{B} \mat{R} \vec{h v},
    \;
    \SURF_{hv}
    =
    \mat{M}[^{-1}] \mat{R}[^T] \mat{B} \left(
      \mat{R} \vec{h v^2} + \frac{1}{2} g \mat{R} \vec{h^2}
    \right).
  \end{equation*}
  The numerical fluxes used for the volume discretisation and as surface flux may
  be combined arbitrarily, as done by \citet{gassner2016split, wintermeyer2016entropy},
  where they have used $f^{-1,1}$ as volume flux and $f^{1,\frac{1}{3}}$ as surface flux.
  Thus, a general semidiscretisation is of the form
  \begin{equation}
  \label{eq:gnum-semidisc}
  \begin{aligned}
    \partial_t \vec{h}
    =&
    - \VOL^{a_1,a_2}_{h}
    + \SURF_{h}
    - \mat{M}[^{-1}] \mat{R}[^T] \mat{B} \vec{f}^{b_1,b_2}_{h},
    \\
    \partial_t \vec{hv}
    =&
    - \VOL^{a_1,a_2}_{hv}
    + \SURF_{hv}
    - \mat{M}[^{-1}] \mat{R}[^T] \mat{B} \vec{f}^{b_1,b_2}_{hv},
  \end{aligned}
  \end{equation}
  with $a_1, a_2, b_1, b_1 \in \R$.
\end{remark}

\begin{remark}
  The volume terms \eqref{eq:gnum-extended-volume-terms} vanish for the lake-at-rest
  initial condition $h+b \equiv \mathrm{const}, hv = 0$. Since $\vec{v} = 0$,
  this is immediately clear for $\VOL^{a_1,a_2}_{h}$. For the discharge term
  $\VOL^{a_1,a_2}_{hv}$, using $\left( \mat{h}+\mat{b} \right) = (h+b) \mat{\I}$,
  $h+b \equiv \mathrm{const}$ = 0,
  \begin{align*}
  \stepcounter{equation}\tag{\theequation}
    &
    \frac{1 + a_{1} }{4} \mat{D} \vec{h^2}
    + \frac{1 - a_{1} }{2} \mat{h} \mat{D} \vec{h}
    + \frac{3-a_1}{4} \mat{h} \mat{D} \vec{b}
    + \frac{1+a_1}{4} \left( \mat{D} \vec{h b} - \mat{b} \mat{D} \vec{h} \right)
    \\
    =&
    \frac{1 + a_{1} }{4} \mat{D} \left( \mat{h}+\mat{b} \right) \vec{h}
    + \frac{1 - a_{1} }{2} \mat{h} \mat{D} \vec{h}
    + \frac{3-a_1}{4} \mat{h} \mat{D} \vec{b}
    - \frac{1+a_1}{4} \mat{b} \mat{D} \vec{h}
    \\
    =&
    \frac{1 + a_{1} }{4} \left( \mat{h}+\mat{b} \right) \mat{D} \vec{h}
    + \frac{1 - a_{1} }{2} \mat{h} \mat{D} \vec{h}
    + \frac{3-a_1}{4} \mat{h} \mat{D} \vec{b}
    - \frac{1+a_1}{4} \mat{b} \mat{D} \vec{h}
    \\
    =&
    \frac{3-a_1}{4} \mat{h} \mat{D} \vec{h+b}
    =
    0.
  \end{align*}
\end{remark}

\section{Extension to general SBP bases}
\label{sec:gauss}

In this section, an extension of the previous result to a nodal DG method using
Gauß nodes instead of Lobatto nodes or more general SBP bases will be
investigated.

Although the volume terms \eqref{eq:gnum-extended-volume-terms} have been
derived in section \ref{sec:source} with the assumption of a diagonal-norm SBP
basis including boundary nodes, they can be easily transferred to the setting of
a general SBP basis.
If the multiplication operators are self-adjoint with respect to the scalar product
induced by $\mat{M}$, e.g. for a nodal basis with diagonal mass matrix, then
the same volume terms \eqref{eq:gnum-extended-volume-terms} can be used.
Otherwise, some multiplication operators $\mat{a}$ have to be replaced by their
$\mat{M}$-adjoints $\mat{a}[^*] = \mat{M}[^{-1}] \mat{a}[^T] \mat{M}$, as proposed
by \citet{ranocha2015extended}. This results in the volume terms
\begin{align*}
\stepcounter{equation}\tag{\theequation}
\label{eq:gnum-extended-volume-terms-general}
  \VOL^{a_1,a_2}_{h}
  =&
    \frac{ 3 - a_{1} }{4} \mat{D} \vec{h v}
  + \frac{ 1 + a_{1} }{4} \left(
    \mat{h}[^*] \mat{D} \vec{v} + \mat{v}[^*] \mat{D} \vec{h} \right)
  \widesmall{}{\\&}
  + \frac{ a_{1} + 3 a_{2} - 2 }{8 g} \left(
    \mat{D} \vec{v^3} - \mat{v}[^*] \mat{D} \vec{v^2} - \mat{v^2}[^*] \mat{D} \vec{v}
    \right),
  \\
  \VOL^{a_1,a_2}_{hv}
  =&
    \frac{1 + a_{1} }{4} g \mat{D} \vec{h^2}
  + \frac{1 - a_{1} }{2} g \mat{h}[^*] \mat{D} \vec{h}
  - \frac{2 a_{1} + 3 a_{2} - 5 }{8} \mat{D} \vec{h v^2}
  \\&
  + \frac{2 a_{1} + 3 a_{2} - 1 }{8} \left( \mat{h}[^*] \mat{D} \vec{v^2}
    + \mat{v^2}[^*] \mat{D} \vec{h} \right)
  + \frac{1}{2} \left( \mat{h v}[^*] \mat{D} \vec{v} + \mat{v}[^*] \mat{D} \vec{h v} \right)
  \\&
  + \frac{ a_{1} + 3 a_{2} - 2 }{16 g} \left( \mat{D} \vec{v^4}
    - 2 \mat{v^2}[^*] \mat{D} \vec{v^2} \right)
  + \frac{3-a_1}{4} g \mat{h}[^*] \mat{D} \vec{b}
  \widesmall{}{\\&}
  + \frac{1+a_1}{4} g \left( \mat{D} \vec{h b} - \mat{b}[^*] \mat{D} \vec{h} \right)
  \\&
  + \frac{ a_{1} + 3 a_{2} - 2 }{8} \left(
        \mat{D} \vec{b v^2}
      - \mat{b}[^*] \mat{D} \vec{v^2}
      - 2 \mat{v}[^*] \mat{D} \vec{b v}
      +   \mat{v^2}[^*] \mat{D} \vec{b}
      + 2 \mat{b v}[^*] \mat{D} \vec{v}
    \right)\!.
\end{align*}
However, the surface terms \eqref{eq:gnum-extended-surface-terms} also have to be
adapted to a general basis. Often, the split form of the volume terms is described
as some correction for the product rule that does not hold discretely. However,
as described by \citet{ranocha2016summation}, it is the multiplication that is
not correct on a discrete level, resulting in an invalid product rule. Moreover,
if no boundary nodes are included in the basis, this inexactness also has to be
compensated in the surface terms. Thus, in the same spirit as the split form of
the volume terms can be seen as corrections to inexact multiplication, some kind
of correction has to be used for the interpolation to the boundaries.

Investigating \emph{conservation} (across elements), the time derivatives of the
conserved variables \eqref{eq:gnum-semidisc} are multiplied with $\vec{1}^T \mat{M}$,
corresponding to integration over an element. This yields for the volume terms
\eqref{eq:gnum-extended-volume-terms-general}
\begin{align*}
\stepcounter{equation}\tag{\theequation}
\label{eq:surf-term-from-vol-cons-h}
  \vec{1}^T \mat{M} \VOL^{a_1,a_2}_{h}
  =&
    \frac{ 3 - a_{1} }{4} \vec{1}^T \mat{M} \mat{D} \vec{h v}
  + \frac{ 1 + a_{1} }{4} \vec{1}^T \mat{M} \left(
    \mat{h}[^*] \mat{D} \vec{v} + \mat{v}[^*] \mat{D} \vec{h} \right)
  \\&
  + \frac{ a_{1} + 3 a_{2} - 2 }{8 g} \vec{1}^T \mat{M} \left(
    \mat{D} \vec{v^3} - \mat{v}[^*] \mat{D} \vec{v^2} - \mat{v^2}[^*] \mat{D} \vec{v}
    \right)
  \\
  =&
    \frac{ 3 - a_{1} }{4} \vec{1}^T \mat{M} \mat{D} \vec{h v}
  + \frac{ 1 + a_{1} }{4} \left(  \vec{h}^T \mat{M} \mat{D} \vec{v}
                                + \vec{v}^T \mat{M} \mat{D} \vec{h} \right)
  \\&
  + \frac{ a_{1} + 3 a_{2} - 2 }{8 g} \left(
    \vec{1}^T \mat{M} \mat{D} \vec{v^3} - \vec{v}^T \mat{M} \mat{D} \vec{v^2}
    - \vec{v^2}^T \mat{M} \mat{D} \vec{v}
    \right)
  \\
  \stackrel{\text{SBP }}{=}&
    \frac{ 3 - a_{1} }{4} \vec{1}^T \mat{R}[^T] \mat{B} \mat{R} \vec{h v}
  + \frac{ 1 + a_{1} }{4} \vec{h}^T \mat{R}[^T] \mat{B} \mat{R} \vec{v}
  \\&
  + \frac{ a_{1} + 3 a_{2} - 2 }{8 g} \left(
    \vec{1}^T \mat{R}[^T] \mat{B} \mat{R} \vec{v^3}
    - \vec{v}^T \mat{R}[^T] \mat{B} \mat{R} \vec{v^2}
    \right).
\end{align*}
Here, $\mat{h} \vec{1} = \vec{h}$, the SBP property $\mat{M} \mat{D} =
\mat{R}[^T] \mat{B} \mat{R} - \mat{D}[^T] \mat{M}$, and $\mat{D} \vec{1} = 0$
have been used. If multiplication and restriction to the boundary commute,
these volume terms are simply $\vec{1}^T \mat{R}[^T] \mat{B} \mat{R} \vec{h v}$
and yield the desired integral form.
Similarly,
\begin{align*}
\stepcounter{equation}\tag{\theequation}
\label{eq:surf-term-from-vol-cons-hv}
  &
  \vec{1}^T \mat{M} \VOL^{a_1,a_2}_{hv}
  \\
  =&
    \frac{1 + a_{1} }{4} g \vec{1}^T \mat{M} \mat{D} \vec{h^2}
  + \frac{1 - a_{1} }{2} g \vec{1}^T \mat{M} \mat{h}[^*] \mat{D} \vec{h}
  - \frac{2 a_{1} + 3 a_{2} - 5 }{8} \vec{1}^T \mat{M} \mat{D} \vec{h v^2}
  \\&
  + \frac{2 a_{1} + 3 a_{2} - 1 }{8} \vec{1}^T \mat{M} \left( \mat{h}[^*] \mat{D} \vec{v^2}
    + \mat{v^2}[^*] \mat{D} \vec{h} \right)
  + \frac{1}{2} \vec{1}^T \mat{M}  \left( \mat{h v}[^*] \mat{D} \vec{v}
    + \mat{v}[^*] \mat{D} \vec{h v} \right)
  \\&
  + \frac{ a_{1} + 3 a_{2} - 2 }{16 g} \vec{1}^T \mat{M} \left( \mat{D} \vec{v^4}
    - 2 \mat{v^2}[^*] \mat{D} \vec{v^2} \right)
  + \frac{3-a_1}{4} g \vec{1}^T \mat{M} \mat{h}[^*] \mat{D} \vec{b}
  \\&
  + \frac{1+a_1}{4} g \vec{1}^T \mat{M} \left( \mat{D} \vec{h b}
    - \mat{b}[^*] \mat{D} \vec{h} \right)
  \\&
  + \frac{ a_{1} + 3 a_{2} - 2 }{8} \vec{1}^T \mat{M} \left(
        \mat{D} \vec{b v^2}
      - \mat{b}[^*] \mat{D} \vec{v^2}
      - 2 \mat{v}[^*] \mat{D} \vec{b v}
      +   \mat{v^2}[^*] \mat{D} \vec{b}
      + 2 \mat{b v}[^*] \mat{D} \vec{v}
    \right)
  \\
  =&
    \frac{1 + a_{1} }{4} g \vec{1}^T \mat{M} \mat{D} \vec{h^2}
  + \frac{1 - a_{1} }{2} g \vec{h}^T \mat{M} \mat{D} \vec{h}
  - \frac{2 a_{1} + 3 a_{2} - 5 }{8} \vec{1}^T \mat{M} \mat{D} \vec{h v^2}
  \\&
  + \frac{2 a_{1} + 3 a_{2} - 1 }{8} \left( \vec{h}^T \mat{M} \mat{D} \vec{v^2}
    + \vec{v^2}^T \mat{M} \mat{D} \vec{h} \right)
  \widesmall{}{\\&}
  + \frac{1}{2} \left( \vec{h v}^T \mat{M} \mat{D} \vec{v}
                      + \vec{v}^T \mat{M} \mat{D} \vec{h v} \right)
  \\&
  + \frac{ a_{1} + 3 a_{2} - 2 }{16 g} \left( \vec{1}^T \mat{M} \mat{D} \vec{v^4}
    - 2 \vec{v^2}^T \mat{M} \mat{D} \vec{v^2} \right)
  + \frac{3-a_1}{4} g \vec{h}^T \mat{M} \mat{D} \vec{b}
  \\&
  + \frac{1+a_1}{4} g \left( \vec{1}^T \mat{M} \mat{D} \vec{h b}
    - \vec{b}^T \mat{M} \mat{D} \vec{h} \right)
  \widesmall{}{\\&}
  + \frac{ a_{1} + 3 a_{2} - 2 }{8} \left(
        \vec{1}^T \mat{M} \mat{D} \vec{b v^2}
      - \vec{b}^T \mat{M}  \mat{D} \vec{v^2}
    \right)
  \\&
  - \frac{ a_{1} + 3 a_{2} - 2 }{8} \left(
        2 \vec{v}^T \mat{M}  \mat{D} \vec{b v}
      -   \vec{v^2}^T \mat{M}  \mat{D} \vec{b}
      - 2 \vec{b v}^T \mat{M}  \mat{D} \vec{v}
    \right)
  \\
  \stackrel{\text{SBP }}{=}&
    \frac{1 + a_{1} }{4} g \vec{1}^T \mat{R}[^T] \mat{B} \mat{R} \vec{h^2}
  + \frac{1 - a_{1} }{4} g \vec{h}^T \mat{R}[^T] \mat{B} \mat{R} \vec{h}
  \widesmall{}{\\&}
  - \frac{2 a_{1} + 3 a_{2} - 5 }{8} \vec{1}^T \mat{R}[^T] \mat{B} \mat{R} \vec{h v^2}
  \widesmall{\\&}{}
  + \frac{2 a_{1} + 3 a_{2} - 1 }{8} \vec{h}^T \mat{R}[^T] \mat{B} \mat{R} \vec{v^2}
  \widesmall{}{\\&}
  + \frac{1}{2} \vec{v}^T \mat{R}[^T] \mat{B} \mat{R} \vec{h v}
  + \frac{ a_{1} + 3 a_{2} - 2 }{16 g} \left[
        \vec{1}^T \mat{R}[^T] \mat{B} \mat{R} \vec{v^4}
      - \vec{v^2}^T \mat{R}[^T] \mat{B} \mat{R} \vec{v^2}
    \right]
  \\&
  + \left\{ g \vec{h}^T \mat{M} \mat{D} \vec{b} \right\}
  - \frac{ a_{1} + 3 a_{2} - 2 }{4} \left\{
        \vec{v}^T \mat{M}  \mat{D} \vec{b v}
      - \vec{v^2}^T \mat{M}  \mat{D} \vec{b}
      - \vec{b v}^T \mat{M}  \mat{D} \vec{v}
    \right\}
  \\&
  + \frac{1+a_1}{4} g \left[ \vec{1}^T \mat{R}[^T] \mat{B} \mat{R} \vec{h b}
    - \vec{b}^T \mat{R}[^T] \mat{B} \mat{R} \vec{h} \right]
  \widesmall{}{\\&}
  + \frac{ a_{1} + 3 a_{2} - 2 }{8} \left[
        \vec{1}^T \mat{R}[^T] \mat{B} \mat{R} \vec{b v^2}
      - \vec{b}^T \mat{R}[^T] \mat{B} \mat{R} \vec{v^2}
    \right].
\end{align*}
Here, the terms in squared brackets $[\cdot]$ vanish if restriction to the boundary and
multiplication commute, i.e. for a basis using Lobatto nodes. However, for other
bases using e.g. Gauß nodes, these contributions are not zero in general.
The first term in curly brackets $\{\cdot\}$ is a consistent discretisation of the
source term $\int g h \partial_x b$. The second term in curly brackets $\{\cdot\}$ 
vanishes, if the product rule is valid, e.g. for constant bottom topography $b$.
However, for general bottom topography, it is not of the desired form for the
source influence $\int g h \partial_x b$ and it might be better to set it to
zero by choosing $a_2 = \frac{2 - a_1}{3}$, i.e. only the one-parameter family
instead of the two-parameter family.

These surface terms obtained in \eqref{eq:surf-term-from-vol-cons-h} and
\eqref{eq:surf-term-from-vol-cons-hv} have to be balanced by the surface terms
of the SBP SAT semidiscretisation \eqref{eq:gnum-semidisc} in order to get the
desired result
\begin{equation}
\label{eq:surf-term-cons-condition}
\begin{aligned}
  \vec{1}^T \mat{M} \partial_t \vec{h}
  =&
  - \vec{1}^T \mat{R}[^T] \mat{B} \vec{f}^{b_1,b_2}_{h},
  \\
  \vec{1}^T \mat{M} \partial_t \vec{h v}
  =&
  - \vec{1}^T \mat{R}[^T] \mat{B} \vec{f}^{b_1,b_2}_{hv}
  + \text{consistent contribution of} - \int g h \partial_x b,
\end{aligned}
\end{equation}
leading to a conservative scheme.

Turning to \emph{stability}, the approximation of
\begin{equation}
  \int \partial_t U
  =
  \int U'(u) \cdot \partial_t u
  =
  \int w \cdot \partial_t u
  ,\qquad
  w = \begin{pmatrix} g(h+b)-\frac{1}{2}v^2,v\end{pmatrix}^T,
\end{equation}
influenced by the volume terms is given by
\begin{align*}
\stepcounter{equation}\tag{\theequation}
\label{eq:surf-term-from-vol-stab}
  &
  \vec{w}_1^T \mat{M} \VOL^{a_1,a_2}_{h} + \vec{w}_2^T \mat{M} \VOL^{a_1,a_2}_{hv}
  \\
  =&
  % insert
    \frac{ 3 - a_{1} }{4} g \left( \vec{h} + \vec{b} \right)^T \mat{M} \mat{D} \vec{h v}
  + \frac{ 1 + a_{1} }{4} g \left( \vec{h} + \vec{b} \right)^T \mat{M} \left(
    \mat{h}[^*] \mat{D} \vec{v} + \mat{v}[^*] \mat{D} \vec{h} \right)
  \\&
  + \frac{ a_{1} + 3 a_{2} - 2 }{8} \left( \vec{h} + \vec{b} \right)^T \mat{M} \left(
    \mat{D} \vec{v^3} - \mat{v}[^*] \mat{D} \vec{v^2} - \mat{v^2}[^*] \mat{D} \vec{v}
    \right)
  \widesmall{}{\\&}
  - \frac{ 3 - a_{1} }{8} \vec{v^2}^T \mat{M} \mat{D} \vec{h v}
  \widesmall{\\&}{}
  - \frac{ 1 + a_{1} }{8} \vec{v^2}^T \mat{M} \left(
    \mat{h}[^*] \mat{D} \vec{v} + \mat{v}[^*] \mat{D} \vec{h} \right)
  \widesmall{}{\\&}
  - \frac{ a_{1} + 3 a_{2} - 2 }{16 g} \vec{v^2}^T \mat{M} \left(
    \mat{D} \vec{v^3} - \mat{v}[^*] \mat{D} \vec{v^2} - \mat{v^2}[^*] \mat{D} \vec{v}
    \right)
  \\&
  + \frac{1 + a_{1} }{4} g \vec{v}^T \mat{M} \mat{D} \vec{h^2}
  + \frac{1 - a_{1} }{2} g \vec{v}^T \mat{M} \mat{h}[^*] \mat{D} \vec{h}
  - \frac{2 a_{1} + 3 a_{2} - 5 }{8} \vec{v}^T \mat{M} \mat{D} \vec{h v^2}
  \\&
  + \frac{2 a_{1} + 3 a_{2} - 1 }{8} \vec{v}^T \mat{M} \left( \mat{h}[^*] \mat{D} \vec{v^2}
    + \mat{v^2}[^*] \mat{D} \vec{h} \right)
  + \frac{1}{2} \vec{v}^T \mat{M} \left( \mat{h v}[^*] \mat{D} \vec{v}
    + \mat{v}[^*] \mat{D} \vec{h v} \right)
  \\&
  + \frac{ a_{1} + 3 a_{2} - 2 }{16 g} \vec{v}^T \mat{M} \left( \mat{D} \vec{v^4}
    - 2 \mat{v^2}[^*] \mat{D} \vec{v^2} \right)
  + \frac{3-a_1}{4} g \vec{v}^T \mat{M} \mat{h}[^*] \mat{D} \vec{b}
  \\&
  + \frac{1+a_1}{4} g \vec{v}^T \mat{M} \left( \mat{D} \vec{h b}
    - \mat{b}[^*] \mat{D} \vec{h} \right)
  \\&
  + \frac{ a_{1} + 3 a_{2} - 2 }{8} \vec{v}^T \mat{M} \left(
        \mat{D} \vec{b v^2}
      - \mat{b}[^*] \mat{D} \vec{v^2}
      - 2 \mat{v}[^*] \mat{D} \vec{b v}
      +   \mat{v^2}[^*] \mat{D} \vec{b}
      + 2 \mat{b v}[^*] \mat{D} \vec{v}
    \right)
  \\
  % expand
  =&
    \frac{ 3 - a_{1} }{4} g \left( \vec{h} + \vec{b} \right)^T \mat{M} \mat{D} \vec{h v}
  + \frac{ 1 + a_{1} }{4} g \left( \vec{h^2} + \vec{bh} \right)^T \mat{M} \mat{D} \vec{v}
  \widesmall{}{\\&}
  + \frac{ 1 + a_{1} }{4} g \left( \vec{h v} + \vec{b v} \right)^T \mat{M} \mat{D} \vec{h}
  \widesmall{\\&}{}
  + \frac{ a_{1} + 3 a_{2} - 2 }{8} \left(\vec{h} + \vec{b}\right)^T \mat{M} \mat{D} \vec{v^3}
  \widesmall{}{\\&}
  - \frac{ a_{1} + 3 a_{2} - 2 }{8} \left(\vec{hv} + \vec{bv}\right)^T \mat{M} \mat{D} \vec{v^2}
  \widesmall{\\&}{}
  - \frac{ a_{1} + 3 a_{2} - 2 }{8} \left(\vec{hv^2} + \vec{bv^2}\right)^T \mat{M} \mat{D} \vec{v}
  \widesmall{}{\\&}
  - \frac{ 3 - a_{1} }{8} \vec{v^2}^T \mat{M} \mat{D} \vec{h v}
  - \frac{ 1 + a_{1} }{8} \vec{h v^2}^T \mat{M} \mat{D} \vec{v}
  \widesmall{\\&}{}
  - \frac{ 1 + a_{1} }{8} \vec{v^3}^T \mat{M} \mat{D} \vec{h}
  \widesmall{}{\\&}
  - \frac{ a_{1} + 3 a_{2} - 2 }{16 g} \vec{v^2}^T \mat{M} \mat{D} \vec{v^3}
  + \frac{ a_{1} + 3 a_{2} - 2 }{16 g} \vec{v^3}^T \mat{M} \mat{D} \vec{v^2}
  \\&
  + \frac{ a_{1} + 3 a_{2} - 2 }{16 g} \vec{v^4}^T \mat{M} \mat{D} \vec{v}
  + \frac{1 + a_{1} }{4} g \vec{v}^T \mat{M} \mat{D} \vec{h^2}
  + \frac{1 - a_{1} }{2} g \vec{h v}^T \mat{M} \mat{D} \vec{h}
  \\&
  - \frac{2 a_{1} + 3 a_{2} - 5 }{8} \vec{v}^T \mat{M} \mat{D} \vec{h v^2}
  + \frac{2 a_{1} + 3 a_{2} - 1 }{8} \vec{h v}^T \mat{M} \mat{D} \vec{v^2}
  \widesmall{}{\\&}
  + \frac{2 a_{1} + 3 a_{2} - 1 }{8} \vec{v^3}^T \mat{M} \mat{D} \vec{h}
  \widesmall{\\&}{}
  + \frac{1}{2} \vec{h v^2}^T \mat{M} \mat{D} \vec{v}
  + \frac{1}{2} \vec{v^2}^T \mat{M} \mat{D} \vec{h v}
  \widesmall{}{\\&}
  + \frac{ a_{1} + 3 a_{2} - 2 }{16 g} \vec{v}^T \mat{M} \mat{D} \vec{v^4}
  \widesmall{\\&}{}
  - \frac{ a_{1} + 3 a_{2} - 2 }{8 g} \vec{v^3}^T \mat{M} \mat{D} \vec{v^2}
  \widesmall{}{\\&}
  + \frac{3-a_1}{4} g \vec{h v}^T \mat{M} \mat{D} \vec{b}
  + \frac{1+a_1}{4} g \vec{v}^T \mat{M} \mat{D} \vec{h b}
  \\&
  - \frac{1+a_1}{4} g \vec{b v}^T \mat{M} \mat{D} \vec{h}
  + \frac{ a_{1} + 3 a_{2} - 2 }{8} \vec{v}^T \mat{M} \mat{D} \vec{b v^2}
  \widesmall{}{\\&}
  - \frac{ a_{1} + 3 a_{2} - 2 }{8} \vec{b v}^T \mat{M} \mat{D} \vec{v^2}
  \widesmall{\\&}{}
  - \frac{ a_{1} + 3 a_{2} - 2 }{4} \vec{v^2 }^T \mat{M} \mat{D} \vec{b v}
  \widesmall{}{\\&}
  + \frac{ a_{1} + 3 a_{2} - 2 }{8} \vec{v^3}^T \mat{M} \mat{D} \vec{b}
  \widesmall{\\&}{}
  + \frac{ a_{1} + 3 a_{2} - 2 }{4} \vec{b v^2}^T \mat{M} \mat{D} \vec{v}
  \\
  % sort
  =&
    \frac{ 3 - a_{1} }{4} g \left(
      \vec{h}^T \mat{M} \mat{D} \vec{h v} + \vec{h v}^T \mat{M} \mat{D} \vec{h}
    \right)
  + \frac{ 3 - a_{1} }{4} g \left(
      \vec{b}^T \mat{M} \mat{D} \vec{h v} + \vec{h v}^T \mat{M} \mat{D} \vec{b}
    \right)
  \\&
  + \frac{ 1 + a_{1} }{4} g \!\left(
      \vec{h^2}^T \mat{M} \mat{D} \vec{v} + \vec{v}^T \mat{M} \mat{D} \vec{h^2}
    \right)
  + \frac{ 1 + a_{1} }{4} g \!\left(
      \vec{b h}^T \mat{M} \mat{D} \vec{v} + \vec{v}^T \mat{M} \mat{D} \vec{b h}
    \right)
  \\&
  + 0 \, \vec{b v}^T \mat{M} \mat{D} \vec{h}
  + \frac{ a_{1} + 3 a_{2} - 2 }{8} \left(
      \vec{h}^T \mat{M} \mat{D} \vec{v^3} + \vec{v^3}^T \mat{M} \mat{D} \vec{h}
    \right)
  \widesmall{}{\\&}
  + \frac{ a_{1} + 1 }{8} \left(
      \vec{h v}^T \mat{M} \mat{D} \vec{v^2} + \vec{v^2}^T \mat{M} \mat{D} \vec{h v}
    \right)
  \\&
  + \frac{ 5 - 2 a_{1} - 3 a_{2} }{8} \left(
      \vec{h v^2}^T \mat{M} \mat{D} \vec{v} + \vec{v}^T \mat{M} \mat{D} \vec{h v^2}
    \right)
  \widesmall{}{\\&}
  - \frac{ a_{1} + 3 a_{2} - 2 }{16 g} \left(
      \vec{v^2}^T \mat{M} \mat{D} \vec{v^3} + \vec{v^3}^T \mat{M} \mat{D} \vec{v^2}
    \right)
  \\&
  + \frac{ a_{1} + 3 a_{2} - 2 }{16 g} \left(
      \vec{v^4}^T \mat{M} \mat{D} \vec{v} + \vec{v}^T \mat{M} \mat{D} \vec{v^4}
    \right)
  \widesmall{}{\\&}
  + \frac{ a_{1} + 3 a_{2} - 2 }{8} \left(
      \vec{b v^2}^T \mat{M} \mat{D} \vec{v} + \vec{v}^T \mat{M} \mat{D} \vec{b v^2}
    \right)
  \\&
  - \frac{ a_{1} + 3 a_{2} - 2 }{4} \left(
      \vec{b v}^T \mat{M} \mat{D} \vec{v^2} + \vec{v^2 }^T \mat{M} \mat{D} \vec{b v}
    \right)
  \widesmall{}{\\&}
  + \frac{ a_{1} + 3 a_{2} - 2 }{8} \left(
      \vec{b}^T \mat{M} \mat{D} \vec{v^3} + \vec{v^3}^T \mat{M} \mat{D} \vec{b}
    \right)
  \\
  %SBP
  \stackrel{\text{SBP }}{=}&
    \frac{ 3 - a_{1} }{4} g \vec{h}^T \mat{R}[^T] \mat{B} \mat{R} \vec{h v}
  + \frac{ 3 - a_{1} }{4} g \vec{b}^T \mat{R}[^T] \mat{B} \mat{R} \vec{h v}
  + \frac{ 1 + a_{1} }{4} g \vec{v}^T \mat{R}[^T] \mat{B} \mat{R} \vec{h^2}
  \\&
  + \frac{ 1 + a_{1} }{4} g \vec{v}^T \mat{R}[^T] \mat{B} \mat{R} \vec{b h}
  + \frac{ a_{1} + 3 a_{2} - 2 }{8} \vec{h}^T \mat{R}[^T] \mat{B} \mat{R} \vec{v^3}
  \widesmall{}{\\&}
  + \frac{ a_{1} + 1 }{8} \vec{h v}^T \mat{R}[^T] \mat{B} \mat{R} \vec{v^2}
  \widesmall{\\&}{}
  + \frac{ 5 - 2 a_{1} - 3 a_{2} }{8} \vec{v}^T \mat{R}[^T] \mat{B} \mat{R} \vec{h v^2}
  \widesmall{}{\\&}
  - \frac{ a_{1} + 3 a_{2} - 2 }{16 g} \vec{v^2}^T \mat{R}[^T] \mat{B} \mat{R} \vec{v^3}
  + \frac{ a_{1} + 3 a_{2} - 2 }{16 g} \vec{v}^T \mat{R}[^T] \mat{B} \mat{R} \vec{v^4}
  \\&
  + \frac{ a_{1} + 3 a_{2} - 2 }{8} \vec{v}^T \mat{R}[^T] \mat{B} \mat{R} \vec{b v^2}
  - \frac{ a_{1} + 3 a_{2} - 2 }{4} \vec{b v}^T \mat{R}[^T] \mat{B} \mat{R} \vec{v^2}
  \widesmall{}{\\&}
  + \frac{ a_{1} + 3 a_{2} - 2 }{8} \vec{b}^T \mat{R}[^T] \mat{B} \mat{R} \vec{v^3}.
\end{align*}
If multiplication and restriction to the boundary commute, these terms simplify
to $\vec{1}^T \mat{R}[^T] \mat{B} \mat{R} \vec{F}$, where
$\vec{F} = g \vec{h^2 v} + g \vec{b h v} + \frac{1}{2} \vec{h v^3}$ is the entropy
flux \eqref{eq:entropy-flux}.

These surface contributions resulting from the volume terms have to be balanced
by the surface terms $\vec{w}_1^T \mat{M} \SURF_{h} + \vec{w}_2^T \mat{M} \SURF_{hv}$,
in order to get an estimate of the form $\jump{w} \cdot \fnum - \jump{\psi}$
for the entropy change influenced by one boundary node (if the bottom topography
is continuous across elements). That is, the simple interpolations
\begin{equation}
  \mat{M}[^{-1}] \mat{R}[^T] \mat{B} \mat{R} \vec{h v},
  \qquad
  \mat{M}[^{-1}] \mat{R}[^T] \mat{B} \mat{R} \vec{h v^2},
  \qquad
  \mat{M}[^{-1}] \mat{R}[^T] \mat{B} \mat{R} \vec{h^2},
\end{equation}
in the surface terms \eqref{eq:gnum-extended-surface-terms} for the method
including boundary nodes have to be adapted.

The following combination of surface term structures proposed by
\citet{ranocha2016summation, ortleb2016kinetic} will be investigated
\begin{align*}
\stepcounter{equation}\tag{\theequation}
\label{eq:gnum-extended-surface-terms-general-ansatz}
  &
  \SURF^{a_1,a_2}_{h}
  \\
  =&
  % h v
    b_1 \mat{M}[^{-1}] \mat{R}[^T] \mat{B} \mat{R} \vec{h v}
  + b_2 \mat{M}[^{-1}] \mat{R}[^T] \mat{B} (\mat{R} \vec{h}) (\mat{R} \vec{v})
  + b_3 \mat{h}[^*] \mat{M}[^{-1}] \mat{R}[^T] \mat{B} \mat{R} \vec{v}
  \widesmall{}{\\&}
  + b_4 \mat{v}[^*] \mat{M}[^{-1}] \mat{R}[^T] \mat{B} \mat{R} \vec{h}
  \\&
  % 0/g v^3
  + \frac{c_1}{g} \mat{M}[^{-1}] \mat{R}[^T] \mat{B} \mat{R} \vec{v^3}
  + \frac{c_2}{g} \mat{M}[^{-1}] \mat{R}[^T] \mat{B} (\mat{R} \vec{v}) (\mat{R} \vec{v^2})
  + \frac{c_3}{g} \mat{v}[^*] \mat{M}[^{-1}] \mat{R}[^T] \mat{B} \mat{R} \vec{v^2}
  \\&
  + \frac{c_4}{g} \mat{v^2}[^*] \mat{M}[^{-1}] \mat{R}[^T] \mat{B} \mat{R} \vec{v},
  \\
  \phantom{\SURF}
  \\
  &
  \SURF^{a_1,a_2}_{hv}
  \\
  =&
  % h v^2
    d_1 \mat{M}[^{-1}] \mat{R}[^T] \mat{B} \mat{R} \vec{h v^2}
  + d_2 \mat{M}[^{-1}] \mat{R}[^T] \mat{B} (\mat{R} \vec{h}) (\mat{R} \vec{v^2})
  \widesmall{}{\\&}
  + d_3 \mat{M}[^{-1}] \mat{R}[^T] \mat{B} (\mat{R} \vec{h v}) (\mat{R} \vec{v})
  \widesmall{\\&}{}
  + d_4 \mat{M}[^{-1}] \mat{R}[^T] \mat{B} (\mat{R} \vec{h}) (\mat{R} \vec{v})^2
  \widesmall{}{\\&}
  + d_5 \mat{v}[^*] \mat{M}[^{-1}] \mat{R}[^T] \mat{B} \mat{R} \vec{h v}
  + d_6 \mat{h v}[^*] \mat{M}[^{-1}] \mat{R}[^T] \mat{B} \mat{R} \vec{v}
  \\&
  + d_7 \mat{v^2}[^*] \mat{M}[^{-1}] \mat{R}[^T] \mat{B} \mat{R} \vec{h}
  + d_8 \mat{h}[^*] \mat{M}[^{-1}] \mat{R}[^T] \mat{B} \mat{R} \vec{v^2}
  \\&
  % 1/2 g h^2
  + e_1 g \mat{M}[^{-1}] \mat{R}[^T] \mat{B} \mat{R} \vec{h^2}
  + e_2 g \mat{M}[^{-1}] \mat{R}[^T] \mat{B} ( \mat{R} \vec{h} )^2
  + e_3 g \mat{h}[^*] \mat{M}[^{-1}] \mat{R}[^T] \mat{B} \mat{R} \vec{h}
  \\&
  % 0/g v^4
  + \frac{k_1}{g} \mat{M}[^{-1}] \mat{R}[^T] \mat{B} \mat{R} \vec{v^4}
  + \frac{k_2}{g} \mat{M}[^{-1}] \mat{R}[^T] \mat{B} (\mat{R} \vec{v}) (\mat{R} \vec{v^3})
  + \frac{k_3}{g} \mat{M}[^{-1}] \mat{R}[^T] \mat{B} (\mat{R} \vec{v^2})^2
  \\&
  + \frac{k_4}{g} \mat{M}[^{-1}] \mat{R}[^T] \mat{B} (\mat{R} \vec{v})^2 (\mat{R} \vec{v^2})
  + \frac{k_5}{g} \mat{M}[^{-1}] \mat{R}[^T] \mat{B} (\mat{R} \vec{v})^4
  \widesmall{}{\\&}
  + \frac{k_6}{g} \mat{v}[^*] \mat{M}[^{-1}] \mat{R}[^T] \mat{B} \mat{R} \vec{v^3}
  \widesmall{\\&}{}
  + \frac{k_7}{g} \mat{v}[^*] \mat{M}[^{-1}] \mat{R}[^T] \mat{B}(\mat{R}\vec{v})(\mat{R}\vec{v^2})
  \widesmall{}{\\&}
  + \frac{k_8}{g} \mat{v}[^*] \mat{M}[^{-1}] \mat{R}[^T] \mat{B} (\mat{R} \vec{v})^3
  + \frac{k_9}{g} \mat{v^2}[^*] \mat{M}[^{-1}] \mat{R}[^T] \mat{B} \mat{R} \vec{v^2}
  \\&
  + \frac{k_{10}}{g} \mat{v^2}[^*] \mat{M}[^{-1}] \mat{R}[^T] \mat{B} (\mat{R} \vec{v})^2
  + \frac{k_{11}}{g} \mat{v^3}[^*] \mat{M}[^{-1}] \mat{R}[^T] \mat{B} \mat{R} \vec{v}
  \\&
  % 0 b v^2
  + l_1 \mat{M}[^{-1}] \mat{R}[^T] \mat{B} \mat{R} \vec{b v^2}
  + l_2 \mat{M}[^{-1}] \mat{R}[^T] \mat{B} (\mat{R} \vec{b}) (\mat{R} \vec{v^2})
  \widesmall{}{\\&}
  + l_3 \mat{M}[^{-1}] \mat{R}[^T] \mat{B} (\mat{R} \vec{b v}) (\mat{R} \vec{v})
  \widesmall{\\&}{}
  + l_4 \mat{M}[^{-1}] \mat{R}[^T] \mat{B} (\mat{R} \vec{b}) (\mat{R} \vec{v})^2
  \widesmall{}{\\&}
  + l_5 \mat{b}[^*] \mat{M}[^{-1}] \mat{R}[^T] \mat{B} \mat{R} \vec{v^2}
  + l_6 \mat{v^2}[^*] \mat{M}[^{-1}] \mat{R}[^T] \mat{B} \mat{R} \vec{b}
  \\&
  + l_7 \mat{b v}[^*] \mat{M}[^{-1}] \mat{R}[^T] \mat{B} \mat{R} \vec{v}
  + l_8 \mat{v}[^*] \mat{M}[^{-1}] \mat{R}[^T] \mat{B} \mat{R} \vec{b v}
  \\&
  + l_9 \mat{b}[^*] \mat{M}[^{-1}] \mat{R}[^T] \mat{B} (\mat{R} \vec{v})^2
  + l_{10} \mat{v}[^*] \mat{M}[^{-1}] \mat{R}[^T] \mat{B} (\mat{R} \vec{b}) (\mat{R} \vec{v})
  \\&
  % 0 g h b
  + m_1 g \mat{M}[^{-1}] \mat{R}[^T] \mat{B} \mat{R} \vec{b h}
  + m_2 g \mat{M}[^{-1}] \mat{R}[^T] \mat{B} (\mat{R} \vec{b}) (\mat{R} \vec{h})
  \\&
  + m_3 g \mat{h}[^*] \mat{M}[^{-1}] \mat{R}[^T] \mat{B} \mat{R} \vec{b}
  + m_4 g \mat{b}[^*] \mat{M}[^{-1}] \mat{R}[^T] \mat{B} \mat{R} \vec{h}
  \\&
  % correction for psi(Rh, Rv)
  - \frac{1}{2} \mat{v}[^*] \mat{M}[^{-1}] \mat{R}[^T] \mat{B} \vec{f}^{a_1,a_2}_{h}
  + \frac{1}{2} \mat{M}[^{-1}] \mat{R}[^T] \mat{B} (\vec{f}^{a_1,a_2}_{h}) (\mat{R} \vec{v}),
\end{align*}
where $b_i, c_i, d_i, e_i, k_i, l_i, m_i \in \R$ are free parameters that have
to be determined.

Considering \emph{conservation} for $h$, the relevant conditions are obtained by
multiplying the surface terms with $\vec{1}^T \mat{M}$.
\begin{align*}
\stepcounter{equation}\tag{\theequation}
  &
  \vec{1}^T \mat{M} \SURF^{a_1,a_2}_{h}
  \\
  =&
  % h v
    b_1 \vec{1}^T \mat{R}[^T] \mat{B} \mat{R} \vec{h v}
  + b_2 \vec{1}^T \mat{R}[^T] \mat{B} (\mat{R} \vec{h}) (\mat{R} \vec{v})
  + b_3 \vec{h}^T \mat{R}[^T] \mat{B} \mat{R} \vec{v}
  + b_4 \vec{v}^T \mat{R}[^T] \mat{B} \mat{R} \vec{h}
  \\&
  % 0/g v^3
  + \frac{c_1}{g} \vec{1}^T \mat{R}[^T] \mat{B} \mat{R} \vec{v^3}
  + \frac{c_2}{g} \vec{1}^T \mat{R}[^T] \mat{B} (\mat{R} \vec{v}) (\mat{R} \vec{v^2})
  + \frac{c_3}{g} \vec{v}^T \mat{R}[^T] \mat{B} \mat{R} \vec{v^2}
  \widesmall{}{\\&}
  + \frac{c_4}{g} \vec{v^2}^T \mat{R}[^T] \mat{B} \mat{R} \vec{v}
  \\
  =&
    b_1 \vec{1}^T \mat{R}[^T] \mat{B} \mat{R} \vec{h v}
  + \left( b_2 + b_3 + b_4 \right)
      \vec{1}^T \mat{R}[^T] \mat{B} (\mat{R} \vec{h}) (\mat{R} \vec{v})
  \\&
  + c_1 \frac{1}{g} \vec{1}^T \mat{R}[^T] \mat{B} \mat{R} \vec{v^3}
  + \left( c_2 + c_3 + c_4 \right)
      \vec{1}^T \mat{R}[^T] \mat{B} (\mat{R} \vec{v}) (\mat{R} \vec{v^2}).
\end{align*}
Here, some manipulations as
$\vec{h}^T \mat{R}[^T] \mat{B} \mat{R} \vec{v} = h_R v_R - h_L v_L
= \vec{1}^T \mat{R}[^T] \mat{B} (\mat{R} \vec{h}) (\mat{R} \vec{v})$
proposed by \citet{ranocha2016summation} have been used.
Thus, comparison with \eqref{eq:surf-term-from-vol-cons-h} yields the conditions
\begin{align*}
\stepcounter{equation}\tag{\theequation}
\label{eq:linear-system-SURF-cons-h}
  b_1
  =&
  \frac{3-a_1}{4},
  \qquad&
  b_2 + b_3 + b_4 
  =&
  \frac{1+a_1}{4},
  \\
  c_1
  =&
  \frac{a_1+3a_2-2}{8},
  \qquad&
  c_2 + c_3 + c_4
  =&
  - \frac{a_1+3a_2-2}{8}.
\end{align*}

Similarly, for $h v$,
\begin{align*}
\stepcounter{equation}\tag{\theequation}
  &
  \vec{1}^T \mat{M} \SURF^{a_1,a_2}_{hv}
  \\
  =&
  % h v^2
    d_1 \vec{1}^T \mat{R}[^T] \mat{B} \mat{R} \vec{h v^2}
  + d_2 \vec{1}^T \mat{R}[^T] \mat{B} (\mat{R} \vec{h}) (\mat{R} \vec{v^2})
  + d_3 \vec{1}^T \mat{R}[^T] \mat{B} (\mat{R} \vec{h v}) (\mat{R} \vec{v})
  \widesmall{}{\\&}
  + d_4 \vec{1}^T \mat{R}[^T] \mat{B} (\mat{R} \vec{h}) (\mat{R} \vec{v})^2
  \widesmall{\\&}{}
  + d_5 \vec{v}^T \mat{R}[^T] \mat{B} \mat{R} \vec{h v}
  + d_6 \vec{h v}^T \mat{R}[^T] \mat{B} \mat{R} \vec{v}
  \widesmall{}{\\&}
  + d_7 \vec{v^2}^T \mat{R}[^T] \mat{B} \mat{R} \vec{h}
  + d_8 \vec{h}^T \mat{R}[^T] \mat{B} \mat{R} \vec{v^2}
  \\&
  % 1/2 g h^2
  + e_1 g \vec{1}^T \mat{R}[^T] \mat{B} \mat{R} \vec{h^2}
  + e_2 g \vec{1}^T \mat{R}[^T] \mat{B} ( \mat{R} \vec{h} )^2
  + e_3 g \vec{h}^T \mat{R}[^T] \mat{B} \mat{R} \vec{h}
  \\&
  % 0/g v^4
  + \frac{k_1}{g} \vec{1}^T \mat{R}[^T] \mat{B} \mat{R} \vec{v^4}
  + \frac{k_2}{g} \vec{1}^T \mat{R}[^T] \mat{B} (\mat{R} \vec{v}) (\mat{R} \vec{v^3})
  + \frac{k_3}{g} \vec{1}^T \mat{R}[^T] \mat{B} (\mat{R} \vec{v^2})^2
  \\&
  + \frac{k_4}{g} \vec{1}^T \mat{R}[^T] \mat{B} (\mat{R} \vec{v})^2 (\mat{R} \vec{v^2})
  + \frac{k_5}{g} \vec{1}^T \mat{R}[^T] \mat{B} (\mat{R} \vec{v})^4
  + \frac{k_6}{g} \vec{v}^T \mat{R}[^T] \mat{B} \mat{R} \vec{v^3}
  \\&
  + \frac{k_7}{g} \vec{v}^T \mat{R}[^T] \mat{B}(\mat{R}\vec{v})(\mat{R}\vec{v^2})
  + \frac{k_8}{g} \vec{v}^T \mat{R}[^T] \mat{B} (\mat{R} \vec{v})^3
  + \frac{k_9}{g} \vec{v^2}^T \mat{R}[^T] \mat{B} \mat{R} \vec{v^2}
  \\&
  + \frac{k_{10}}{g} \vec{v^2}^T \mat{M}[^{-1}] \mat{R}[^T] \mat{B} (\mat{R} \vec{v})^2
  + \frac{k_{11}}{g} \vec{v^3}^T \mat{M}[^{-1}] \mat{R}[^T] \mat{B} \mat{R} \vec{v}
  \\&
  % 0 b v^2
  + l_1 \vec{1}^T \mat{R}[^T] \mat{B} \mat{R} \vec{b v^2}
  + l_2 \vec{1}^T \mat{R}[^T] \mat{B} (\mat{R} \vec{b}) (\mat{R} \vec{v^2})
  + l_3 \vec{1}^T \mat{R}[^T] \mat{B} (\mat{R} \vec{b v}) (\mat{R} \vec{v})
  \widesmall{}{\\&}
  + l_4 \vec{1}^T \mat{R}[^T] \mat{B} (\mat{R} \vec{b}) (\mat{R} \vec{v})^2
  \widesmall{\\&}{}
  + l_5 \vec{b}^T \mat{R}[^T] \mat{B} \mat{R} \vec{v^2}
  + l_6 \vec{v^2}^T \mat{R}[^T] \mat{B} \mat{R} \vec{b}
  + l_7 \vec{b v}^T \mat{R}[^T] \mat{B} \mat{R} \vec{v}
  \\&
  + l_8 \vec{v}^T \mat{R}[^T] \mat{B} \mat{R} \vec{b v}
  + l_9 \vec{b}^T \mat{R}[^T] \mat{B} (\mat{R} \vec{v})^2
  + l_{10} \vec{v}^T \mat{R}[^T] \mat{B} (\mat{R} \vec{b}) (\mat{R} \vec{v})
  \\&
  % 0 g h b
  + m_1 g \vec{1}^T \mat{R}[^T] \mat{B} \mat{R} \vec{b h}
  + m_2 g \vec{1}^T \mat{R}[^T] \mat{B} (\mat{R} \vec{b}) (\mat{R} \vec{h})
  \widesmall{}{\\&}
  + m_3 g \vec{h}^T \mat{R}[^T] \mat{B} \mat{R} \vec{b}
  + m_4 g \vec{b}^T \mat{R}[^T] \mat{B} \mat{R} \vec{h}.
\end{align*}
Analogously, comparing this with \eqref{eq:surf-term-from-vol-cons-hv} results
in the conditions
\begin{align*}
\stepcounter{equation}\tag{\theequation}
\label{eq:linear-system-SURF-cons-hv}
  d_1
  =&
  \frac{5-2a_1-3a_2}{8},
  \qquad&
  d_2 + d_7 + d_8
  =&
  \frac{2a_1+3a_2-1}{8},
  \\
  d_3 + d_5 + d_6
  =&
  \frac{1}{2},
  \qquad&
  d_4
  =&
  0,
  \\
  e_1
  =&
  \frac{1+a_1}{4},
  \qquad&
  e_2 + e_3
  =&
  \frac{1-a_1}{4},
  \\
  % v^4
  k_1
  =&
  \frac{a_1+3a_2-2}{16},
  \qquad&
  % v, v^3
  k_2 + k_6 + k_{11}
  =&
  0,
  \\
  % v^2, v^2
  k_3 + k_9
  =&
  -\frac{a_1+3a_2-2}{16},
  \qquad&
  % v, v, v^2
  k_4 + k_7 + k_{10}
  =&
  0,
  \\
  % v, v, v, v
  k_5 + k_8
  =&
  0,
  \qquad&
  l_1
  =&
  \frac{a_1+3a_2-2}{8},
  \\
  l_2 + l_5 + l_6
  =&
  -\frac{a_1+3a_2-2}{8},
  \qquad&
  l_3 + l_7 + l_8
  =&
  0,
  \\
  l_4 + l_9 + l_{10}
  =&
  0,
  \qquad&
  m_1
  =&
  \frac{1+a_1}{4},
  \\
  m_2 + m_3 + m_4
  =&
  -\frac{1+a_1}{4}.
\end{align*}

Considering \emph{stability}, the surface terms
\eqref{eq:gnum-extended-surface-terms-general-ansatz} yield
\begin{align*}
\stepcounter{equation}\tag{\theequation}
  &
  \vec{w}_1^T \mat{M} \SURF^{a_1,a_2}_{h} + \vec{w}_2^T \mat{M} \SURF^{a_1,a_2}_{hv}
  \\
  =&
    b_1 g \vec{h}^T \mat{R}[^T] \mat{B} \mat{R} \vec{h v}
  + b_2 g \vec{h}^T \mat{R}[^T] \mat{B} (\mat{R} \vec{h}) (\mat{R} \vec{v})
  + b_3 g \vec{h^2}^T \mat{R}[^T] \mat{B} \mat{R} \vec{v}
  \widesmall{}{\\&}
  + b_4 g \vec{h v}^T \mat{R}[^T] \mat{B} \mat{R} \vec{h}
  \widesmall{\\&}{}
  + c_1 \vec{h}^T \mat{R}[^T] \mat{B} \mat{R} \vec{v^3}
  + c_2 \vec{h}^T \mat{R}[^T] \mat{B} (\mat{R} \vec{v}) (\mat{R} \vec{v^2})
  \widesmall{}{\\&}
  + c_3 \vec{h v}^T \mat{R}[^T] \mat{B} \mat{R} \vec{v^2}
  + c_4 \vec{h v^2}^T \mat{R}[^T] \mat{B} \mat{R} \vec{v}
  \widesmall{\\&}{}
  + b_1 g \vec{b}^T \mat{R}[^T] \mat{B} \mat{R} \vec{h v}
  \widesmall{}{\\&}
  + b_2 g \vec{b}^T \mat{R}[^T] \mat{B} (\mat{R} \vec{h}) (\mat{R} \vec{v})
  + b_3 g \vec{b h}^T \mat{R}[^T] \mat{B} \mat{R} \vec{v}
  + b_4 g \vec{b v}^T \mat{R}[^T] \mat{B} \mat{R} \vec{h}
  \\&
  + c_1 \vec{b}^T \mat{R}[^T] \mat{B} \mat{R} \vec{v^3}
  + c_2 \vec{b}^T \mat{R}[^T] \mat{B} (\mat{R} \vec{v}) (\mat{R} \vec{v^2})
  + c_3 \vec{b v}^T \mat{R}[^T] \mat{B} \mat{R} \vec{v^2}
  \widesmall{}{\\&}
  + c_4 \vec{b v^2}^T \mat{R}[^T] \mat{B} \mat{R} \vec{v}
  \widesmall{\\&}{}
  -\frac{1}{2} b_1 \vec{v^2}^T \mat{R}[^T] \mat{B} \mat{R} \vec{h v}
  -\frac{1}{2} b_2 \vec{v^2}^T \mat{R}[^T] \mat{B} (\mat{R} \vec{h}) (\mat{R} \vec{v})
  \widesmall{}{\\&}
  -\frac{1}{2} b_3 \vec{h v^2}^T \mat{R}[^T] \mat{B} \mat{R} \vec{v}
  -\frac{1}{2} b_4 \vec{v^3}^T \mat{R}[^T] \mat{B} \mat{R} \vec{h}
  \\&
  - \frac{c_1}{2 g} \vec{v^2}^T \mat{R}[^T] \mat{B} \mat{R} \vec{v^3}
  - \frac{c_2}{2 g} \vec{v^2}^T \mat{R}[^T] \mat{B} (\mat{R} \vec{v}) (\mat{R} \vec{v^2})
  - \frac{c_3}{2 g} \vec{v^3}^T \mat{R}[^T] \mat{B} \mat{R} \vec{v^2}
  \widesmall{}{\\&}
  - \frac{c_4}{2 g} \vec{v^4}^T \mat{R}[^T] \mat{B} \mat{R} \vec{v}
  \widesmall{\\&}{}
  + d_1 \vec{v}^T \mat{R}[^T] \mat{B} \mat{R} \vec{h v^2}
  + d_2 \vec{v}^T \mat{R}[^T] \mat{B} (\mat{R} \vec{h}) (\mat{R} \vec{v^2})
  \widesmall{}{\\&}
  + d_3 \vec{v}^T \mat{R}[^T] \mat{B} (\mat{R} \vec{h v}) (\mat{R} \vec{v})
  + d_4 \vec{v}^T \mat{R}[^T] \mat{B} (\mat{R} \vec{h}) (\mat{R} \vec{v})^2
  \\&
  + d_5 \vec{v^2}^T \mat{R}[^T] \mat{B} \mat{R} \vec{h v}
  + d_6 \vec{h v^2}^T \mat{R}[^T] \mat{B} \mat{R} \vec{v}
  + d_7 \vec{v^3}^T \mat{R}[^T] \mat{B} \mat{R} \vec{h}
  + d_8 \vec{h v}^T \mat{R}[^T] \mat{B} \mat{R} \vec{v^2}
  \\&
  + e_1 g \vec{v}^T \mat{R}[^T] \mat{B} \mat{R} \vec{h^2}
  + e_2 g \vec{v}^T \mat{R}[^T] \mat{B} ( \mat{R} \vec{h} )^2
  + e_3 g \vec{h v}^T \mat{R}[^T] \mat{B} \mat{R} \vec{h}
  \\&
  + \frac{k_1}{g} \vec{v}^T \mat{R}[^T] \mat{B} \mat{R} \vec{v^4}
  + \frac{k_2}{g} \vec{v}^T \mat{R}[^T] \mat{B} (\mat{R} \vec{v}) (\mat{R} \vec{v^3})
  + \frac{k_3}{g} \vec{v}^T \mat{R}[^T] \mat{B} (\mat{R} \vec{v^2})^2
  \\&
  + \frac{k_4}{g} \vec{v}^T \mat{R}[^T] \mat{B} (\mat{R} \vec{v})^2 (\mat{R} \vec{v^2})
  + \frac{k_5}{g} \vec{v}^T \mat{R}[^T] \mat{B} (\mat{R} \vec{v})^4
  + \frac{k_6}{g} \vec{v^2}^T \mat{R}[^T] \mat{B} \mat{R} \vec{v^3}
  \\&
  + \frac{k_7}{g} \vec{v^2}^T \mat{R}[^T] \mat{B}(\mat{R}\vec{v})(\mat{R}\vec{v^2})
  + \frac{k_8}{g} \vec{v^2}^T \mat{R}[^T] \mat{B} (\mat{R} \vec{v})^3
  + \frac{k_9}{g} \vec{v^3}^T \mat{R}[^T] \mat{B} \mat{R} \vec{v^2}
  \\&
  + \frac{k_{10}}{g} \vec{v^3}^T \mat{M}[^{-1}] \mat{R}[^T] \mat{B} (\mat{R} \vec{v})^2
  + \frac{k_{11}}{g} \vec{v^4}^T \mat{M}[^{-1}] \mat{R}[^T] \mat{B} \mat{R} \vec{v}
  \\&
  + l_1 \vec{v}^T \mat{R}[^T] \mat{B} \mat{R} \vec{b v^2}
  + l_2 \vec{v}^T \mat{R}[^T] \mat{B} (\mat{R} \vec{b}) (\mat{R} \vec{v^2})
  + l_3 \vec{v}^T \mat{R}[^T] \mat{B} (\mat{R} \vec{b v}) (\mat{R} \vec{v})
  \widesmall{}{\\&}
  + l_4 \vec{v}^T \mat{R}[^T] \mat{B} (\mat{R} \vec{b}) (\mat{R} \vec{v})^2
  \widesmall{\\&}{}
  + l_5 \vec{b v}^T \mat{R}[^T] \mat{B} \mat{R} \vec{v^2}
  + l_6 \vec{v^3}^T \mat{R}[^T] \mat{B} \mat{R} \vec{b}
  \widesmall{}{\\&}
  + l_7 \vec{b v^2}^T \mat{R}[^T] \mat{B} \mat{R} \vec{v}
  \widesmall{\\&}{}
  + l_8 \vec{v^2}^T \mat{R}[^T] \mat{B} \mat{R} \vec{b v}
  + l_9 \vec{b v}^T \mat{R}[^T] \mat{B} (\mat{R} \vec{v})^2
  \widesmall{}{\\&}
  + l_{10} \vec{v^2}^T \mat{R}[^T] \mat{B} (\mat{R} \vec{b}) (\mat{R} \vec{v})
  \\&
  + m_1 g \vec{v}^T \mat{R}[^T] \mat{B} \mat{R} \vec{b h}
  + m_2 g \vec{v}^T \mat{R}[^T] \mat{B} (\mat{R} \vec{b}) (\mat{R} \vec{h})
  + m_3 g \vec{h v}^T \mat{R}[^T] \mat{B} \mat{R} \vec{b}
  \widesmall{}{\\&}
  + m_4 g \vec{b v}^T \mat{R}[^T] \mat{B} \mat{R} \vec{h}
  \widesmall{\\&}{}
  - \frac{1}{2} \vec{v^2}^T \mat{R}[^T] \mat{B} \vec{f}^{a_1,a_2}_{h}
  + \frac{1}{2} \vec{v}^T \mat{R}[^T] \mat{B} (\vec{f}^{a_1,a_2}_{h}) (\mat{R} \vec{v})
  \\
  =& 
    \left( b_1 + b_4 + e_3 \right)
      g \vec{h}^T \mat{R}[^T] \mat{B} \mat{R} \vec{h v}
  + \left( b_2 \right)
      g \vec{h}^T \mat{R}[^T] \mat{B} (\mat{R} \vec{h}) (\mat{R} \vec{v})
  \widesmall{}{\\&}
  + \left( b_3 + e_1 \right)
      g \vec{h^2}^T \mat{R}[^T] \mat{B} \mat{R} \vec{v}
  \widesmall{\\&}{}
  \widesmall{}{\\&}
  + \left( c_1 -\frac{1}{2} b_4 + d_7 \right)
      \vec{h}^T \mat{R}[^T] \mat{B} \mat{R} \vec{v^3}
  \widesmall{}{\\&}
  + \left( c_2 -\frac{1}{2} b_2 + d_2 \right)
      \vec{h}^T \mat{R}[^T] \mat{B} (\mat{R} \vec{v}) (\mat{R} \vec{v^2})
  \\&
  + \left( c_3 -\frac{1}{2} b_1 + d_5 \right)
      \vec{h v}^T \mat{R}[^T] \mat{B} \mat{R} \vec{v^2}
  + \left( c_4 -\frac{1}{2} b_3 + d_1 + d_6 \right)
      \vec{h v^2}^T \mat{R}[^T] \mat{B} \mat{R} \vec{v}
  \\&
  + \left( b_1 + m_3 \right)
      g \vec{b}^T \mat{R}[^T] \mat{B} \mat{R} \vec{h v}
  + \left( b_2 + m_2 \right)
      g \vec{b}^T \mat{R}[^T] \mat{B} (\mat{R} \vec{h}) (\mat{R} \vec{v})
  \widesmall{}{\\&}
  + \left( b_3 + m_1 \right)
      g \vec{b h}^T \mat{R}[^T] \mat{B} \mat{R} \vec{v}
  \widesmall{\\&}{}
  + \left( b_4 + m_4 \right)
      g \vec{b v}^T \mat{R}[^T] \mat{B} \mat{R} \vec{h}
  \widesmall{}{\\&}
  + \left( c_1 + l_6 \right)
      \vec{b}^T \mat{R}[^T] \mat{B} \mat{R} \vec{v^3}
  \widesmall{}{\\&}
  + \left( c_2 + l_2 + l_{10} \right)
      \vec{b}^T \mat{R}[^T] \mat{B} (\mat{R} \vec{v}) (\mat{R} \vec{v^2})
  \\&
  + \left( c_3 + l_5 + l_8 \right)
      \vec{b v}^T \mat{R}[^T] \mat{B} \mat{R} \vec{v^2}
  + \left( c_4 + l_1 + l_7 \right)
      \vec{b v^2}^T \mat{R}[^T] \mat{B} \mat{R} \vec{v}
  \\&
  + \left( - \frac{1}{2} c_1 - \frac{1}{2} c_3 + k_6 + k_9 \right)
      \frac{1}{g} \vec{v^2}^T \mat{R}[^T] \mat{B} \mat{R} \vec{v^3}
  \widesmall{}{\\&}
  + \left( - \frac{1}{2} c_2 + k_3 + k_7 \right)
      \frac{1}{g} \vec{v^2}^T \mat{R}[^T] \mat{B} (\mat{R} \vec{v}) (\mat{R} \vec{v^2})
  \\&
  + \left( - \frac{1}{2} c_4 + k_1 + k_{11} \right)
      \frac{1}{g} \vec{v^4}^T \mat{R}[^T] \mat{B} \mat{R} \vec{v}
  + d_3 \vec{v}^T \mat{R}[^T] \mat{B} (\mat{R} \vec{h v}) (\mat{R} \vec{v})
  \widesmall{}{\\&}
  + d_4 \vec{v}^T \mat{R}[^T] \mat{B} (\mat{R} \vec{h}) (\mat{R} \vec{v})^2
  \widesmall{\\&}{}
  + d_8 \vec{h v}^T \mat{R}[^T] \mat{B} \mat{R} \vec{v^2}
  + e_2 g \vec{v}^T \mat{R}[^T] \mat{B} ( \mat{R} \vec{h} )^2
  \\&
  + \left( k_2 + k_{10} \right)
      \frac{1}{g} \vec{v}^T \mat{R}[^T] \mat{B} (\mat{R} \vec{v}) (\mat{R} \vec{v^3})
  + \left( k_4 + k_8 \right)
      \frac{1}{g} \vec{v}^T \mat{R}[^T] \mat{B} (\mat{R} \vec{v})^2 (\mat{R} \vec{v^2})
  \\&
  + k_5 \frac{1}{g} \vec{v}^T \mat{R}[^T] \mat{B} (\mat{R} \vec{v})^4
  + \left( l_3 + l_9 \right)
      \vec{v}^T \mat{R}[^T] \mat{B} (\mat{R} \vec{b v}) (\mat{R} \vec{v})
  + l_4 \vec{v}^T \mat{R}[^T] \mat{B} (\mat{R} \vec{b}) (\mat{R} \vec{v})^2
  \\&
  - \frac{1}{2} \vec{v^2}^T \mat{R}[^T] \mat{B} \vec{f}^{a_1,a_2}_{h}
  + \frac{1}{2} \vec{v}^T \mat{R}[^T] \mat{B} (\vec{f}^{a_1,a_2}_{h}) (\mat{R} \vec{v}).
\end{align*}
Comparing this with \eqref{eq:surf-term-from-vol-stab} yields the conditions
\begin{align*}
\stepcounter{equation}\tag{\theequation}
\label{eq:linear-system-SURF-stab}
  % (h, hv)
  b_1 + b_4 + e_3
  =&
  \frac{3-a_1}{4},
  \widesmall{\qquad&}{&}
  % (h, h, v)
  b_2 + e_2
  =&
  \frac{1}{2}, %because of \psi !
  \\
  % (h^2, v)
  b_3 + e_1
  =&
  \frac{1+a_1}{4},
  \widesmall{\qquad&}{&}
  % (h, v^3)
  c_1 - \frac{b_4}{2} + d_7
  =&
  \frac{a_1+3a_2-2}{8},
  \\
  % (h, v, v^2)
  c_2 - \frac{b_2}{2} + d_2
  =&
  0,
  \widesmall{\qquad&}{&}
  % (hv, v^2)
  c_3 - \frac{b_1}{2} + d_5 + d_8
  =&
  \frac{a_1+1}{8},
  \\
  % (hv^2, v)
  c_4 - \frac{b_3}{2} + d_1 + d_6
  =&
  \frac{5-2a_1-3a_2}{8},
  \widesmall{\qquad&}{&}
  % (b, hv)
  b_1 + m_3
  =&
  \frac{3-a_1}{4},
  \\
  % (b, h, v)
  b_2 + m_2
  =&
  0,
  \widesmall{\qquad&}{&}
  % (bh, v)
  b_3 + m_1
  =&
  \frac{a_1+1}{4},
  \\
  % (bv, h)
  b_4 + m_4
  =&
  0,
  \widesmall{\qquad&}{&}
  % (b, v^3)
  c_1 + l_6
  =&
  \frac{a_1+3a_2-2}{8},
  \\
  % (b, v, v^2)
  c_2 + l_2 + l_{10}
  =&
  0,
  \widesmall{\qquad&}{&}
  % (bv, v^2)
  c_3 + l_5 + l_8
  =&
  -\frac{a_1+3a_2-2}{4},
  \\
  % (bv^2, v)
  c_4 + l_1 + l_7
  =&
  \frac{a_1+3a_2-2}{8},
  \widesmall{\qquad&}{&}
  % (b v, v, v)
  l_4
  =&
  0,
  \\
  % (v^2, v^3)
  - \frac{c_1}{2} - \frac{c_3}{2} + k_6 + k_9
  =&
  -\frac{a_1+3a_2-2}{16},
  \widesmall{\qquad&}{&}
  % (v, v^2, v^2)
  - \frac{c_2}{2} + k_3 + k_7
  =&
  0,
  \\
  % (v, v^4)
  - \frac{c_4}{2} + k_1 + k_{11}
  =&
  \frac{a_1+3a_2-2}{16},
  \widesmall{\qquad&}{&}
  % (hv, v, v)
  d_3
  =&
  0,
  \\
  % (h, v, v)
  d_4
  =&
  0,
  \widesmall{\qquad&}{&}
  % (v, v, v^3)
  k_2 + k_{10}
  =&
  0,
  \\
  % (v, v, v, v^2)
  k_4 + k_8
  =&
  0,
  \widesmall{\qquad&}{&}
  % (v, v, v, v, v)
  k_5
  =&
  0
  \\
  % (bv, v, v)
  l_3 + l_9
  =&
  0,
  \widesmall{\qquad&}{&}
  % (b, v, v, v)
  l_4
  =&
  0.
\end{align*}
Solving the linear system given by \eqref{eq:linear-system-SURF-cons-h},
\eqref{eq:linear-system-SURF-cons-hv}, and \eqref{eq:linear-system-SURF-stab}
with SymPy \citep{sympy1} results in the free parameters
$m_{4}, k_{9}, k_{10}, k_{11}, l_{10}$ for any given parameters $a_{1}, a_{2}$:
\begin{align*}
\stepcounter{equation}\tag{\theequation}
\label{eq:linear-system-SURF-solution}
  b_{1} =& - \frac{a_{1}}{4} + \frac{3}{4},
  \widesmall{\qquad&}{&\!}
  b_{2} =& \frac{a_{1}}{4} + m_{4} + \frac{1}{4},
  \\
  b_{3} =& 0,
  \widesmall{\qquad&}{&\!}
  b_{4} =& - m_{4},
  \\
  c_{1} =& \frac{a_{1}}{8} + \frac{3 a_{2}}{8} - \frac{1}{4},
  \widesmall{\qquad&}{&\!}
  c_{2} =& - \frac{a_{1}}{8} - \frac{3 a_{2}}{8} - 2 k_{10} - 2 k_{9} + \frac{1}{4},
  \\
  c_{3} =& 2 k_{10} - 2 k_{11} + 2 k_{9},
  \widesmall{\qquad&}{&\!}
  c_{4} =& 2 k_{11},
  \\
  d_{1} =& - \frac{a_{1}}{4} - \frac{3 a_{2}}{8} + \frac{5}{8},
  \widesmall{\qquad&}{&\!}
  d_{2} =& \frac{2 a_{1} + 3 a_{2} - 1}{8} + 2 k_{10} + 2 k_{9} + \frac{m_{4}}{2},
  \\
  d_{3} =& 0,
  \widesmall{\qquad&}{&\!}
  d_{4} =& 0,
  \\
  d_{5} =& 2 k_{11} + \frac{1}{2},
  \widesmall{\qquad&}{&\!}
  d_{6} =& - 2 k_{11},
  \\
  d_{7} =& - \frac{m_{4}}{2},
  \widesmall{\qquad&}{&\!}
  d_{8} =& - 2 k_{10} - 2 k_{9},
  \\
  e_{1} =& \frac{a_{1}}{4} + \frac{1}{4},
  \widesmall{\qquad&}{&\!}
  e_{2} =& - \frac{a_{1}}{4} - m_{4} + \frac{1}{4},
  \\
  e_{3} =& m_{4},
  \widesmall{\qquad&}{&\!}
  k_{1} =& \frac{a_{1}}{16} + \frac{3 a_{2}}{16} - \frac{1}{8},
  \\
  k_{2} =& - k_{10},
  \widesmall{\qquad&}{&\!}
  k_{3} =& - \frac{a_{1}}{16} - \frac{3 a_{2}}{16} - k_{9} + \frac{1}{8},
  \\
  k_{4} =& 0,
  \widesmall{\qquad&}{&\!}
  k_{5} =& 0,
  \\
  k_{6} =& k_{10} - k_{11},
  \widesmall{\qquad&}{&\!}
  k_{7} =& - k_{10},
  \\
  k_{8} =& 0,
  \widesmall{\qquad&}{&\!}
  l_{1} =& \frac{a_{1}}{8} + \frac{3 a_{2}}{8} - \frac{1}{4},
  \\
  \widesmall{
    l_{2} =& \frac{a_{1} + 3 a_{2} - 2}{8} + 2 k_{10} + 2 k_{9} - l_{10},
    }{
    l_{3} =& l_{10},
    }
  \widesmall{\qquad&}{&\!}
  \widesmall{
    l_{3} =& l_{10},
    }{
    l_{2} =& \frac{a_{1} + 3 a_{2} - 2}{8} + 2 k_{10} + 2 k_{9} - l_{10},
    }
  \\
  l_{4} =& 0,
  \widesmall{\qquad&}{&\!}
  l_{5} =& l_{10} - \frac{a_{1} + 3 a_{2} - 2}{4} - 2 k_{10} - 2 k_{9},
  \\
  l_{6} =& 0,
  \widesmall{\qquad&}{&\!}
  l_{7} =& - 2 k_{11},
  \\
  l_{8} =& 2 k_{11} - l_{10},
  \widesmall{\qquad&}{&\!}
  l_{9} =& - l_{10},
  \\
  m_{1} =& \frac{a_{1}}{4} + \frac{1}{4},
  \widesmall{\qquad&}{&\!}
  m_{2} =& - \frac{a_{1}}{4} - m_{4} - \frac{1}{4},
  \\
  m_{3} =& 0.
\end{align*}
This proves the following
\begin{theorem}
\label{thm:gnum-extended-semidisc}
  For $a_1, a_2, \alpha_1, \alpha_2 \in \R$, using a general SBP operator,
  the semidiscretisation
  \begin{equation}
  \label{eq:gnum-extended-semidisc}
  \begin{aligned}
    \partial_t \vec{h}
    =&
    - \VOL^{a_1,a_2}_{h}
    + \SURF^{a_1,a_2}_{h}
    - \mat{M}[^{-1}] \mat{R}[^T] \mat{B} \vec{f}^{\alpha_1,\alpha_2}_{h},
    \\
    \partial_t \vec{hv}
    =&
    - \VOL^{a_1,a_2}_{hv}
    + \SURF^{a_1,a_2}_{hv}
    - \mat{M}[^{-1}] \mat{R}[^T] \mat{B} \vec{f}^{\alpha_1,\alpha_2}_{hv},
  \end{aligned}
  \end{equation}
  with volume terms \eqref{eq:gnum-extended-volume-terms-general}
  and surface terms \eqref{eq:gnum-extended-surface-terms-general-ansatz},
  where the parameters are chosen according to \eqref{eq:linear-system-SURF-solution}
  with free parameters $m_{4}, k_{9}, k_{10}, k_{11}, l_{10} \in \R$,
  \begin{enumerate}
    \item 
    conserves the total mass $\int h$. Additionally, it conserves the total
    momentum $\int hv$, if the bottom topography is constant. Otherwise, the
    rate of change is consistent with the source term $- g h \partial_x b$.
    
    \item
    conserves the total entropy / energy $\int U$.
    
    \item
    handles the lake-at-rest condition correctly.
  \end{enumerate}
  That is, this semidiscretisation is \emph{conservative} (across elements),
  \emph{stable} (entropy conservative), and \emph{well-balanced}.
\end{theorem}

\begin{remark}
  For an SBP operator using a nodal basis with diagonal mass matrix $\mat{M}$,
  the volume terms can be equivalently expressed in the flux difference form
  corresponding to the extended numerical fluxes \eqref{eq:gnum-extended}
  in Lemma \ref{lem:gnum-extended}.
\end{remark}

\begin{remark}
  For the special case $a_2 = \frac{2 - a_1}{3}$, i.e. the one-parameter family,
  the volume terms \eqref{eq:gnum-extended-volume-terms-general} can be considerably
  simplified. Accordingly, the ansatz \eqref{eq:gnum-extended-surface-terms-general-ansatz}
  can be simplified by setting the coefficients $c_i$, $k_i$, and $l_i$ to zero.
  In this case, only $m_4$ remains as free parameter. In this case, the surface terms
  \eqref{eq:gnum-extended-surface-terms-general-ansatz} become
  \begin{equation}
  \label{eq:gnum-extended-surface-terms-general-ansatz-one-parameter-family}
  \begin{aligned}
    \SURF^{a_1,\frac{2-a_1}{3}}_{h}
    =&
    % h v
      \frac{3-a_1}{4} \mat{M}[^{-1}] \mat{R}[^T] \mat{B} \mat{R} \vec{h v}
    + \frac{a_1+1+4m_4}{4} \mat{M}[^{-1}] \mat{R}[^T] \mat{B} (\mat{R} \vec{h}) (\mat{R} \vec{v})
    \\&
    - m_4 \mat{v}[^*] \mat{M}[^{-1}] \mat{R}[^T] \mat{B} \mat{R} \vec{h},
    \\
    \phantom{\SURF}
    \\
    \SURF^{a_1,\frac{2-a_1}{3}}_{hv}
    =&
    % h v^2
      \frac{3-a_1}{8} \mat{M}[^{-1}] \mat{R}[^T] \mat{B} \mat{R} \vec{h v^2}
    + \frac{a_1+1+4m_4}{8} \mat{M}[^{-1}] \mat{R}[^T] \mat{B} (\mat{R} \vec{h}) (\mat{R} \vec{v^2})
    \\&
    + \frac{1}{2} \mat{v}[^*] \mat{M}[^{-1}] \mat{R}[^T] \mat{B} \mat{R} \vec{h v}
    - \frac{m_4}{2} \mat{v^2}[^*] \mat{M}[^{-1}] \mat{R}[^T] \mat{B} \mat{R} \vec{h}
    % 1/2 g h^2
    \widesmall{}{\\&}
    + \frac{a_1+1}{4} g \mat{M}[^{-1}] \mat{R}[^T] \mat{B} \mat{R} \vec{h^2}
    \widesmall{\\&}{}
    + \frac{1-a_1-4m_4}{4} g \mat{M}[^{-1}] \mat{R}[^T] \mat{B} ( \mat{R} \vec{h} )^2
    \widesmall{}{\\&}
    + m_4 g \mat{h}[^*] \mat{M}[^{-1}] \mat{R}[^T] \mat{B} \mat{R} \vec{h}
    % 0 g h b
    \widesmall{\\&}{}
    + \frac{a_1+1}{4} g \mat{M}[^{-1}] \mat{R}[^T] \mat{B} \mat{R} \vec{b h}
    \widesmall{}{\\&}
    - \frac{a_1+1+4m_4}{4} g \mat{M}[^{-1}] \mat{R}[^T] \mat{B} (\mat{R} \vec{b}) (\mat{R} \vec{h})
    \widesmall{\\&}{}
    + m_4 g \mat{b}[^*] \mat{M}[^{-1}] \mat{R}[^T] \mat{B} \mat{R} \vec{h}
    % correction for psi(Rh, Rv)
    \widesmall{}{\\&}
    - \frac{1}{2} \mat{v}[^*] \mat{M}[^{-1}] \mat{R}[^T] \mat{B} \vec{f}^{a_1,a_2}_{h}
    + \frac{1}{2} \mat{M}[^{-1}] \mat{R}[^T] \mat{B} (\vec{f}^{a_1,a_2}_{h}) (\mat{R} \vec{v}).
  \end{aligned}
  \end{equation}
  Here, the choice $m_4 = 0$ results in the fewest terms. Additionally, choosing
  $a_1 = -1$ cancels some other terms and results in the skew-symmetric form of
  \citet{gassner2016well}.
\end{remark}

\begin{remark}
  In this one-dimensional setting using a linear grid, the surface correction terms
  are more complicated than the corresponding terms for Burgers' equation
  \citep{ranocha2016summation}. Thus, although an extension to Gauß nodes on
  two-dimensional curvilinear grids may seem possible, it would result in a vast
  amount of surface and volume terms. As described (in an earlier version) by
  \citet{wintermeyer2015entropy}, the use of the flux differencing form is
  clearly superior for the implementation. Although the usage of Gauß nodes may
  yield improved accuracy, the additional correction terms render the method
  inferior compared to the usage of Lobatto nodes and the flux differencing form.
\end{remark}

\begin{remark}
  The two-parameter family of fluxes \eqref{eq:gnum-w} has been derived in entropy
  variables $w$ and translated to primitive variables $h, v$ \eqref{eq:gnum-p}
  in the following calculations. Hence, it may seem natural to consider a splitting
  similar to \eqref{eq:gnum-extended-volume-terms-general} and
  \eqref{eq:gnum-extended-surface-terms-general-ansatz} but expressed using entropy
  instead of primitive variables. The volume terms can be obtained similarly to
  \eqref{eq:gnum-extended-volume-terms-general}, but the surface terms are more
  delicate. A general ansatz similar to \eqref{eq:gnum-extended-surface-terms-general-ansatz}
  without the correction $- \mat{v}[^*] \mat{M}[^{-1}] \mat{R}[^T] \mat{B} \vec{f}^{a_1,a_2}_{h}
  + \mat{M}[^{-1}] \mat{R}[^T] \mat{B} (\vec{f}^{a_1,a_2}_{h}) (\mat{R} \vec{v})$
  has not been successful, i.e. the resulting linear system was not solvable.
  However, it might be possible to get the desired result using another ansatz
  for the surface terms.
\end{remark}

\section{Numerical surface fluxes and positivity preservation}
\label{sec:fluxes}

In this section, several numerical surface fluxes will be presented, which can
be used to get entropy stable, positivity preserving, and well-balanced schemes.
Therefore, the framework of positivity preservation by \citet{zhang2011maximum}
will be presented and adapted to the setting of a nodal SBP method.

\subsection{Positivity preservation}
\label{sec:positivity}

The argumentation of \citet{zhang2011maximum} can be summarised as
\begin{enumerate}
  \item 
  Since the method should be conservative, ensure a non-negative mean value of
  $h$ in each cell.
  
  \item
  If the height becomes negative somewhere but the mean value is non-negative,
  use a suitable limiter to enforce the non-negativity as needed.
\end{enumerate}
In a pure finite volume framework using an explicit Euler step, in cell $i$
\begin{equation}
\label{eq:FV-mean-evolution}
  \overline{h}_i^+
  =
  \overline{h}_i - \frac{\Delta t}{\Delta x} \left(
    \fnum(\overline{u}_{i}, \overline{u}_{i+1}) - \fnum(\overline{u}_{i-1}, \overline{u}_{i})
   \right).
\end{equation}
If there are numerical fluxes such that under a suitable CFL condition
$\Delta t \leq c \Delta x$ the non-negativity can be guaranteed, this can be
extended to the DG setting as proved by \citet{zhang2011maximum} and references
cited by them as follows.

Choose $q+1$ Lobatto-Legendre nodes $x_k$ in the element $i$. Since Lobatto
quadrature with weight $\omega_k$ at $x_k$ is exact for polynomials of degree
$2q-1$, ensure $2q-1 \geq p$, i.e. $q \geq \frac{p-1}{2}$. If the finite
volume method \eqref{eq:FV-mean-evolution} using the numerical flux $\fnum$ is
positivity preserving under the CFL condition $\Delta t \leq c \Delta x$, the
DG method is positivity preserving under the scaled CFL condition
$\Delta t \leq \omega_q c \Delta x$, if the values of $h_{i,k}$ at the quadrature
nodes $x_0, \dots, x_1$ are non-negative.

Then, after the Euler step, the mean value in each cell $i$ is ensured to be
non-negative. Applying the simple linear scaling limiter
\begin{equation}
\label{eq:scaling-limiter}
  h_i
  \mapsto
  \tilde h_i
  :=
  \theta h_i + (1-\theta) \overline{h}_i
  =
  \overline{h}_i + \theta \left( h_i - \overline{h}_i \right),
  \quad
  \theta
  =
  \begin{cases}
    1, & \min h_i \geq 0,
    \\
    \frac{\overline{h}_i}{\overline{h}_i - \min h_i}, & \min h_i < 0,
  \end{cases}
\end{equation}
of \citet[Section 2.3]{liu1996nonoscillatory} ensures $\min h_i \geq 0$.
Here, $\min h_i$ can be computed as the minimum of the polynomial $h_i$ of degree
$\leq p$ in the complete cell, or as the minimum of $h_i$ at the nodes
$x_0, \dots, x_q$ used to guarantee the non-negativity of the cell mean.
Applying this limiter in each Euler sub-step, this procedure can be extended to
SSP methods consisting of convex combinations of Euler steps.

In a nodal collocation framework such as nodal DG or FD, it would be more natural
to enforce the non-negativity of the water height at the collocation nodes
$\xi_0, \dots, \xi_p$ in the standard element. Thus, if Lobatto nodes are used,
the framework of \citet{zhang2011maximum} described above extends directly,
where the non-negativity of the mean value $\overline{h}_i^+$ can be guaranteed
under the possibly worse CFL condition $\Delta t \leq \omega_p c \Delta x$,
where $p \geq q$, when the limiter \eqref{eq:scaling-limiter} is applied
with $\min h_i = \min \set{h_{i}(\xi_0), \dots, h_i(\xi_p)}$. Alternatively, to
get a better CFL condition $\Delta t \leq \omega_q c \Delta x$, the water height
can be interpolated to the other nodes $x_0, \dots, x_q$ and the limiter
\eqref{eq:scaling-limiter} can be applied with
$\min h_i = \min \set{h_{i}(x_0), \dots, h_{i}(x_q), h_{i}(\xi_0), \dots, h_i(\xi_p)}$.

Similarly, if Gauß nodes are considered, only one possibility is obvious:
Interpolate to suitable Lobatto nodes $x_0, \dots, x_q$ and use the limiter
\eqref{eq:scaling-limiter}, again with the choice of the minimum as minimum
over solution nodes and interpolation points
$\min h_i = \min \set{h_{i}(x_0), \dots, h_{i}(x_q), h_{i}(\xi_0), \dots, h_i(\xi_p)}$.

If this limiter should be coupled with an entropy conservative / stable method,
a natural question is whether the limiter \eqref{eq:scaling-limiter} is entropy
stable. This is indeed true, since for a convex entropy $U$
\begin{equation}
\begin{aligned}
  \overline{ U(\tilde h_i) }
  =&
  \overline{ \theta U(h_i) + (1-\theta) U(\overline{h}_i) }
  \\
  \stackrel{U \text{ convex}}{\leq}&
  \theta \overline{ U(h_i) } + (1-\theta) \overline{ U(\overline{h}_i) }
  =
  \theta \overline{ U(h_i) } + (1-\theta) U(\overline{h}_i)
  \\
  \stackrel{\text{Jensen}}{\leq}&
  \theta \overline{ U(h_i) } + (1-\theta) \overline{ U(h_i) }
  =
  \overline{ U(h_i) },
\end{aligned}
\end{equation}
where the monotonicity of the mean value, the convexity of $U$, and Jensen's
inequality have been used. This proves
\begin{lemma}
  The scaling limiter \eqref{eq:scaling-limiter} is entropy stable for
  $\theta \in [0,1]$.
\end{lemma}

\subsection{Entropy conservative fluxes for constant bottom topography \texorpdfstring{$b$}{b}}

Here, the one-parameter entropy conservative numerical flux $f_h^{a_1} =
\frac{1-a_1}{2} \mean{h v} + \frac{1+a_1}{2} \mean{h} \mean{v}$  \eqref{eq:fnum}
will be considered. Inserting this in a FV evolution equation \eqref{eq:FV-mean-evolution}
(and dropping the bars $\overline{\cdot}$),
\begin{align*}
\stepcounter{equation}\tag{\theequation}
\label{eq:fnum-EC-positivity}
  h_i^+
  =&
  h_i - \frac{\Delta t}{\Delta x} \left(
    f_h^{a_1}(u_{i}, u_{i+1}) - f_h^{a_1}(u_{i-1}, u_{i}) \right)
  \\
  =&
  h_i - \frac{\Delta t}{\Delta x} \Bigg(
      \frac{1-a_1}{2} \frac{h_{i} v_{i} + h_{i+1} v_{i+1}}{2}
      + \frac{1+a_1}{2} \frac{h_{i} + h_{i+1}}{2} \frac{v_{i} + v_{i+1}}{2}
    \\&\qquad\qquad\quad
    - \frac{1-a_1}{2} \frac{h_{i} v_{i} + h_{i-1} v_{i-1}}{2}
      - \frac{1+a_1}{2} \frac{h_{i} + h_{i-1}}{2} \frac{v_{i} + v_{i-1}}{2}
  \Bigg)
  \\
  =&
  h_{i} \left[
    1 - \frac{\Delta t}{\Delta x} \frac{1+a_1}{8} (v_{i+1} - v_{i-1})
  \right]
  - h_{i+1} \frac{\Delta t}{\Delta x} \left[
    \frac{3-a_1}{8} v_{i+1} + \frac{1+a_1}{8} v_{i}
  \right]
  \\&
  + h_{i-1} \frac{\Delta t}{\Delta x} \left[
    \frac{3-a_1}{8} v_{i-1} + \frac{1+a_1}{8} v_{i}
  \right].
\end{align*}
Considering non-negative water height $h$, if $h_i = 0$, $h_{i-1}, h_{i+1} > 0$,
$v_i = 0$, $(3-a_1) v_{i-1} < 0$, $(3-a_1) v_{i+1} > 0$, the height $h_i^+$
becomes negative for $\Delta t > 0$.

If only positive water height $h$ should be considered, using the same conditions
as before but $h_i > 0$, the new water height can be guaranteed to be positive,
but only under a CFL condition depending on $h_i$ with allowed $\Delta t \to 0$
as $h_i \to 0$.

\begin{lemma}
\label{lem:positivity-EC_flux}
  The one-parameter family of entropy conservative fluxes \eqref{eq:fnum} for
  the shallow water equations with constant bottom topography $b$ is not positivity
  preserving under a CFL condition not blowing up as $h \to 0$.
\end{lemma}

If the full two-parameter family \eqref{eq:gnum-p} of entropy conservative fluxes
is considered, terms of the form $v_{i\pm1}^3 - v_{i\pm1}^2 v_i - v_{i\pm1} v_i^2 + v_i^3$
have to be added. Since these terms may be of arbitrary size and sign and do not
contain any multiple of the water height $h$, they can render the water height
negative. Thus, this family is not positivity preserving, too.

\subsection{Adding dissipation for constant bottom topography \texorpdfstring{$b$}{b}}

Classically, adding "dissipation" to some kind of central flux $f^\mathrm{cent}$
would result in $f^\mathrm{cent} - P \jump{u}$,
where $P$ is positive semi-definite. However, this might not result in an entropy
stable scheme. Therefore, the amount of dissipation is better chosen as
proportional to the jump of entropy variables $\jump{w}$ instead of $\jump{u}$,
as used inter alia by \citet{barth1999numerical, fjordholm2011well, gassner2016well}.
In the simplest case, $P = \lambda \I$ with $\lambda \geq 0$. Then
$\lambda \jump{u} \approx \lambda \frac{\partial u}{\partial w} \jump{w}$.
Since $\partial_w u = \left( \partial_u w \right)^{-1}$, $\partial_w u = U''(u)$,
and $U$ is convex, this dissipation matrix is positive semi-definite.
Thus, the resulting numerical flux
\begin{equation}
\label{eq:fnum-LLF-type}
  \fnum
  =
  f^{a_1}
  - 
  \frac{\lambda}{2} \overline{\partial_w u} \jump{w}
\end{equation}
is entropy stable, where $f^{a_1}$ is the one-parameter entropy conservative
flux \eqref{eq:fnum} and $\overline{\partial_w u}$ is a suitable positive semi-definite
approximation of the entropy Jacobian $\partial_w u$, e.g. the Jacobian evaluated 
at some mean value $\overline{u}$.

For the shallow water equations \eqref{eq:SWE} (cf. \eqref{eq:entropy-jacobian}),
\begin{equation}
  \partial_u w
  =
  \begin{pmatrix}
    g + \frac{v^2}{h} & - \frac{v}{h}
    \\
    - \frac{v}{h}     & \frac{1}{h}
  \end{pmatrix},
  \qquad
  \partial_w u
  =
  \begin{pmatrix}
    \frac{1}{g} & \frac{v}{g}
    \\
    \frac{v}{g} & h + \frac{v^2}{g}
  \end{pmatrix}.
\end{equation}
Thus, adding dissipation in the form \eqref{eq:fnum-LLF-type}, the following
additional contribution has to be added to the right hand side of
\eqref{eq:fnum-EC-positivity} for the entropy variables $w = \begin{pmatrix}
g h - \frac{1}{2} v^2, v \end{pmatrix}^T$, i.e. if the bottom topography $b$
is continuous across cell boundaries
\begin{align*}
\stepcounter{equation}\tag{\theequation}
  &
  \frac{\Delta t}{2 \Delta x} \left[
    \lambda_{i+\frac{1}{2}} \overline{\partial_w u}_{i+\frac{1}{2}} (w_{i+1} - w_{i})
    - \lambda_{i-\frac{1}{2}} \overline{\partial_w u}_{i-\frac{1}{2}} (w_{i} - w_{i-1})
  \right]_{h}
  \\
  =&
  \frac{\Delta t}{\Delta x} \frac{\lambda_{i+\frac{1}{2}}}{2} \left(
      h_{i+1} - \frac{v_{i+1}^2}{2 g}
    - h_{i} + \frac{v_{i}^2}{2 g}
    + \frac{\overline{v}_{i+\frac{1}{2}}}{g} (v_{i+1} - v_{i})
    \right)
  \\&
  + \frac{\Delta t}{\Delta x} \frac{\lambda_{i-\frac{1}{2}}}{2} \left(
    - h_{i} + \frac{v_{i}^2}{2 g}
    + h_{i-1} - \frac{v_{i-1}^2}{2 g}
    - \frac{\overline{v}_{i-\frac{1}{2}}}{g} (v_{i} - v_{i-1})
  \right).
\end{align*}
The additional positive values of $h_{i\pm1}$ are crucial to obtain the positivity
preserving property under a suitable CFL condition, if $\lambda$ is large enough.
The negative values of $h_i$ are weighted by $\Delta t$ and can thus be bounded
by the positive terms in \eqref{eq:fnum-EC-positivity} if $\Delta t$ is small
enough. The remaining terms containing only the velocity but not the height of
the water may be problematic. However, if the Jacobian $\overline{\partial_u w}$
is evaluated at the arithmetic mean value, 
\begin{equation}
  \overline{v}_{i+\frac{1}{2}}
  =
  \mean{v}_{i,i+1}
  =
  \frac{v_{i+1} + v_{i}}{2}
  \quad \Rightarrow \quad
  \overline{v}_{i+\frac{1}{2}} (v_{i+1} - v_{i})
  =
  \frac{1}{2} v_{i+1}^2 - \frac{1}{2} v_{i}^2,
\end{equation}
and similarly for $\overline{v}_{i-\frac{1}{2}}$. Thus, these additional terms
vanish and (cf. equation \eqref{eq:fnum-EC-positivity})
\begin{align*}
\stepcounter{equation}\tag{\theequation}
\label{eq:fnum-LLF-type-positivity}
  h_i^+
  =&
  h_{i} \left[
    1 - \frac{\Delta t}{\Delta x} \frac{1+a_1}{8} (v_{i+1} - v_{i-1})
    - \frac{\lambda_{i+\frac{1}{2}} + \lambda_{i-\frac{1}{2}}}{2} \frac{\Delta t}{\Delta x}
  \right]
  \\&
  + h_{i+1} \frac{\Delta t}{\Delta x} \left[
    \frac{\lambda_{i+\frac{1}{2}}}{2}
    - \frac{3-a_1}{8} v_{i+1} - \frac{1+a_1}{8} v_{i}
  \right]
  \widesmall{}{\\&}
  + h_{i-1} \frac{\Delta t}{\Delta x} \left[
    \frac{\lambda_{i-\frac{1}{2}}}{2}
    + \frac{3-a_1}{8} v_{i-1} + \frac{1+a_1}{8} v_{i}
  \right].
\end{align*}
Hence, for $\lambda_{i \pm \frac{1}{2}} \geq \max \set{ \abs{v_{i}}, \abs{v_{i\pm1}} }$,
the coefficients of $h_{i\pm1}$ are non-negative, and under the CFL condition
\begin{equation}
\label{eq:fnum-LLF-type-positivity-CFL}
  \frac{\Delta t}{\Delta x} \left(
    \abs{\frac{1+a_1}{8}} (\abs{v_{i+1}} + \abs{v_{i-1}})
    + \frac{\lambda_{i+\frac{1}{2}} + \lambda_{i-\frac{1}{2}}}{2}
  \right)
  \leq
  1,
\end{equation}
the new height $h_i^+$ is non-negative. Note that this CFL condition does not blow
up as $h_i \to 0$. With this estimate, the choice $a_1 = -1$ seems to be optimal
in order to get the least restrictive CFL condition.
Again, considering the full two-parameter family \eqref{eq:gnum-p} instead of
the one-parameter family \eqref{eq:fnum-p} results in additional terms of order
$v^3$ not containing any contribution of the water height $h$. Thus, these terms
are of arbitrary size and sign and destroy the positivity preservation as in
\begin{lemma}
\label{lem:fnum-LLF-type-positivity}
  The one-parameter local Lax-Friedrichs type flux \eqref{eq:fnum-LLF-type} is entropy stable,
  if $\lambda \overline{\partial_u w}$ is positive semi-definite. Additionally,
  it is positivity preserving under the CFL condition
  \eqref{eq:fnum-LLF-type-positivity-CFL}, if
  \begin{equation}
    \overline{\partial_u w}
    =
    \partial_u w \left( \mean{h}, \mean{v} \right),
    \qquad
    \lambda \geq \max \set{ \abs{v_-}, \abs{v_+} },
  \end{equation}
  and the bottom topography $b$ is continuous across cell boundaries.
\end{lemma}
In the implementation, $\lambda = \max\set{ \abs{v_-} + \sqrt{g h_-},
\abs{v_+} + \sqrt{g h_+}}$ is chosen, as in the classical local Lax-Friedrichs flux.

However, if the bottom topography $b$ is discontinuous across cell boundaries,
additional terms have to be considered, since the entropy variables are
$w = \begin{pmatrix} g (h+b) - \frac{1}{2} v^2, v \end{pmatrix}^T$. These terms are
\begin{equation}
  \frac{\Delta t}{\Delta x} \frac{\lambda_{i+\frac{1}{2}}}{2} \left(
      b_{i+1} - b_{i} \right)
  + \frac{\Delta t}{\Delta x} \frac{\lambda_{i-\frac{1}{2}}}{2} \left(
    - b_{i} + b_{i-1} \right).
\end{equation}
Adding these terms to the right hand side of equation \eqref{eq:fnum-LLF-type-positivity},
the CFL condition \eqref{eq:fnum-LLF-type-positivity-CFL} gets lost. If all
heights are positive, it is possible to guarantee $h_i^+ \geq 0$, but the
corresponding CFL condition blows up as $h \to 0$, since $b_{i\pm1} - b_{i}$
may be of arbitrary size and has to be balanced by positive contributions of
the $h$ terms.

With the choice of $\overline{\partial_u w}$ as in Lemma \ref{lem:fnum-LLF-type-positivity},
the dissipation term becomes
\begin{align*}
\stepcounter{equation}\tag{\theequation}
  &
  - \frac{\lambda}{2} \partial_u w \left( \mean{u} \right) \jump{w}
  \\
  =&
  - \frac{\lambda}{2}
  \begin{pmatrix}
    \frac{1}{g}         &  \frac{\mean{v}}{g}
    \\
    \frac{\mean{v}}{g}  &  \mean{h} + \frac{\mean{v}^2}{g}
  \end{pmatrix}
  \begin{pmatrix}
    g \jump{h+b} - \frac{1}{2} \jump{v^2}
    \\
    \jump{v}
  \end{pmatrix}
  \\
  =&
  - \frac{\lambda}{2}
  \begin{pmatrix}
    \jump{h+b} - \frac{\mean{v}}{g} \jump{v} + \frac{\mean{v}}{g} \jump{v}
    \\
    \mean{v} \jump{h+b} -\frac{\mean{v}^2}{g} \jump{v}
    + \mean{h} \jump{v} + \frac{\mean{v}^2}{g} \jump{v}
  \end{pmatrix}
  \\
  =&
  - \frac{\lambda}{2}
  \begin{pmatrix}
    \jump{h+b}
    \\
    \jump{h v} + \mean{v} \jump{b}
  \end{pmatrix},
\end{align*}
where the product rule \eqref{eq:product-2} has been used. Thus, if $b$ is continuous
across cell boundaries, $\jump{b} = 0$ and the dissipation term is simply the
classical local Lax-Friedrichs dissipation term $- \frac{\lambda}{2} \jump{u}$.

\begin{remark}
  The same numerical flux can also be obtained as a Rusanov type flux. Choosing
  a dissipation approximately as
  $- \abs{ f'(u) } \jump{u} \approx - \abs{ \partial_u f } \, \partial_w u \jump{w}$
  and using the scaling of \citet[Theorem 4]{barth1999numerical} for the
  eigenvectors, resulting in $\abs{ \partial_u f } = R \abs{ \Lambda } R^{-1}$
  and $\partial_w u = R R^T$, where $\Lambda$ is the diagonal matrix of
  eigenvalues of $f'(u)$, yields
  $- \abs{ f'(u) } \jump{u} \approx - R \abs{\Lambda} R^T \jump{w}$.
  Thus, the Rusanov choice of dissipation with $\abs{\Lambda} = \lambda \I$, where
  $\lambda > 0$ is the largest eigenvalue, yields exactly the same Lax Friedrichs
  type dissipative flux \eqref{eq:fnum-LLF-type}.
\end{remark}

\begin{remark}
  A Roe type dissipation operator can be constructed by choosing $\abs \Lambda =
  \diag{\abs\lambda_i}$, where $\lambda_i$ are the eigenvalues of $f'(u)$.
  However, some attempts to prove the preservation of positivity have not been
  successful.
\end{remark}

\subsection{Classical numerical fluxes for constant bottom topography \texorpdfstring{$b$}{b}}

Not only the local Lax Friedrichs type numerical flux \eqref{eq:fnum-LLF-type}
in the previous section is positivity preserving and entropy stable, but also
its classical variant
\begin{equation}
\label{eq:fnum-LLF}
  \fnum = \mean{f} - \frac{\lambda}{2} \jump{u}.
\end{equation}
This can be established using general results of \citet{frid2001maps, frid2004correction,
bouchut2003entropy}, but also by direct calculation.

Using the FV update procedure \eqref{eq:FV-mean-evolution}, the water height after
one time step using the local Lax Friedrichs flux \eqref{eq:fnum-LLF} becomes
\begin{align*}
\stepcounter{equation}\tag{\theequation}
  h_i^+
  =&
  h_i - \frac{\Delta t}{\Delta x} \left(
    \fnum(u_{i}, u_{i+1}) - \fnum(u_{i-1}, u_{i})
  \right)
  \\
  =&
  h_i - \frac{\Delta t}{\Delta x} \Bigg(
      \frac{h_{i+1} v_{i+1} + h_{i} v_{i}}{2}
    - \frac{\lambda_{i+\frac{1}{2}}}{2} (h_{i+1} - h_{i} )
    \widesmall{}{\\&\qquad\qquad\quad}
    - \frac{h_{i} v_{i} + h_{i-1} v_{i-1}}{2}
    + \frac{\lambda_{i-\frac{1}{2}}}{2} (h_{i} - h_{i-1} )
  \Bigg)
  \\
  =&
  h_i \left[ 1
    - \frac{\Delta t}{\Delta x} \frac{\lambda_{i-\frac{1}{2}}+\lambda_{i+\frac{1}{2}}}{2} \right]
  + h_{i+1} \frac{ \lambda_{i+\frac{1}{2}} - v_{i+1} }{2}
  + h_{i-1} \frac{ \lambda_{i-\frac{1}{2}} - v_{i-1} }{2}.
\end{align*}
Thus, if $\lambda_{i\pm\frac{1}{2}} > v_{i\pm1}$ and
$\Delta t < \frac{2 \Delta x}{\lambda_{i-\frac{1}{2}}+\lambda_{i+\frac{1}{2}}}$,
the new water height $h_i^+$ is non-negative, if the previous water heights
$h_{i-1}, h_i, h_{i+1}$ are non-negative.

The entropy stability in the semidiscrete setting with vanishing bottom topography
$b$ can be established by
\begin{align*}
\stepcounter{equation}\tag{\theequation}
  &
  \jump{w} \cdot \fnum - \jump{\psi}
  \\
  =&
  \left( g \jump{h} - \frac{1}{2} \jump{v^2} \right)
    \left( \mean{hv} - \frac{\lambda}{2} \jump{h} \right)
  \widesmall{}{\\&\quad}
  + \jump{v} \left( \mean{h v^2} + \frac{1}{2} g \mean{h^2} - \frac{\lambda}{2} \jump{hv} \right)
  - \frac{1}{2} g \jump{h^2 v}
  \\
  =&
  g \mean{h v} \jump{h}
  - \mean{h v} \mean{v} \jump{v}
  - \frac{\lambda}{2} g \jump{h}^2
  + \frac{\lambda}{2} \mean{v} \jump{h} \jump{v}
  + \mean{h v^2} \jump{v}
  \\&
  + \frac{1}{2} g \mean{h^2} \jump{v}
  - \frac{\lambda}{2} \mean{h} \jump{v}^2
  - \frac{\lambda}{2} \mean{v} \jump{h} \jump{v}
  - \frac{1}{2} g \jump{h^2 v}
  \\
  =&
  - \frac{\lambda}{2} g \jump{h}^2
  - \frac{\lambda}{2} \mean{h} \jump{v}^2
  + g \mean{h v} \jump{h}
  - \mean{h v} \mean{v} \jump{v}
  + \mean{h v^2} \jump{v}
  \widesmall{}{\\&}
  - g \mean{h} \mean{v} \jump{h}
  ,
\end{align*}
where the product rule \eqref{eq:product-2} has been used. Direct calculation yields
\begin{align*}
\stepcounter{equation}\tag{\theequation}
  &
  \jump{w} \cdot \fnum - \jump{\psi}
  \\
  =&
  - \frac{1}{2} g \underbrace{ (h_+ - h_-)^2 }_{\geq 0}
    \underbrace{ \left( \lambda - \frac{v_+ - v_-}{2} \right) }_{\geq 0}
  - \frac{1}{4} \underbrace{ h_+ }_{\geq 0} \underbrace{ (\lambda - v_+) }_{\geq 0}
    \underbrace{ (v_+ - v_-)^2 }_{\geq 0}
  \widesmall{}{\\&}
  - \frac{1}{4} \underbrace{ h_- }_{\geq 0} \underbrace{ (\lambda - v_-) }_{\geq 0}
    \underbrace{ (v_+ - v_-)^2 }_{\geq 0}
  \\
  \leq& 0,
\end{align*}
if $h_+, h_- \geq 0$ and $\lambda \geq \abs{v_+}, \abs{v_-}$. This proves
\begin{lemma}
\label{lem:LLF}
  The local Lax-Friedrichs flux \eqref{eq:fnum-LLF} used in a simple FV method
  for the shallow water equations with constant bottom topography $b$
  is positivity preserving, if the water height is non-negative and the CFL
  condition
  \begin{equation}
  \label{eq:fnum-LLF-CFL}
    1 - \frac{\Delta t}{\Delta x} \frac{\lambda_{i-\frac{1}{2}}+\lambda_{i+\frac{1}{2}}}{2}
    \geq 0
  \end{equation}
  is fulfilled. Additionally, it is entropy stable in the semidiscrete sense, if
  the water heights are non-negative $(h_+, h_- \geq 0)$ and
  $\lambda \geq \abs{v_+}, \abs{v_-}$.
\end{lemma}
In the implementation, $\lambda = \max\set{ \abs{v_-} + \sqrt{g h_-},
\abs{v_+} + \sqrt{g h_+}}$ is chosen.

\begin{remark}
  Using a Godunov scheme with exact solution of the Riemann problem results in
  an entropy stable and positivity preserving scheme, since the exact solution
  has these properties. However, some nonlinear root finding algorithm has to be
  applied to compute the exact solution of the Riemann problem. Therefore, this
  will not be considered here in more detail. The Riemann problem for the shallow
  water equations with vanishing bottom topography is described in detail inter alia
  by \citet[Chapter 5]{holden2002front}.
\end{remark}

\begin{remark}
  \citet{harten1983upstream} proposed an approximate Riemann solver using only
  one intermediate state, known as HLL Riemann solver, using estimates
  of the slowest and fastest wave speeds $s_-, s_+$. If these estimates are lower
  and upper bounds of the wave speeds in the solution of the Riemann problem,
  this flux is entropy stable, as remarked by \citet{harten1983upstream}.
  Additionally, by the right choice of wave speed estimates as in the famous
  HLLE version for Euler equations proposed by \citet{einfeldt1988godunov}, the
  numerical flux is positivity preserving, similar to the results for gas dynamics
  established by \citet{einfeldt1991godunov}, since the intermediate state in the
  approximating solution of the Riemann problem satisfies this condition.
\end{remark}

\subsection{Suliciu relaxation solver for constant bottom topography \texorpdfstring{$b$}{b}}

The Suliciu relaxation solver described by \citet[Section 2.4]{bouchut2004nonlinear}
for the shallow water equations has been implemented with some technical modifications
to allow vanishing water height $h$.
This numerical flux is entropy stable and positivity preserving under a CFL
condition not blowing up as $h \to 0$.

\subsection{Kinetic solver for constant bottom topography \texorpdfstring{$b$}{b}}

Using a kinetic approach, \citet{perthame2001kinetic} proposed a finite volume
method for the shallow water equations with general bottom topography $b$ that
has the three desired properties, i.e. that is entropy stable, positivity preserving,
and well-balanced under a suitable CFL condition not blowing up as $h \to 0$. 
However, the corresponding numerical flux has to be evaluated by some quadrature
if the bottom topography varies. This will not be used here. However, in the case
of a constant topography, all integrals can be evaluated analytically.

\subsection{Hydrostatic reconstruction approach for general bottom topography
            \texorpdfstring{$b$}{b}}
\label{sec:hydrostatic-reconstruction}

The hydrostatic reconstruction has been introduced by \citet{audusse2004fast}
as a means to extend a numerical flux for the shallow water equations with
constant bottom topography $b$ to varying $b$, preserving good properties of
the numerical flux, notably entropy stability and positivity preservation. In
addition, the resulting method is well-balanced for the lake-at-rest initial
condition.

In the context of extended numerical fluxes incorporating the source term,
the flux using hydrostatic reconstruction can be described as follows:
Compute at first the limited values
$\tilde h_i = \max\set{0, h_i+b_i-\max\set{b_i,b_k}}$,
$\tilde h_k = \max\set{0, h_k+b_k-\max\set{b_i,b_k}}$ and use the extended numerical
fluxes with new arguments
\begin{equation}
  \fnum_{h}(\tilde h_i, v_i, \tilde h_k, v_k),
  \qquad
  \fnum_{hv}(\tilde h_i, v_i, \tilde h_k, v_k) + \frac{g}{2} (h_i^2 - \tilde h_i^2),
\end{equation}
for the water height $h$ and the discharge $hv$ to compute the rate of change
in cell $i$ influenced by cell $k$.

This results in a consistent and well-balanced numerical flux that is positivity
preserving and entropy stable, if the given fluxes $\fnum_{h}, \fnum_{hv}$ have
these properties for the shallow water equations with constant bottom $b$.

However, this hydrostatic reconstruction has some disadvantages for some combinations
of bottom slope, mesh size, and water height, as described by \citet{delestre2012limitation},
at least if used for a first order FV scheme.

\begin{remark}
  By an approach based on relaxation, \citet{berthon2016fully} constructed an
  approximate Riemann solver that has the three desired properties, i.e. that is
  entropy stable, positivity preserving, and well-balanced under a suitable CFL
  condition. However, the existence of some parameters and a suitable time step
  are based on asymptotic arguments and can therefore not be implemented directly.

  However, the shallow water equations are derived based on the assumption of low
  variations in the bottom topography. Hence, a discretisation of it that is
  continuous across elements seems to be natural.
\end{remark}

\section{Finite volume subcells}
\label{sec:subcell-FV}

Although the analysis of the previous sections suggests that the semidiscretisation
of Theorem \ref{thm:gnum-extended-semidisc} with appropriate positivity preserving
and entropy stable fluxes of section \ref{sec:fluxes} and the positivity preserving
limiter of \citet{zhang2011maximum}, described also in section \ref{sec:positivity},
is stable, there are problems at wet-dry fronts in the practical implementation.
These problems can be handled by some appropriate limiting strategy, e.g.
TVB limiters used by \citet{xing2010positivity} or the slope limiter used by
\citet{duran2014recent}. However, since the high order of the approximation is
lost in these cases, the approach of finite volume subcells used similarly by
\citet{meister2016positivity} in the context of the shallow water equations will
be pursued. Applications of finite volume subcells have also been proposed inter
alia by \citet{huerta2012simple, dumbser2014posteriori, sonntag2014shock}.
Additionally, given the interpretation of SBP methods with diagonal
norm as subcell flux differencing methods by \citet{fisher2013discretely, fisher2013high},
the usage of FV subcells seems to be quite natural.

In order to use finite volume subcells to compute the time derivative, the
general procedure can be described as follows:
\begin{enumerate}
  \item 
  Decide, whether the high-order discretisation or FV subcells of first order
  should be used.
  
  \item
  Project the polynomial of degree $\leq p$ onto a piecewise constant solution.
  
  \item
  Compute the classical FV time derivative.
\end{enumerate}
As a detector to use FV subcells, the water height in the element or adjacent
elements will be used, as described in section \ref{sec:numerical-tests}.

The projection in step 2 is done for a diagonal-norm nodal SBP basis simply by
taking subcells of length $M_{i,i}$ with corresponding value $u_i$. This is not
an exact projection for the polynomial $u$ in general, but is very simple and
fits to the subcell flux differencing framework of \citet{fisher2013discretely,
fisher2013high}. It has also been used by \citet{sonntag2014shock} in the context
of the Euler equations.

Thus, for an SBP SAT semidiscretisation
\begin{equation}
  \partial_t u = - \VOL + \SURF - \mat{M}[^{-1}] \mat{R}[^T] \mat{B} \vecfnum,
\end{equation}
the surface terms $\SURF$ are set to zero and the numerical flux $\vecfnum$ is computed
using the outer values $u_0, u_p$ instead of a higher order interpolation $\mat{R} \vec{u}$
-- for nodal bases including boundary points, this makes no difference.
The volume terms $\VOL$ are computed via FV subcells as
\begin{gather}
  \left[ \VOL \right]_0
  =
  \frac{ \fnum_{0,1} }{ M_{0,0} },
  \qquad
  \left[ \VOL \right]_p
  =
  - \frac{ \fnum_{p,p-1} }{ M_{p,p} },
  \\
  \left[ \VOL \right]_i
  =
  \frac{ \fnum_{i,i+1} - \fnum_{i,i-1} }{ M_{i,i} }
  \quad \text{ for } i \in \set{1, \dots, p-1}.
\end{gather}
Therefore, the numerical flux terms $\mat{M}[^{-1}] \mat{R}[^T] \mat{B} \vecfnum$ and
the volume terms $\VOL$ form together a finite volume discretisation of the subcells.

\section{Numerical tests}
\label{sec:numerical-tests}

In this section, some numerical tests of the proposed schemes will be performed.
For the time integration, the third-order, three stage SSP method SSPRK(3,3) given
by \citet{gottlieb1998total} will be used, see also \citet{gottlieb2011strong}.

\subsection{Well-balancedness and entropy conservation}
\label{sec:NUM_TEST_WB_EC}

In this test case, the well-balancedness and entropy conservation properties of
the two-parameter family of fluxes \eqref{eq:gnum-p} and corresponding volume
\eqref{eq:gnum-extended-volume-terms-general} and surface terms
\eqref{eq:gnum-extended-surface-terms-general-ansatz}, where the parameters are
chosen according to \eqref{eq:linear-system-SURF-solution} and the free parameters
$m_{4}, k_{9}, k_{10}, k_{11}, l_{10}$ have been set to zero, are investigated.

The domain $[-1,1]$ is equipped with periodic boundary conditions and the solution
is evolved in the time interval $[0,1]$ using \num{1000} steps of the SSPRK(3,3)
method. The bottom topography and initial condition are given by
\begin{equation}
\label{eq:lake-at-rest-initial-condition}
  b(x) = \sin \frac{\pi x}{4},
  \qquad
  h_0(x) = 1 - b(x)
  ,\qquad
  hv_0(x) = 0
\end{equation}
for the lake-at-rest test case. Using $N = \frac{120}{p+1} = 15$ elements with
polynomials of degree $\leq p = 7$, represented using Gauß nodes, the simulations
have been performed for $(a_1,a_2) \in \set{-3 + \frac{k}{10},\; 0 \leq k \leq 60}^2$
with gravitational constant $g = 1$.

\begin{figure}[!phtb]
\centering
  \includegraphics[width=0.9\textwidth]{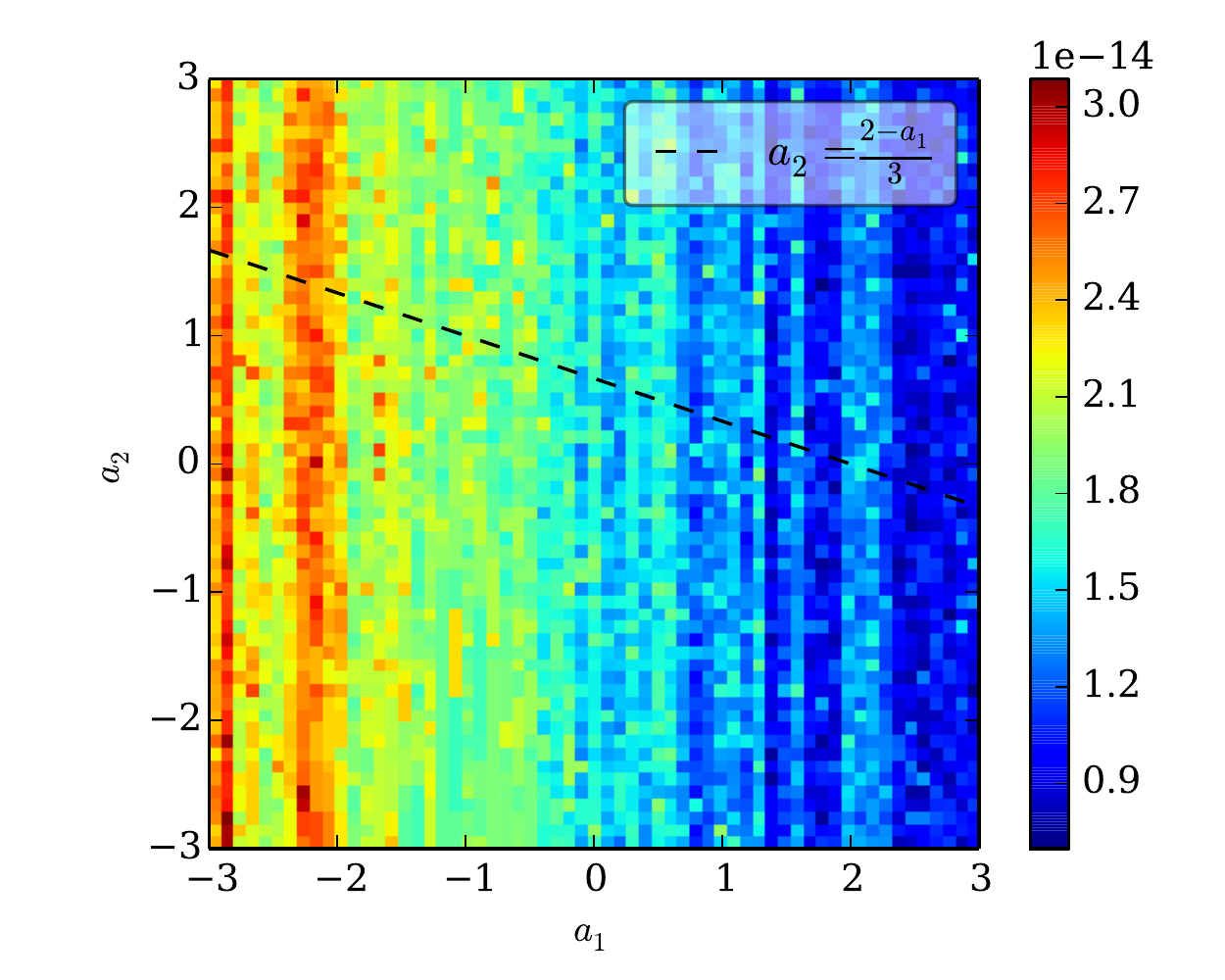}
  \caption{Maximum norm error $\max\set{\norm{h(1)-h_0}_\infty, \norm{hv(1)-hv_0}_\infty}$
           of solutions computed using the entropy conservative fluxes with varying
           parameters $a_1, a_2$ for the lake-at-rest initial condition
           \eqref{eq:lake-at-rest-initial-condition}.}          
  \label{fig:NUM_TEST_WB_EC_lake_at_rest}
\end{figure}

The results, visualised in Figure \ref{fig:NUM_TEST_WB_EC_lake_at_rest} show the
excellent well-balancedness of the methods. The maximum errors
$\max\set{\norm{h(1)-h_0}_\infty, \norm{hv(1)-hv_0}_\infty}$ (computed at the
nodes) are of order $10^{-14}$ using Float64 (i.e. double precision) in
Julia 0.4.7 \citep{bezanson2014julia}.

\begin{figure}[!phtb]
\centering
  \begin{subfigure}[b]{0.495\textwidth}
    \includegraphics[width=\textwidth]{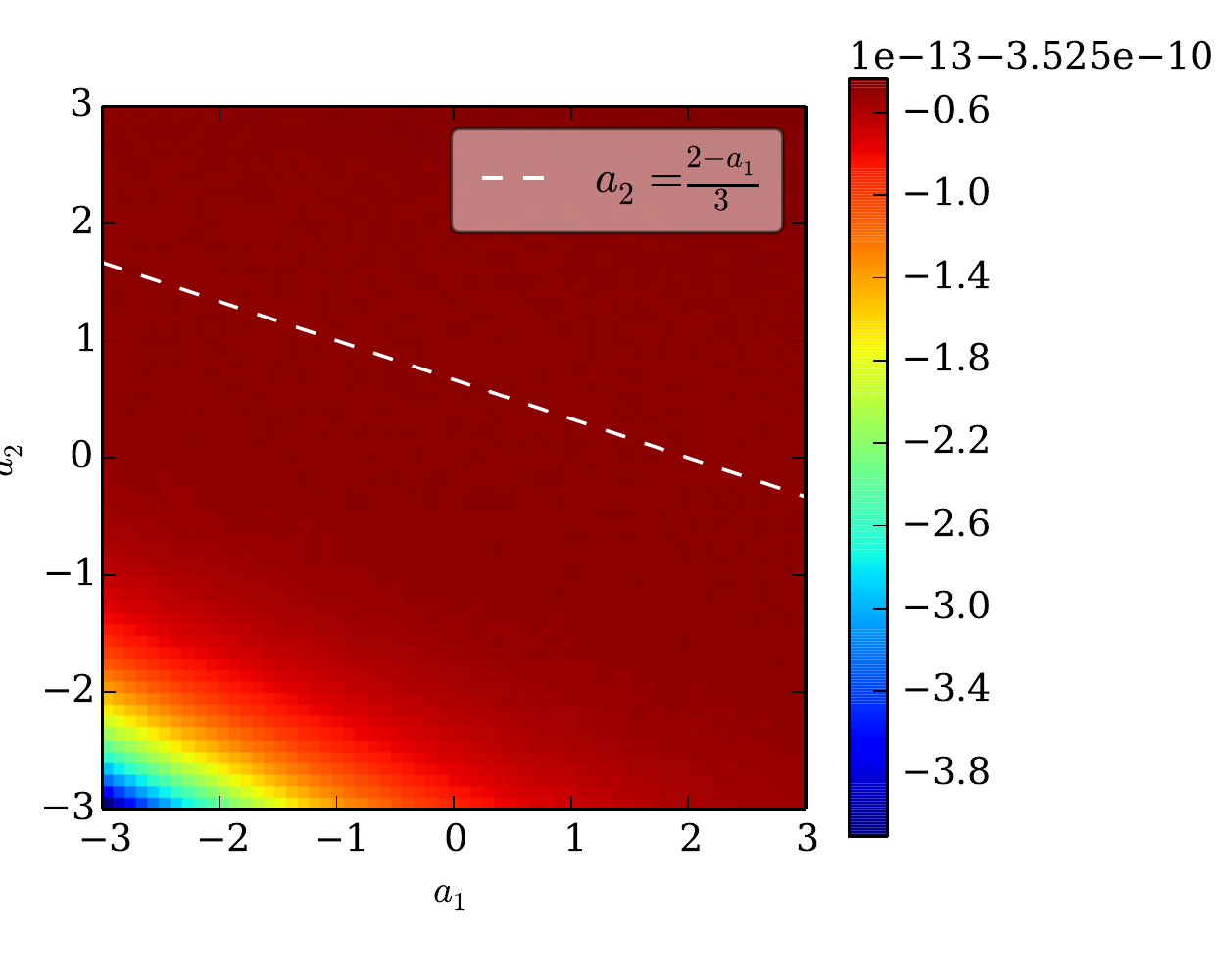}
    \caption{\num{1000} steps.}
  \end{subfigure}%
  ~
  \begin{subfigure}[b]{0.495\textwidth}
    \includegraphics[width=\textwidth]{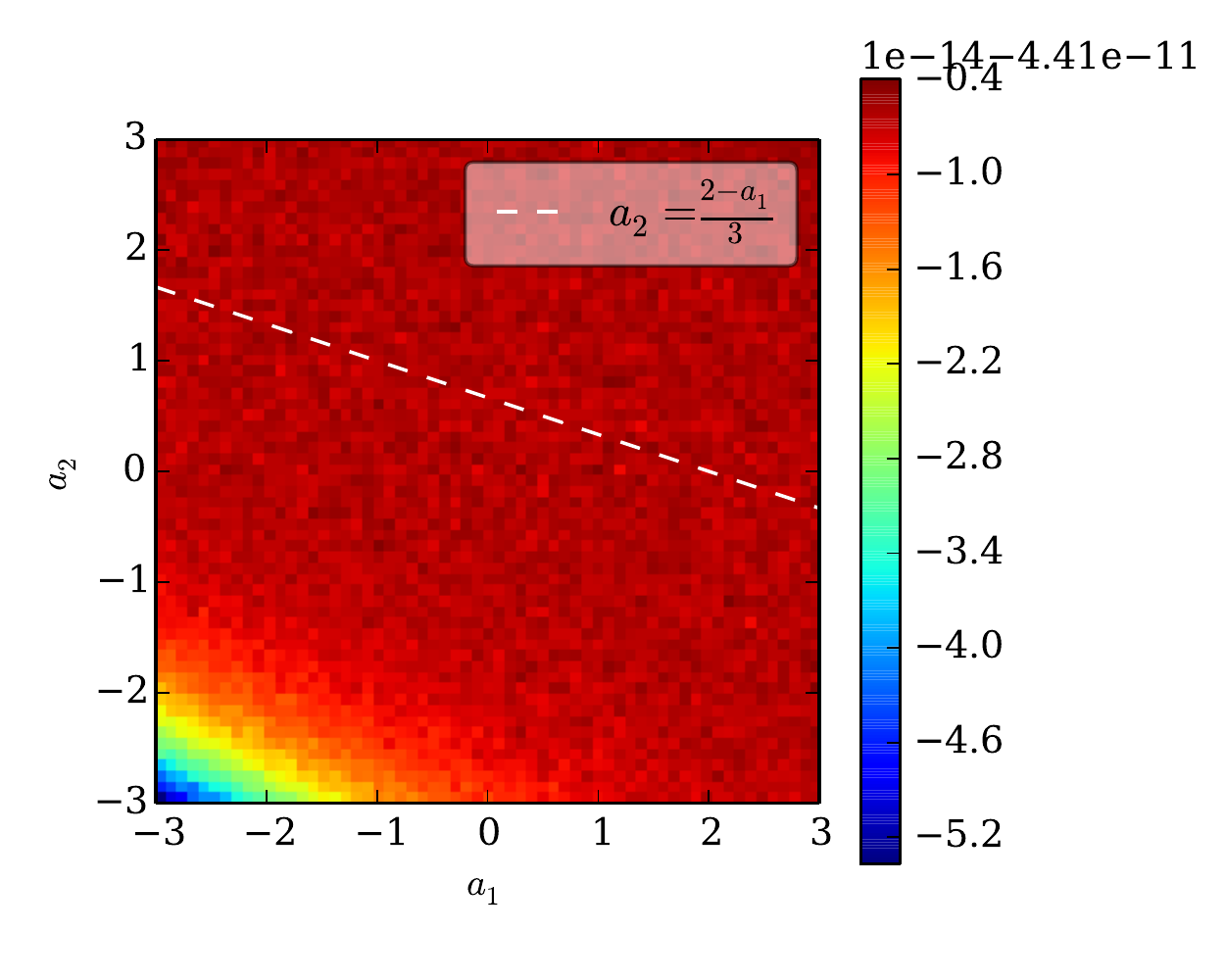}
    \caption{\num{2000} steps.}
  \end{subfigure}
  \caption{Relative entropy dissipation $\left( \int U(1) - \int U(0) \right) / \int U(0)$
           of solutions computed using the entropy conservative fluxes with varying
           parameters $a_1, a_2$ for the initial condition
           \eqref{eq:smooth-solution-initial-condition}.}          
  \label{fig:NUM_TEST_WB_EC_smooth_solution}
\end{figure}

\begin{figure}[!htb]
\centering
  \begin{subfigure}[b]{0.495\textwidth}
    \includegraphics[width=\textwidth]{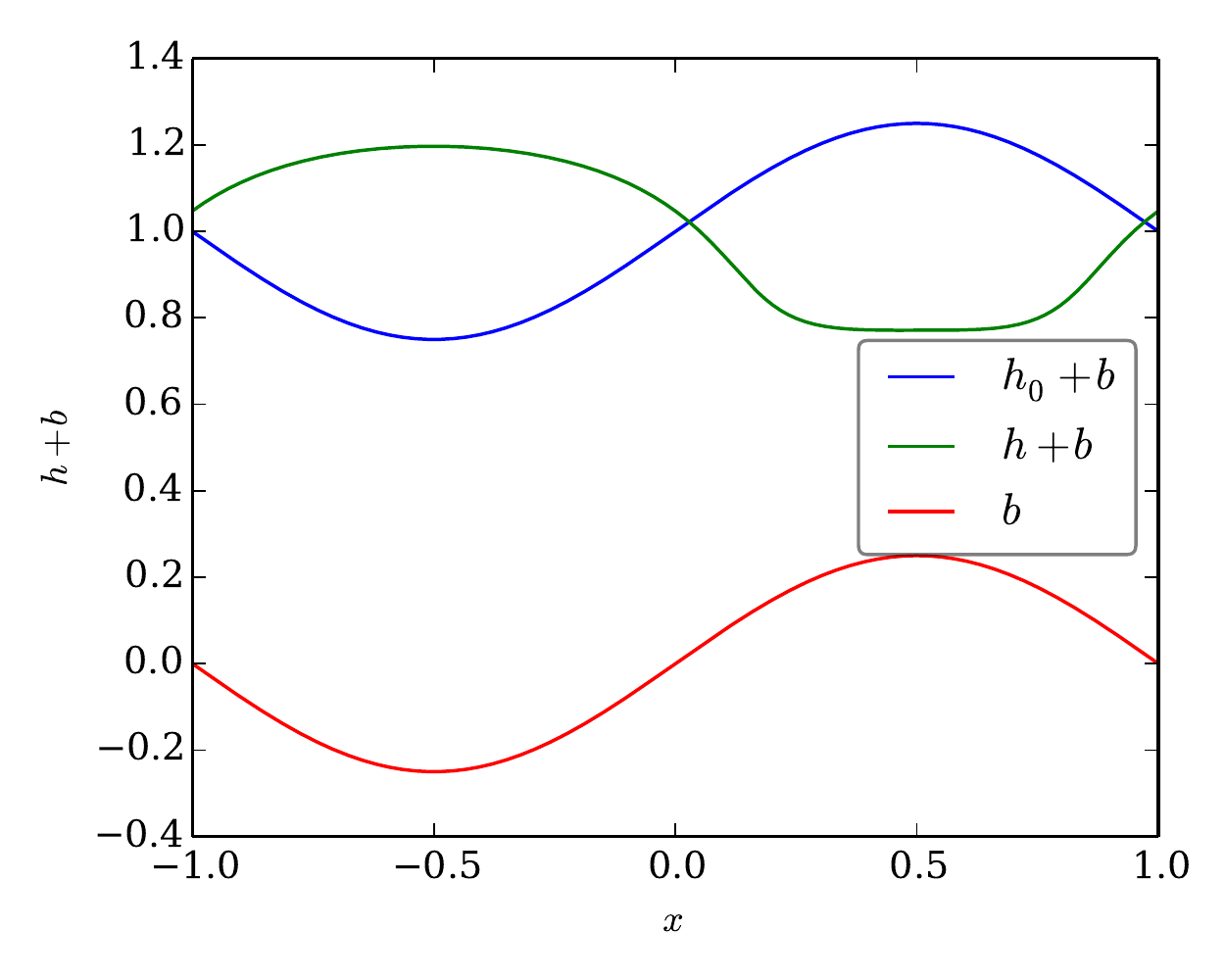}
    \caption{Water height $h(1)$.}
  \end{subfigure}%
  ~
  \begin{subfigure}[b]{0.495\textwidth}
    \includegraphics[width=\textwidth]{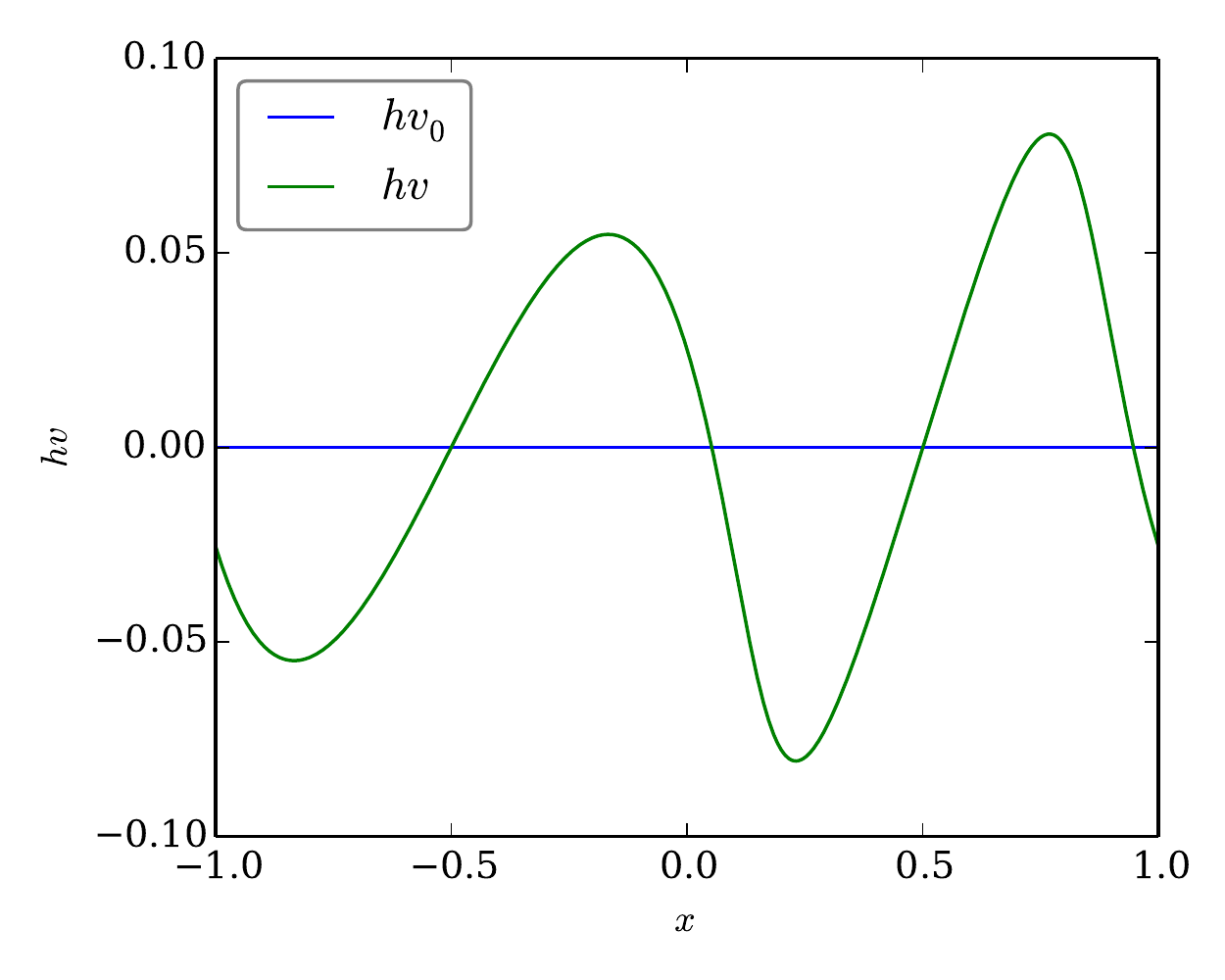}
    \caption{Discharge $hv(1)$.}
  \end{subfigure}
  \caption{Solutions computed using the entropy conservative fluxes with parameters
            $a_1 = -1, a_2 = \frac{2-a_1}{3} = 1$ for the initial condition
           \eqref{eq:smooth-solution-initial-condition}.}          
  \label{fig:NUM_TEST_WB_EC_smooth_solution__h_hv}
\end{figure}

Choosing the initial condition
\begin{equation}
\label{eq:smooth-solution-initial-condition}
  h_0(x) = 1
  ,\qquad
  hv_0(x) = 0,
\end{equation}
$h+b$ is no longer constant, but the solution remains smooth at first. Computing
again until $t = 1$ with the same parameters as before, the loss of entropy is
visualised in Figure \ref{fig:NUM_TEST_WB_EC_smooth_solution}. Using \num{1000}
steps of SSPRK(3,3), the relative entropy dissipation
$\left( \int U(1) - \int U(0) \right) / \int U(0)$ is of order $10^{-10}$, with
variations of order $10^{-13}$ for different parameters $a_1, a_2$. This loss
of entropy is caused by the time integrator, as can be seen by refining the time
step. Using \num{2000} steps, the relative dissipation is order order $10^{-11}$
with variations of order $10^{-14}$.
The smooth solutions for $a_1 = -1, a_2 = \frac{2-a_1}{3} = 1$ are plotted in
Figure \ref{fig:NUM_TEST_WB_EC_smooth_solution__h_hv}.

\begin{figure}[!htb]
\centering
  \begin{subfigure}[b]{0.495\textwidth}
    \includegraphics[width=\textwidth]{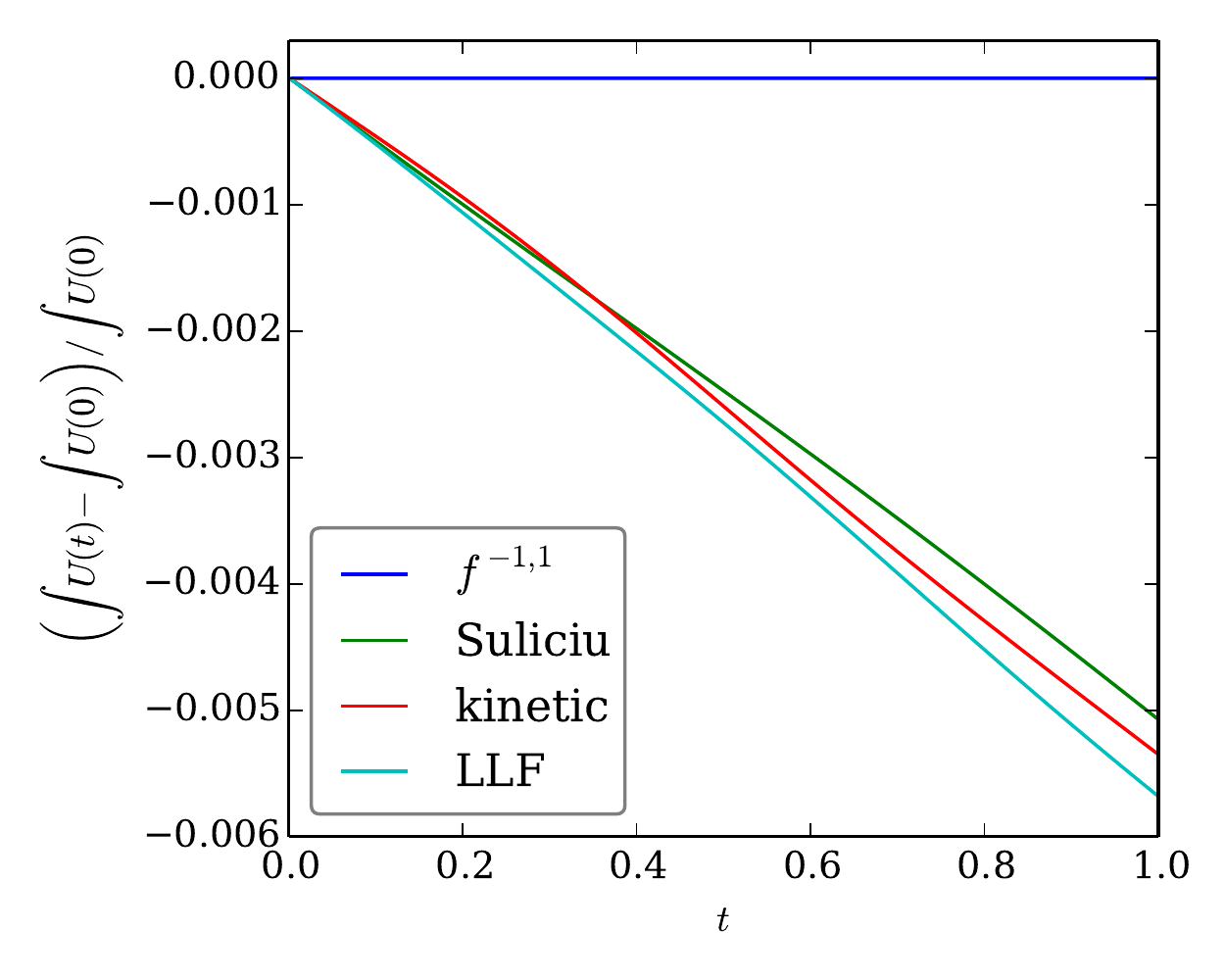}
    \caption{$p = 0$, $N = \frac{120}{p+1} = 120$.}
  \end{subfigure}%
  ~
  \begin{subfigure}[b]{0.495\textwidth}
    \includegraphics[width=\textwidth]{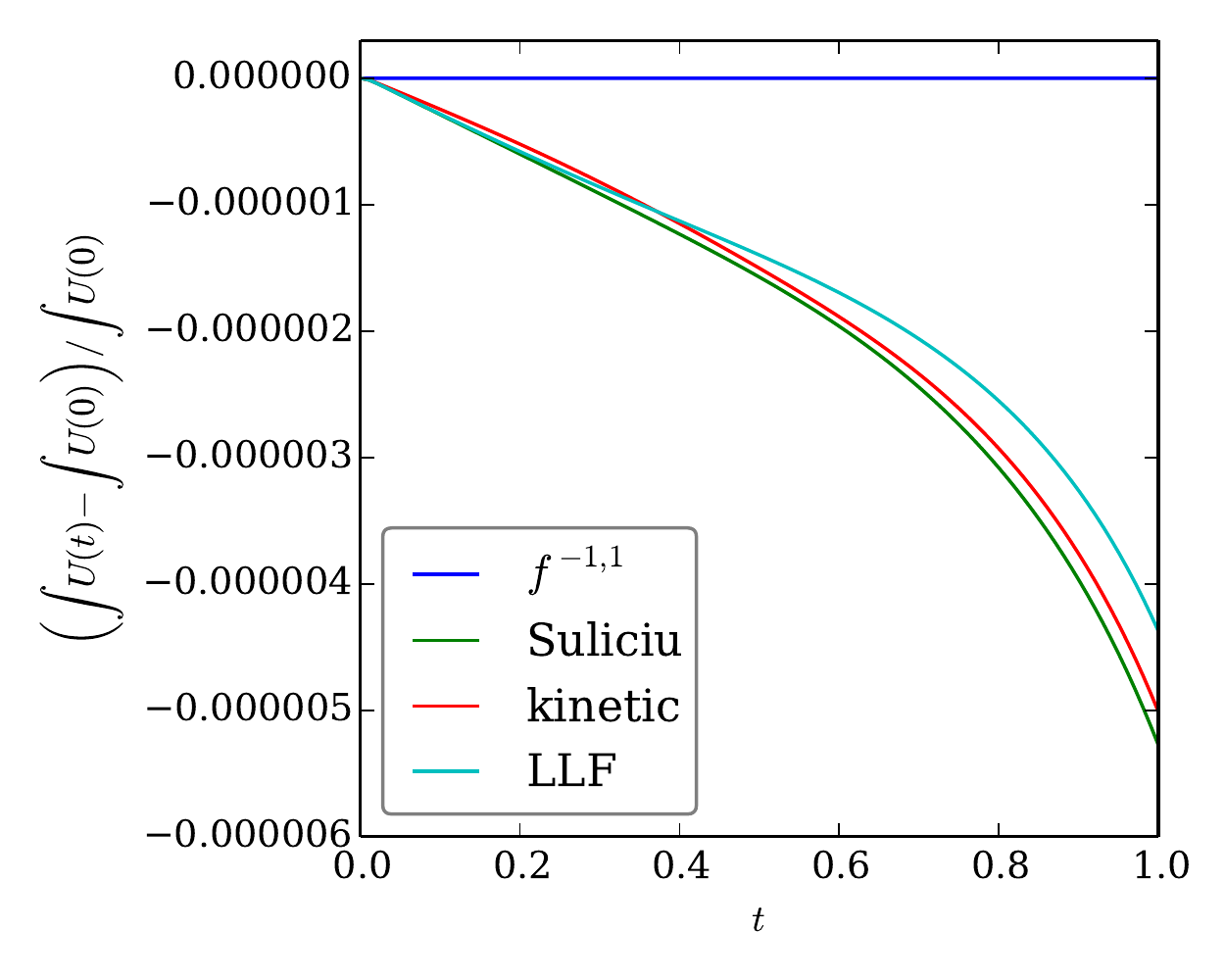}
    \caption{$p = 1$, $N = \frac{120}{p+1} = 60$.}
  \end{subfigure}%
  \\
  \begin{subfigure}[b]{0.495\textwidth}
    \includegraphics[width=\textwidth]{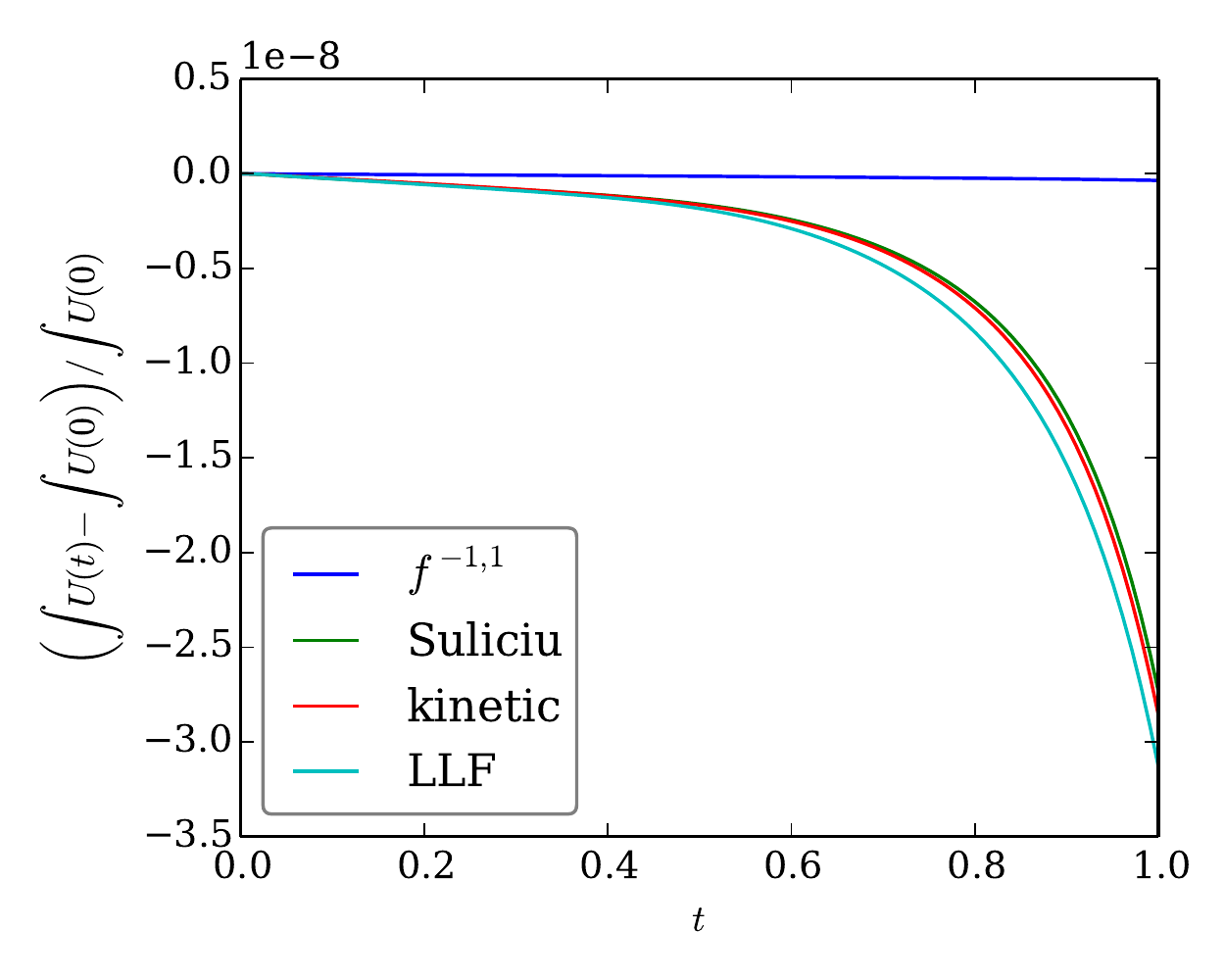}
    \caption{$p = 2$, $N = \frac{120}{p+1} = 40$.}
  \end{subfigure}%
  ~
  \begin{subfigure}[b]{0.495\textwidth}
    \includegraphics[width=\textwidth]{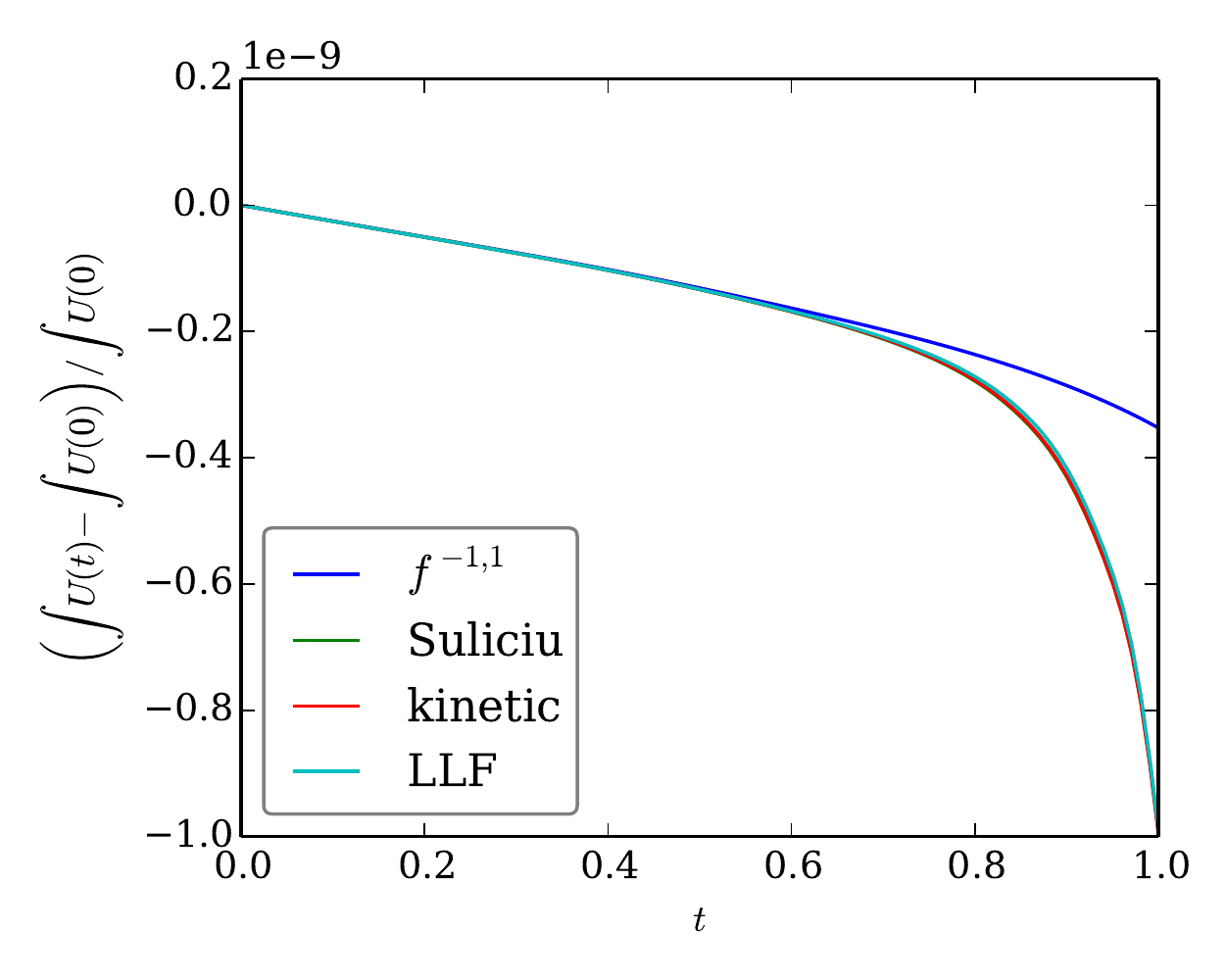}
    \caption{$p = 3$, $N = \frac{120}{p+1} = 30$.}
  \end{subfigure}%
  \\
  \begin{subfigure}[b]{0.495\textwidth}
    \includegraphics[width=\textwidth]{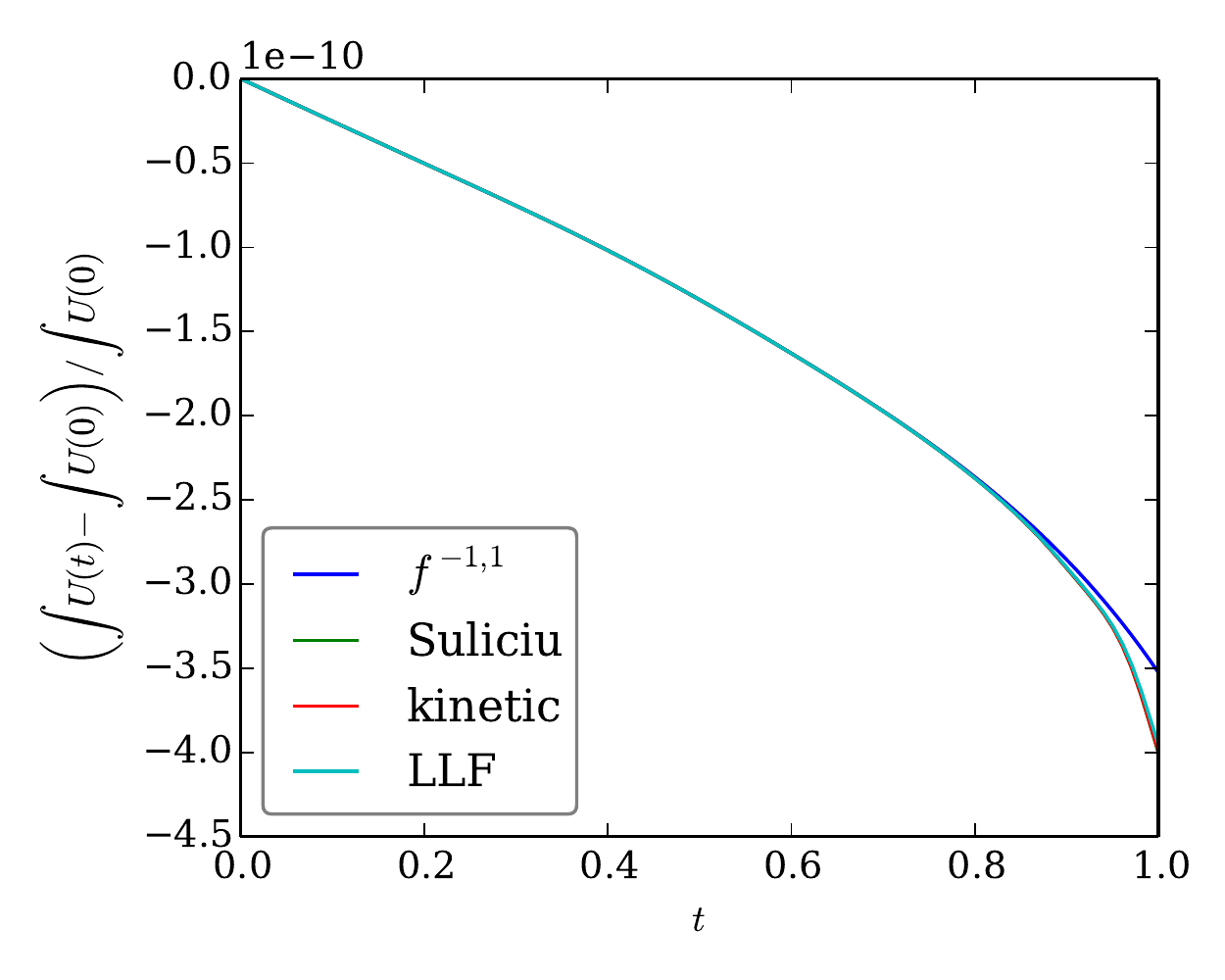}
    \caption{$p = 4$, $N = \frac{120}{p+1} = 20$.}
  \end{subfigure}%
  ~
  \begin{subfigure}[b]{0.495\textwidth}
    \includegraphics[width=\textwidth]{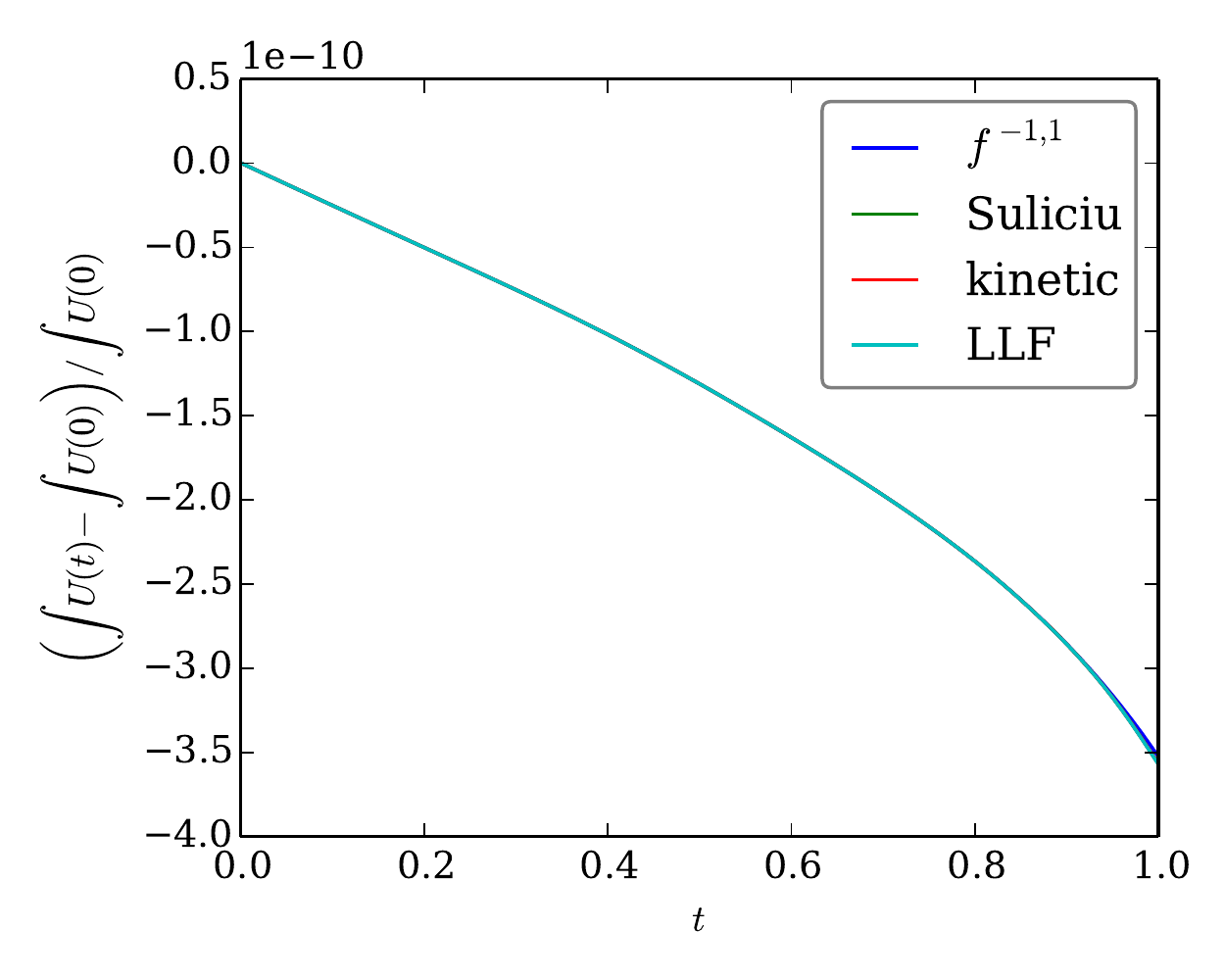}
    \caption{$p = 5$, $N = \frac{120}{p+1} = 15$.}
  \end{subfigure}%
  \caption{Relative entropy dissipation $\left( \int U(1) - \int U(0) \right) / \int U(0)$
           of solutions computed using different surface fluxes for the initial
           condition \eqref{eq:smooth-solution-initial-condition} with varying
           degree $p$ and number of elements $N = \frac{120}{p+1}$.}          
  \label{fig:NUM_TEST_WB_EC_smooth_solution_1k__U_diff_fluxes}
\end{figure}

The influence of the numerical (surface) flux is visualised in Figure
\ref{fig:NUM_TEST_WB_EC_smooth_solution_1k__U_diff_fluxes}. There, the number of
degrees of freedom $N \cdot (p+1)$ has been kept constant, while the polynomial
degree $p$ varies between $0$ (first order FV scheme) and $5$. As can be seen there,
the entropy conservative flux $f^{-1,1}$ is indeed entropy conservative, while
the entropy stable fluxes are a bit dissipative. The dissipation increases
from the Suliciu flux \citep[Section 2.4]{bouchut2004nonlinear} over the kinetic
flux \citep{perthame2001kinetic}
to the local Lax-Friedrichs flux \eqref{eq:fnum-LLF}. All three fluxes have
been implemented using the hydrostatic reconstruction of \citet{audusse2004fast}
described in section \ref{sec:hydrostatic-reconstruction}.

In the finite volume setting $p = 0$, the dissipation is of order $10^{-3}$ and
decreases with increasing polynomial degree $p$. For $p \geq 3$, the curves for
the dissipative fluxes become visually indistinguishable and for
$p \geq 5$ they coincide with the entropy conservative flux $f^{-1,1}$ for this
smooth solution.

\subsection{Lake at rest with emerged bump}
\label{sec:NUM_TEST_emerged_bump}

Here, the lake-at-rest initial condition of SWASHES \citep[Section 3.1.2]{delestre2013swashes}
\begin{equation}
\label{eq:emerged-bump-initial-condition}
\begin{aligned}
  b(x)
  =&
  \begin{cases}
    0.2 - 0.05 (x-10)^2
    ,\text{ if } 8 < x < 12,
    \\
    0,
    \text{ else,}
  \end{cases}
  \\
  h_0(x) =& \max\set{0.1, b(x)} - b(x),
  \quad
  hv_0(x) = 0,
\end{aligned}
\end{equation}
will be used in the domain $[0, 25]$ with periodic boundary conditions for simulations
in the time interval $[0, 1]$ with gravitational constant $g = 9.81$.

\begin{figure}[!htb]
\centering
  \begin{subfigure}[b]{0.495\textwidth}
    \includegraphics[width=\textwidth]{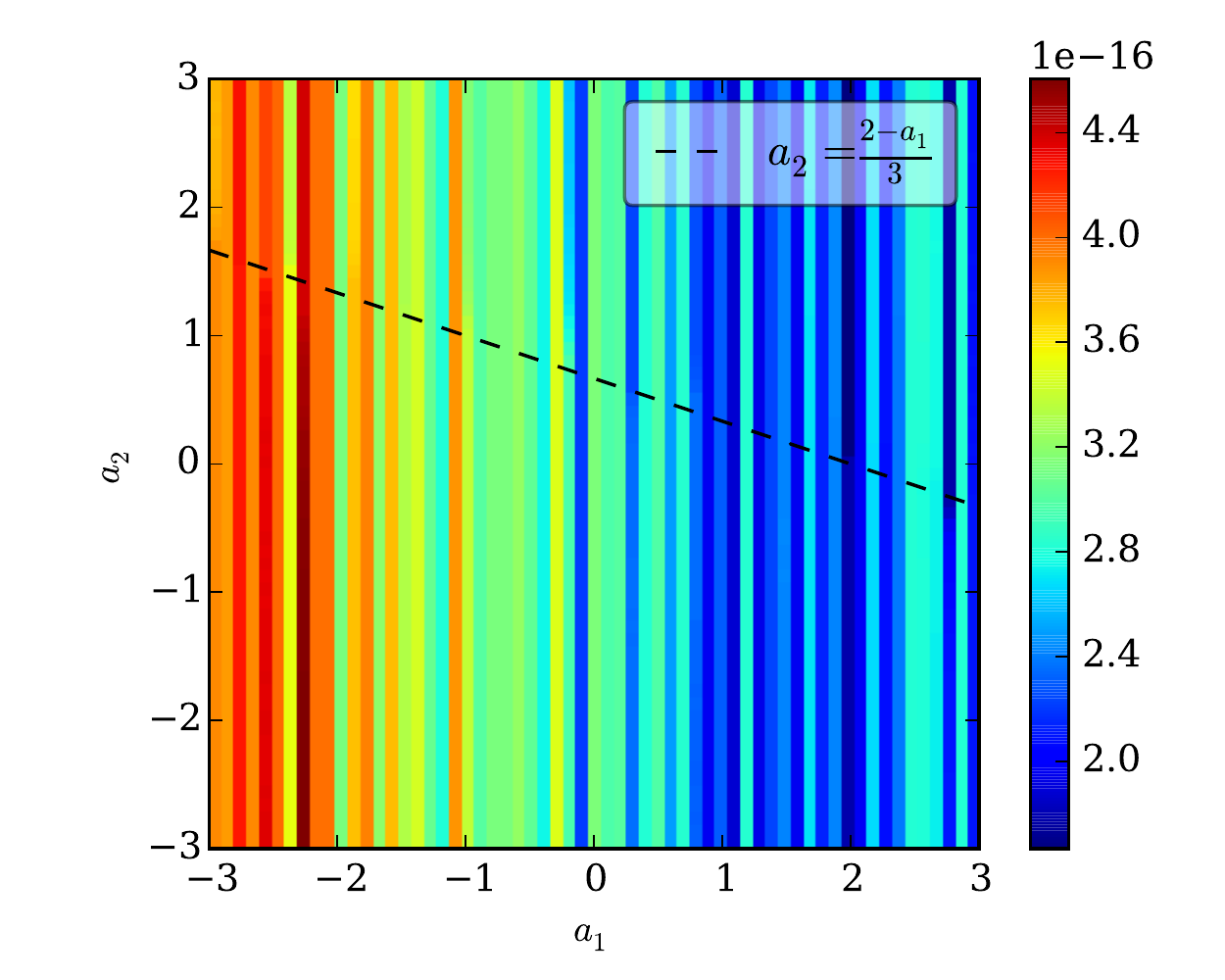}
    \caption{Maximum error.}
  \end{subfigure}%
  ~
  \begin{subfigure}[b]{0.495\textwidth}
    \includegraphics[width=\textwidth]{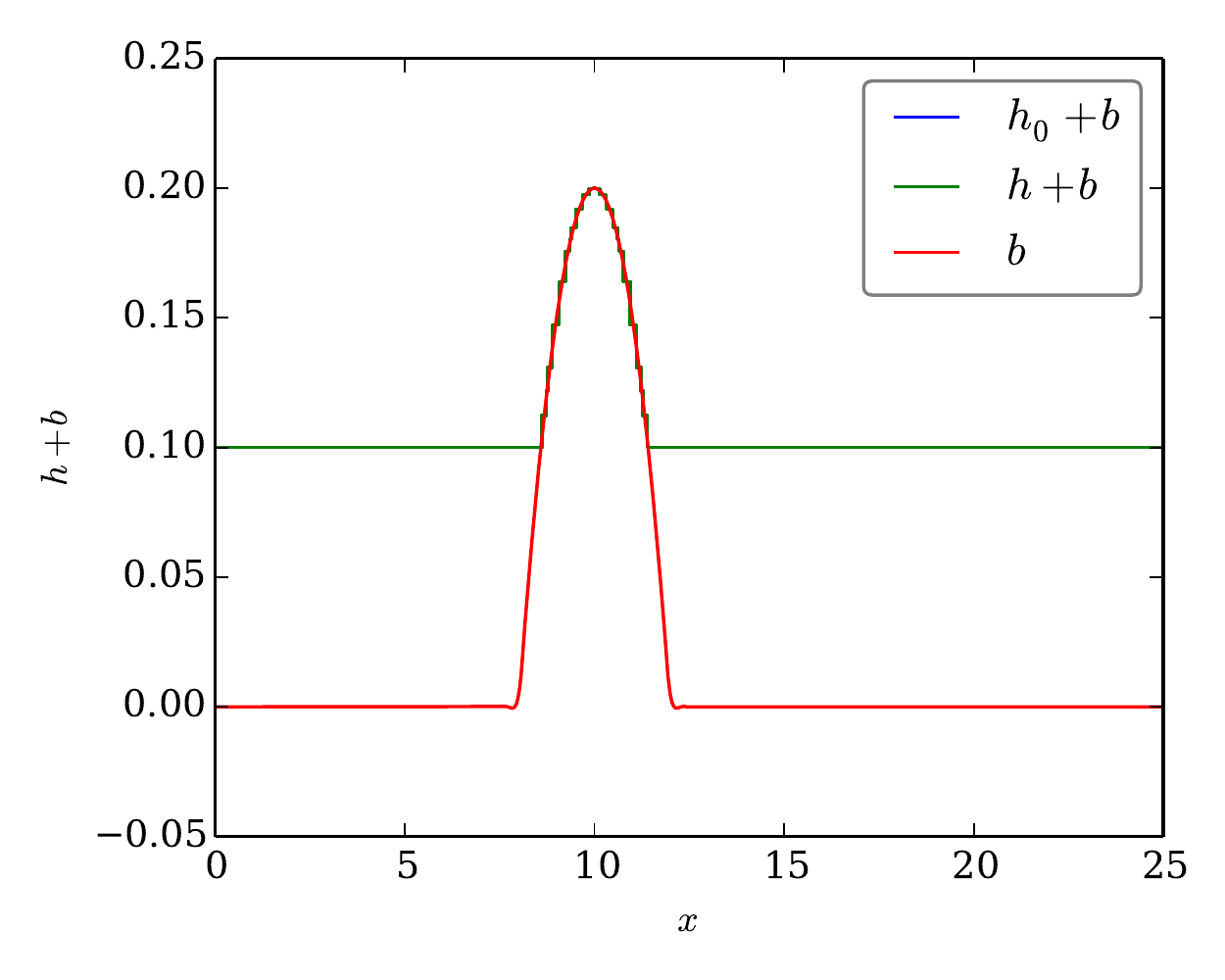}
    \caption{Water height $h(1)$ for $a_1=-1,a_2=1$.}
  \end{subfigure}
  \caption{Maximum norm error $\max\set{\norm{h(1)-h_0}_\infty, \norm{hv(1)-hv_0}_\infty}$
           of solutions computed using the entropy conservative fluxes with varying
           parameters $a_1, a_2$ for the lake-at-rest with emerged bump initial
           condition \eqref{eq:emerged-bump-initial-condition} and solution at
           $t=1$ for $a_1=-1, a_2=1$.}          
  \label{fig:NUM_TEST_emerged_bump}
\end{figure}

If no FV subcells are used, the result strongly depends on the resolution of the
shore and needs in general some additional dissipation to be stable
near the wet-dry front. However, activating FV subcells if the water height $h$
at some node in the element is smaller than $10^{-5}$, the simulation is stable.

These results are shown in Figure \ref{fig:NUM_TEST_emerged_bump} for $N = 40$
elements of polynomials of degree $\leq p = 5$ and the local Lax-Friedrichs
flux \eqref{eq:fnum-LLF} with hydrostatic reconstruction as numerical flux.

The maximum error norm $\max\set{\norm{h(1)-h_0}_\infty, \norm{hv(1)-hv_0}_\infty}$
(computed at the nodes) is of order of magnitude $10^{-16}$ for varying parameters
$(a_1,a_2) \in \set{-3 + \frac{k}{10},\; 0 \leq k \leq 60}^2$ used for the
volume terms \eqref{eq:gnum-extended-volume-terms-general}. Again, Gauß nodes
and corresponding surface terms have been used, where the additional free parameters
have been set to zero. Additionally, the water height for the choice $a_1 = -1,
a_2 = \frac{2-a_1}{3} = 1$ is visualised there.

\subsection{Moving water equilibrium with varying bottom \texorpdfstring{$b$}{b}}
\label{sec:NUM_TEST_moving_water}

Here, a moving water equilibrium of the shallow water equations with gravitational
constant $g = 9.81$ given by 
\begin{equation}
\label{eq:movin-water-eq}
  hv\equiv m = \mathrm{const},
  \qquad
  \frac{1}{2} v^2 + g(h + b) \equiv E = \mathrm{const}
\end{equation}
is considered. The bottom topography is
\begin{equation}
  b(x)
  =
  \begin{cases}
    \frac{1}{4} \cos \left( 10\pi (x+1) \right) + \frac{1}{4}
    , \text{ if } -0.1 < x < 0.1,
    \\
    0, \text{ else,}
  \end{cases}
\end{equation}
and the initial condition is computed by solving the second equation of
\eqref{eq:movin-water-eq} for $h$, inserting $v^2 = \frac{(hv)^2}{h^2} = \frac{m^2}{h^2}$.
Two initial conditions $m = 1, E = 25$ and $m = 3, E = \frac{3}{2} (m g)^{2/3}
+ \frac{g}{2} = 19.203311922761937$ are considered, similar to \citet{audusse2015simple}.

Computing the maximum error $\max\set{\norm{h(1)-h_0}_\infty, \norm{hv(1)-hv_0}_\infty}$
at the nodes yields identical results for both initial conditions with polynomial
degrees $\leq p \in \set{0,\dots,9}$, and parameters $(a_1,a_2) \in
\set{-3 + \frac{k}{10},\; 0 \leq k \leq 60}^2$ for the volume terms, while the
local Lax-Friedrichs flux \eqref{eq:fnum-LLF} has been used as numerical flux.
The domain is divided into $N = 40$ elements using Gauß nodes.

\begin{figure}[!htb]
\centering
  \begin{subfigure}[b]{0.495\textwidth}
    \includegraphics[width=\textwidth]{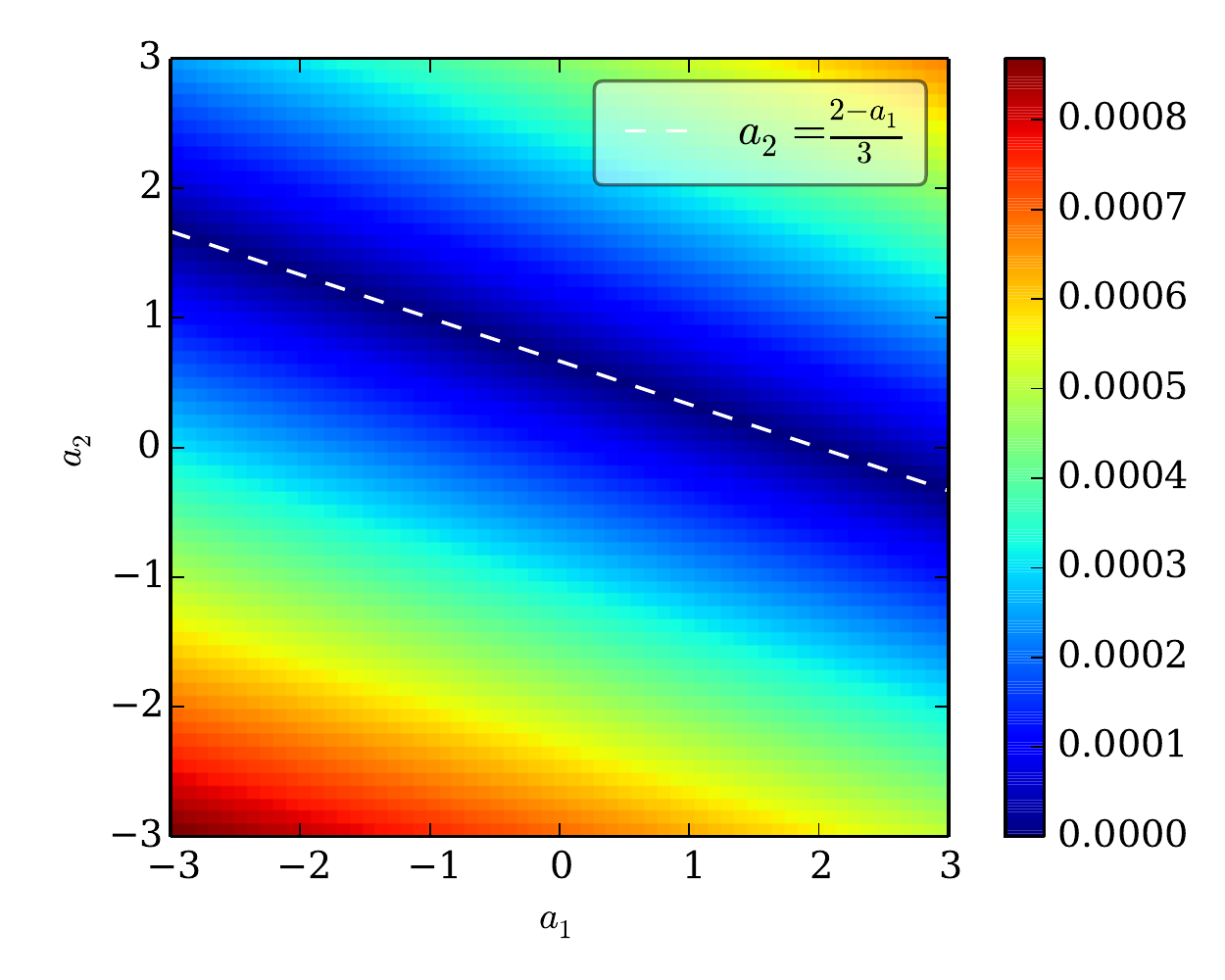}
    \caption{Maximum error for $p=9$.}
  \end{subfigure}%
  ~
  \begin{subfigure}[b]{0.495\textwidth}
    \includegraphics[width=\textwidth]{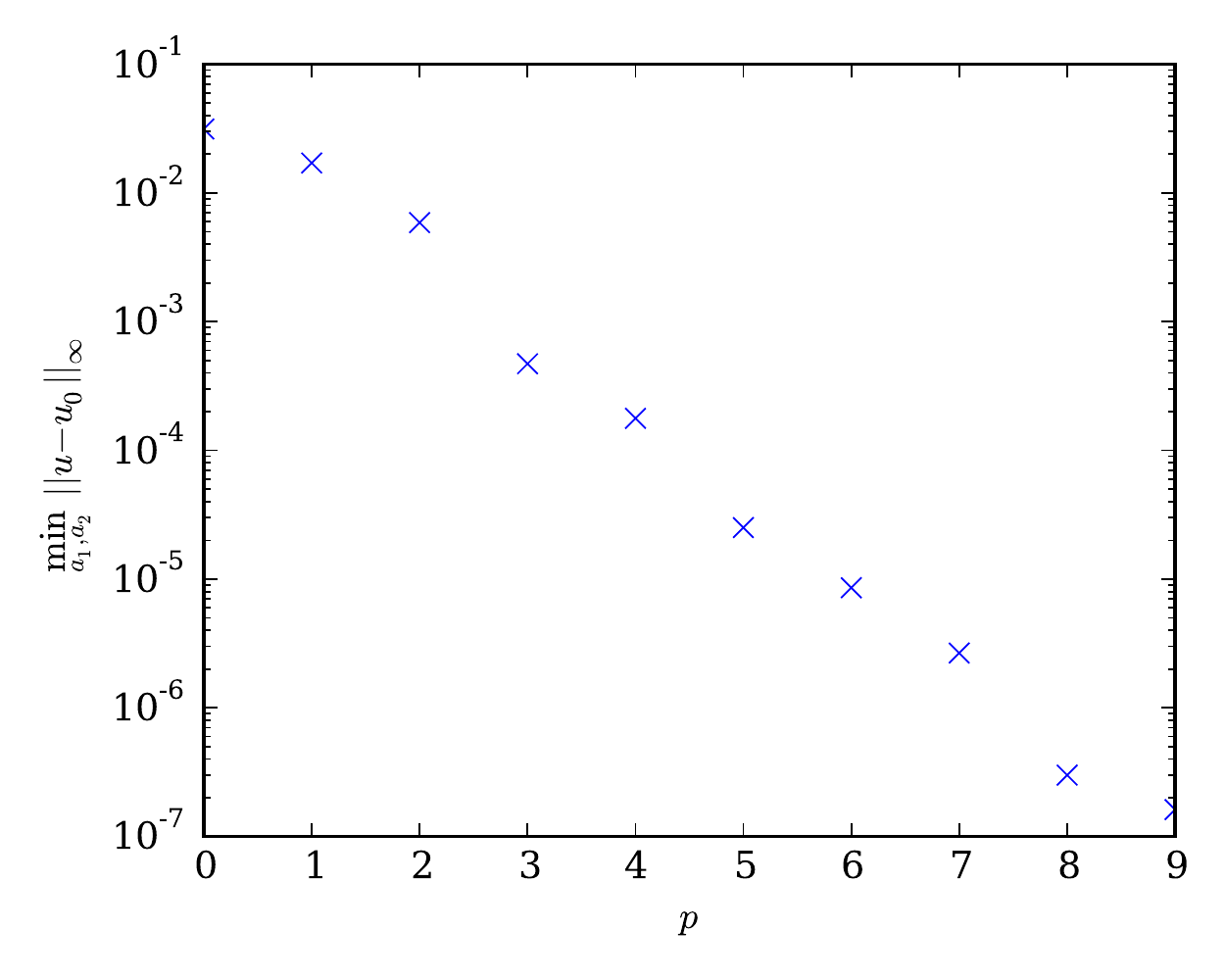}
    \caption{Minimal errors for varying $p$.}
  \end{subfigure}
  \caption{Maximum norm error $\max\set{\norm{h(1)-h_0}_\infty, \norm{hv(1)-hv_0}_\infty}$
           of solutions computed using the entropy conservative fluxes with varying
           parameters $a_1, a_2$ for moving water equilibrium \eqref{eq:movin-water-eq}
           with $m=1, E=25$ for polynomials of degree $\leq p = 9$ and minimal
           values of the maximum errors over $a_1,a_2$.}          
  \label{fig:NUM_TEST_moving_water}
\end{figure}

These results with $m=1, E=25$ are shown in Figure \ref{fig:NUM_TEST_moving_water}.
As can be seen there, the choice $a_2 = \frac{2-a_1}{3}$ is optimal for this problem,
while the choice of $a_1$ does not seem to be critical. This can be explained by
the additional term in $v$ for $a_2 \neq \frac{2-a_1}{3}$ in the numerical flux
\eqref{eq:gnum-p} and the corresponding volume terms
\eqref{eq:gnum-extended-volume-terms-general}. The results for $m = 3, E =
\frac{3}{2} (m g)^{2/3}$ are visually indistinguishable.

Additionally, the minimal values of the maximum error over the parameters
$a_1, a_2$ are plotted in Figure \ref{fig:NUM_TEST_moving_water} for $m=1, E=25$.
The usual superior properties of odd polynomial degrees $p$ as well as
exponential convergence can be seen there.

\subsection{Dam break}
\label{sec:NUM_TEST_dam_break}

\begin{figure}[!htb]
\centering
  \begin{subfigure}[b]{0.495\textwidth}
    \includegraphics[width=\textwidth]{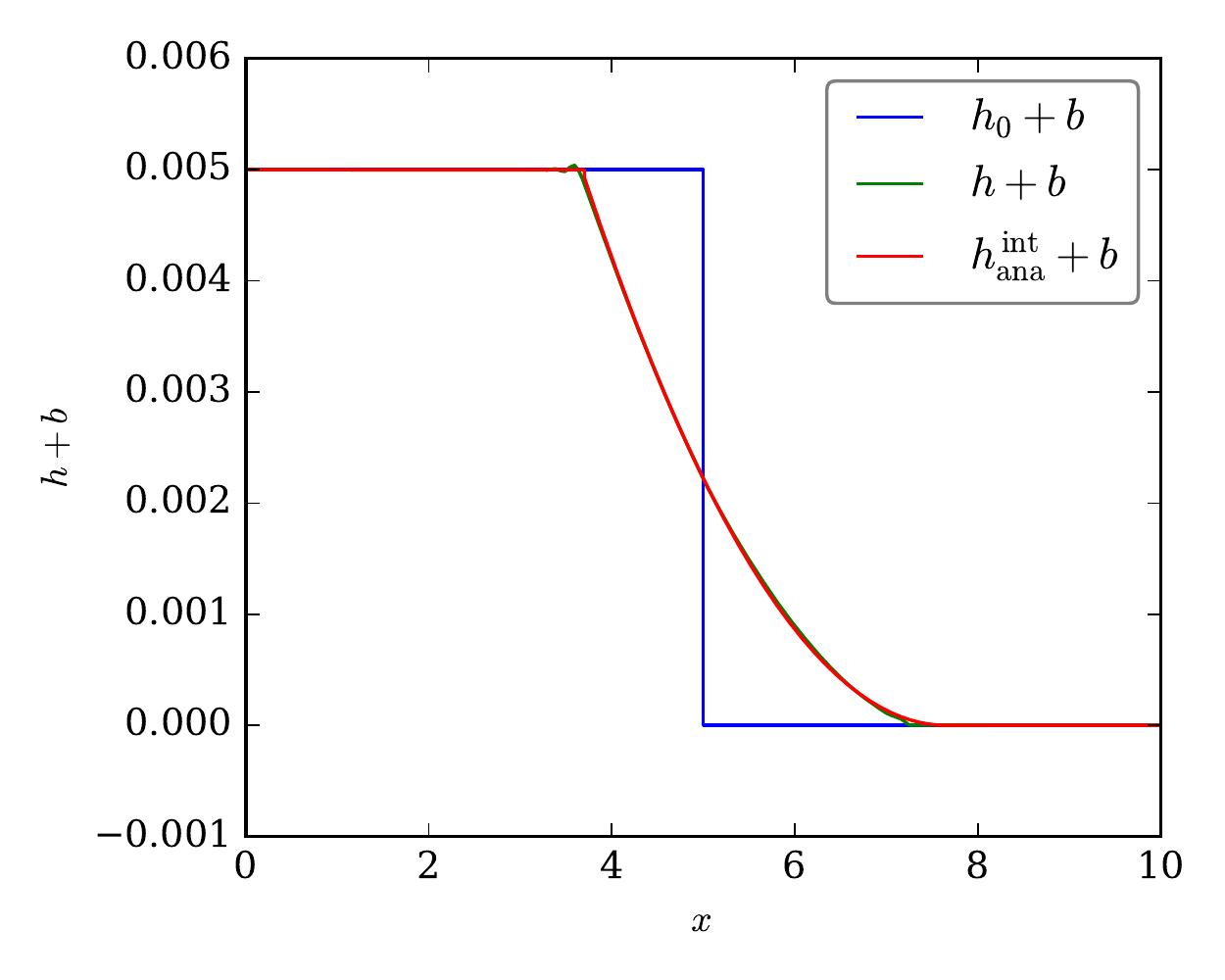}
    \caption{Water height $h$.}
  \end{subfigure}%
  ~
  \begin{subfigure}[b]{0.495\textwidth}
    \includegraphics[width=\textwidth]{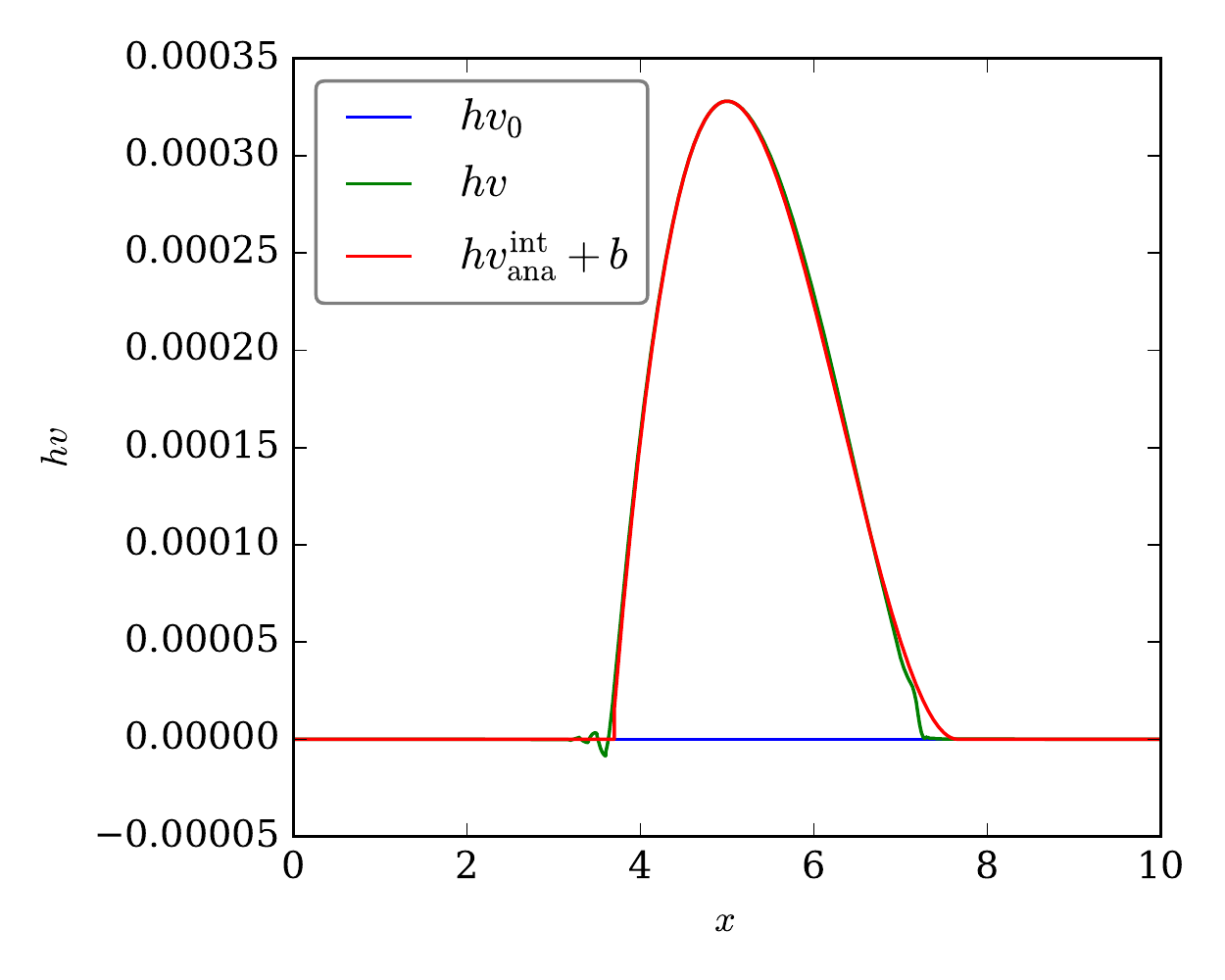}
    \caption{Discharge $hv$.}
  \end{subfigure}
  \caption{Numerical solution computed using $N = 100$ elements of polynomials
           of degree $\leq p = 2$ on Gauß nodes for the dam break problem
           \eqref{eq:dam-break} with local Lax-Friedrichs numerical flux and
           $a_1 = -1, a_2 = 1$.}          
  \label{fig:NUM_TEST_dam_breaks_solution}
\end{figure}

Here, the dam break problem with dry domain and an analytical solution
described by \citet[Section 4.1.2]{delestre2013swashes} will be considered.
The initial condition
\begin{equation}
\label{eq:dam-break}
  h_0(x)
  =
  \begin{cases}
    0.005, \text{ if } x < 5,
    \\
    0, \text{ else},
  \end{cases}
  \qquad
  hv_0(x) = 0,
  \qquad
  b(x) = 0,
\end{equation}
is evolved in the domain $[0,10]$ until $t = 6$, and the gravitational constant
is again $g = 9.81$.
The results of a simulation using $N = 100$ elements
with polynomials of degree $\leq p = 2$ and the local Lax-Friedrichs numerical
flux \eqref{eq:fnum-LLF} are plotted in Figure \ref{fig:NUM_TEST_dam_breaks_solution}.
Here, FV subcells are used in a cell if the water height in the cell itself or
adjacent cells is less then $10^{-6}$, and the parameters are chosen as
$a_1 = -1, a_2 = \frac{2-a_1}{3}= 1$.

\begin{figure}[!htb]
\centering
  \begin{subfigure}[b]{0.495\textwidth}
    \includegraphics[width=\textwidth]{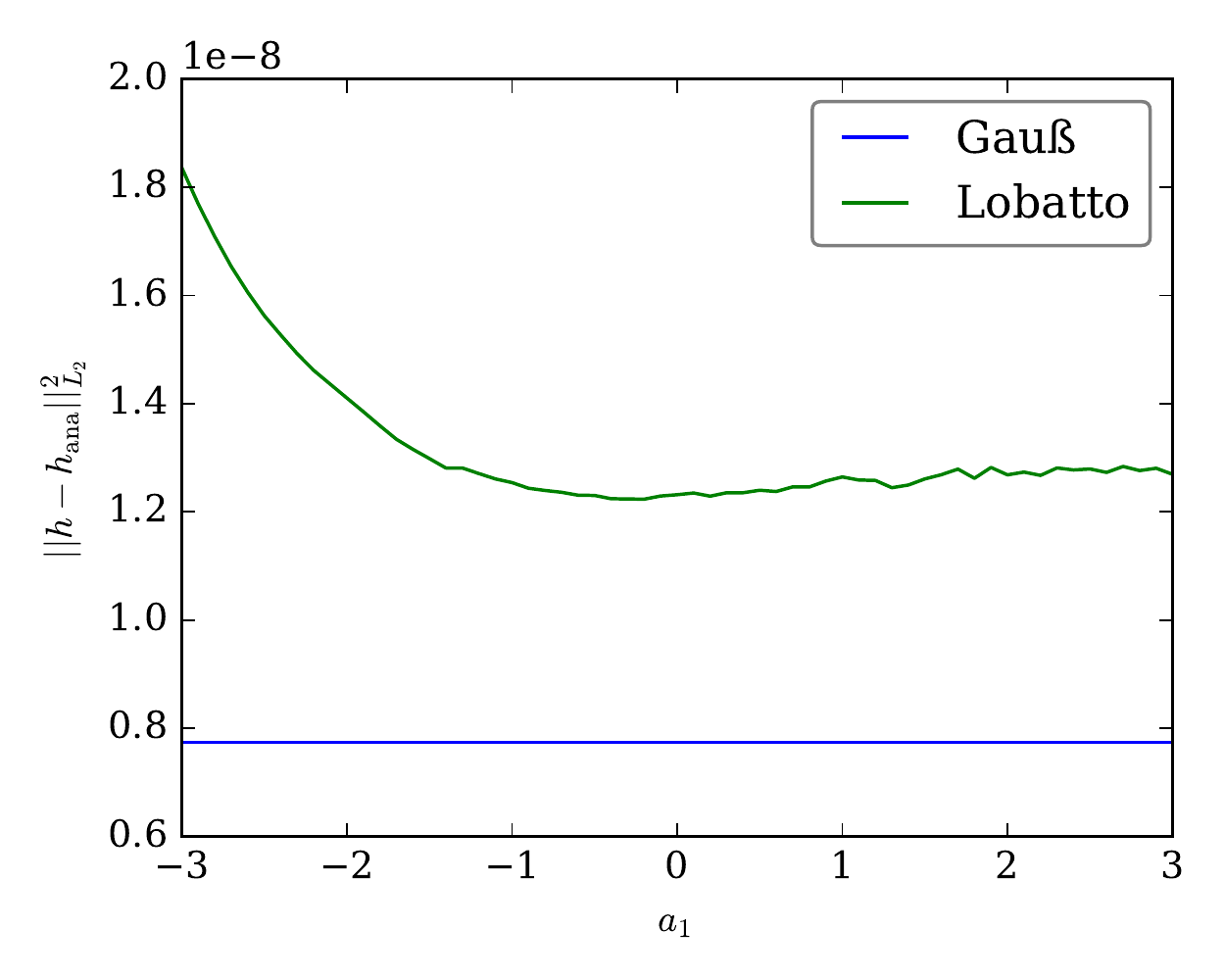}
    \caption{Error $\norm{h - h_\mathrm{ana}}_{L_2}^2$.}
  \end{subfigure}%
  ~
  \begin{subfigure}[b]{0.495\textwidth}
    \includegraphics[width=\textwidth]{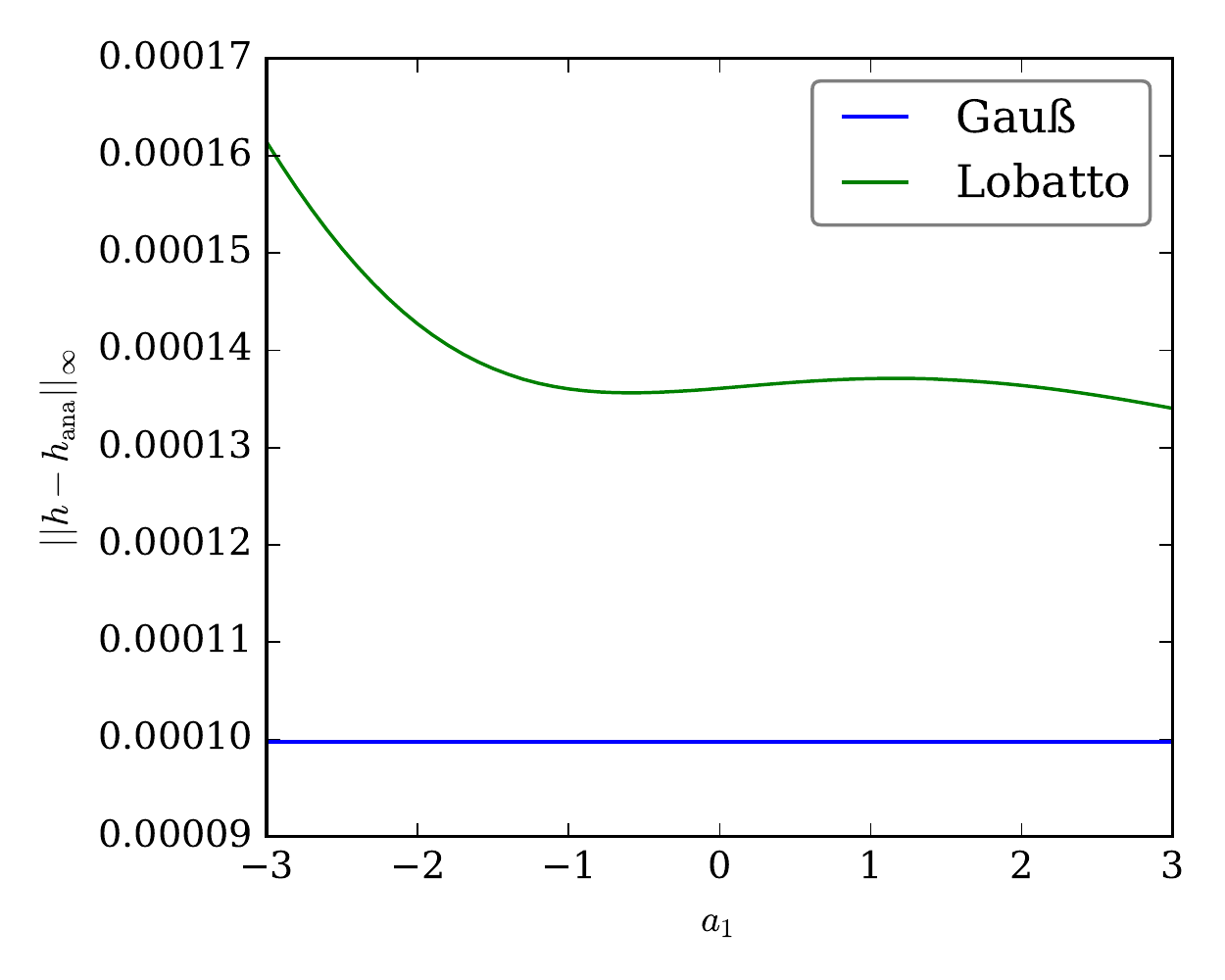}
    \caption{Error $\norm{h - h_\mathrm{ana}}_\infty$.}
  \end{subfigure}%
  \\
  \begin{subfigure}[b]{0.495\textwidth}
    \includegraphics[width=\textwidth]{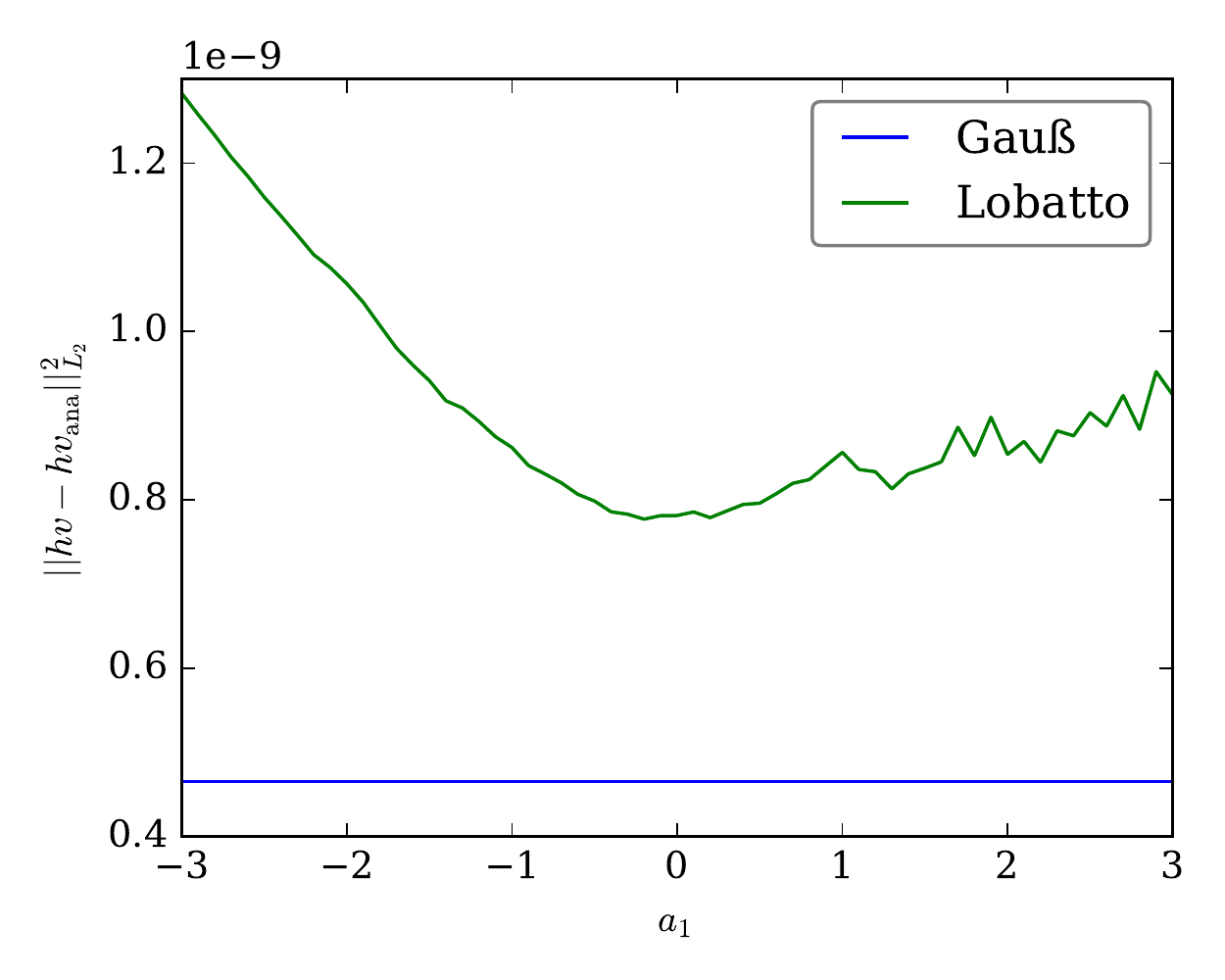}
    \caption{Error $\norm{hv - hv_\mathrm{ana}}_{L_2}^2$.}
  \end{subfigure}%
  ~
  \begin{subfigure}[b]{0.495\textwidth}
    \includegraphics[width=\textwidth]{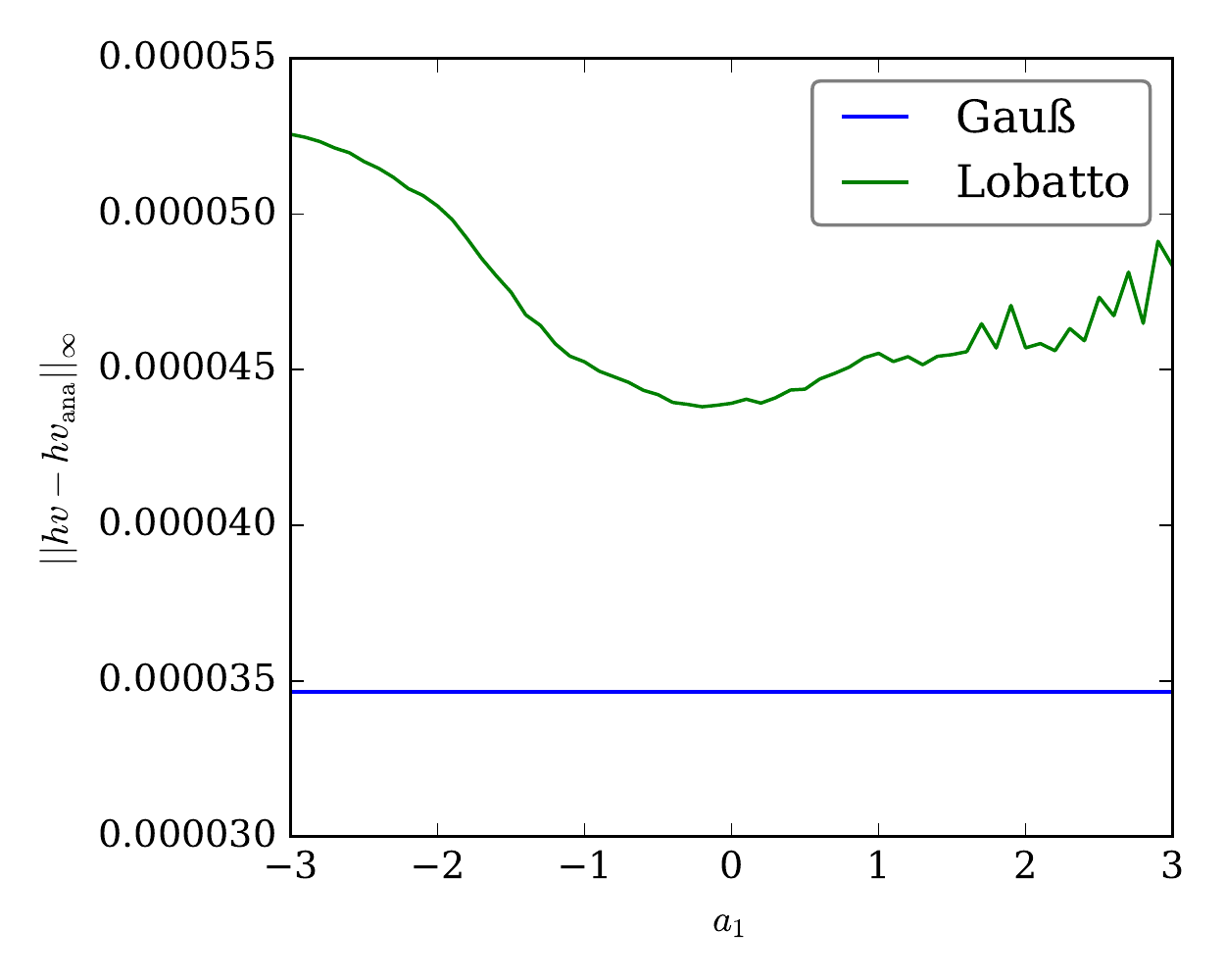}
    \caption{Error $\norm{hv - hv_\mathrm{ana}}_\infty$.}
  \end{subfigure}%
  \caption{Errors of the numerical solution for varying parameter $a_1$ and
           $N = 50$ elements of polynomials of degree $\leq p = 2$ using Gauß and
           Lobatto nodes.}          
  \label{fig:NUM_TEST_dam_breaks_p_2_N_50}
\end{figure}

\begin{figure}[!htb]
\centering
  \begin{subfigure}[b]{0.495\textwidth}
    \includegraphics[width=\textwidth]{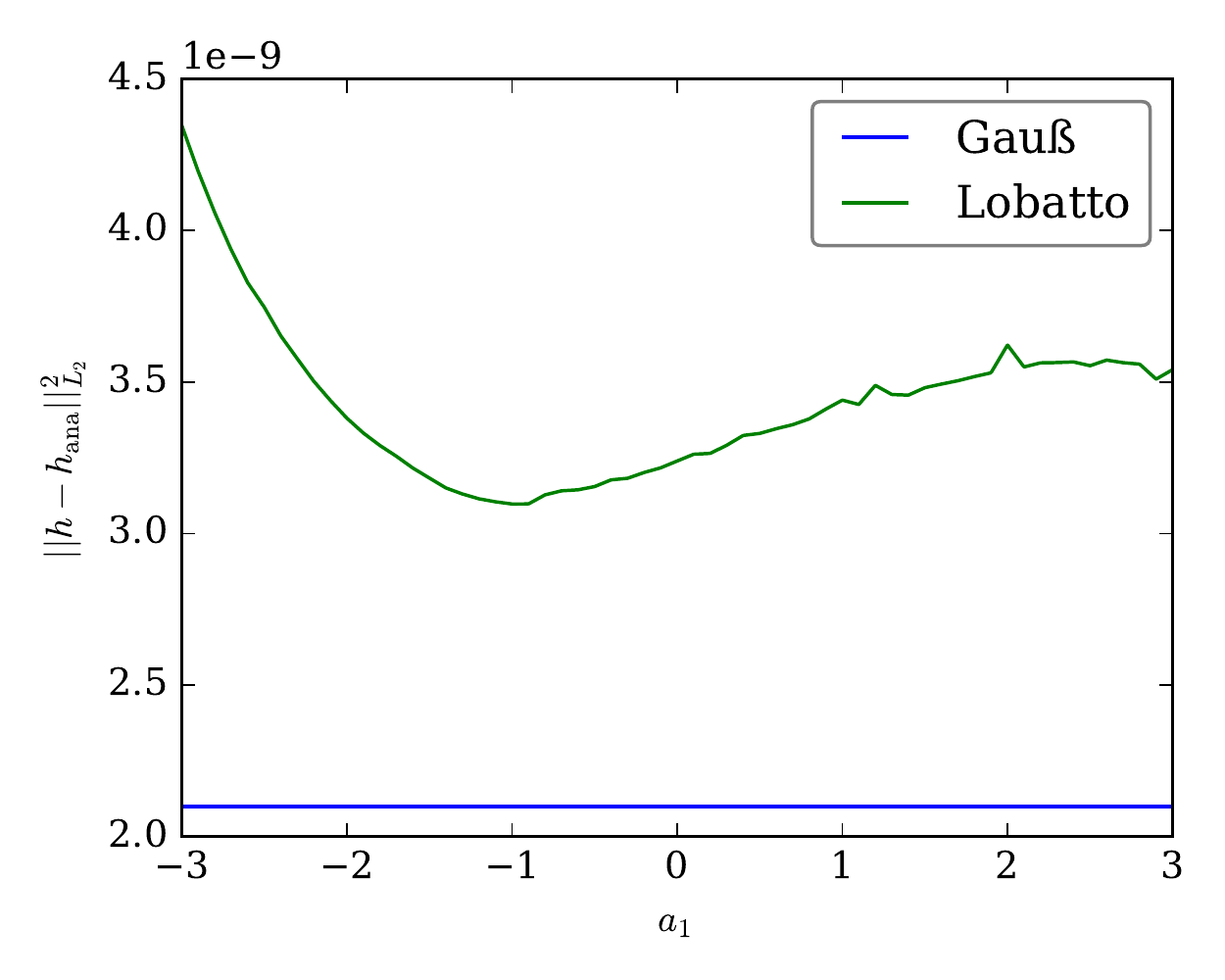}
    \caption{Error $\norm{h - h_\mathrm{ana}}_{L_2}^2$.}
  \end{subfigure}%
  ~
  \begin{subfigure}[b]{0.495\textwidth}
    \includegraphics[width=\textwidth]{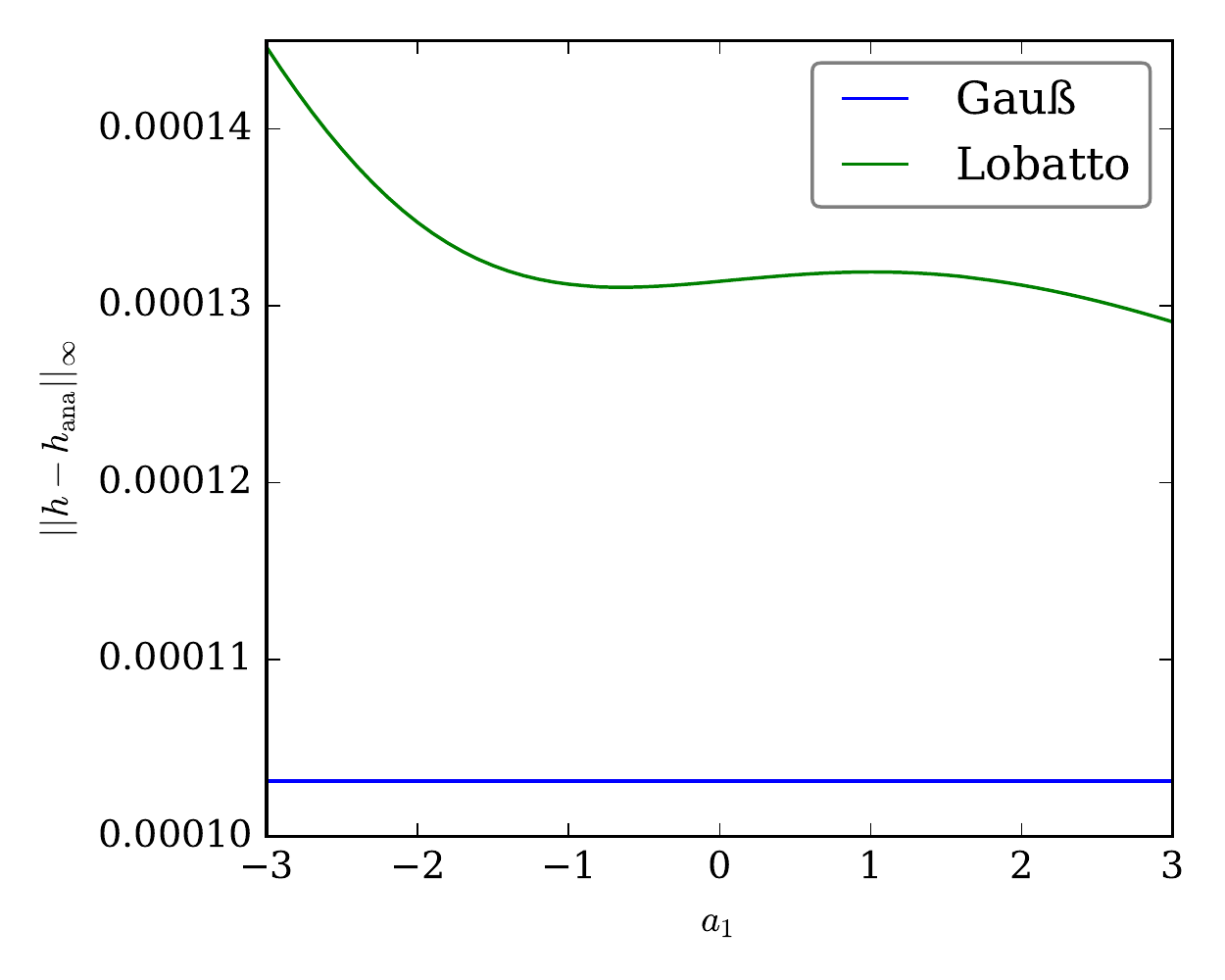}
    \caption{Error $\norm{h - h_\mathrm{ana}}_\infty$.}
  \end{subfigure}%
  \\
  \begin{subfigure}[b]{0.495\textwidth}
    \includegraphics[width=\textwidth]{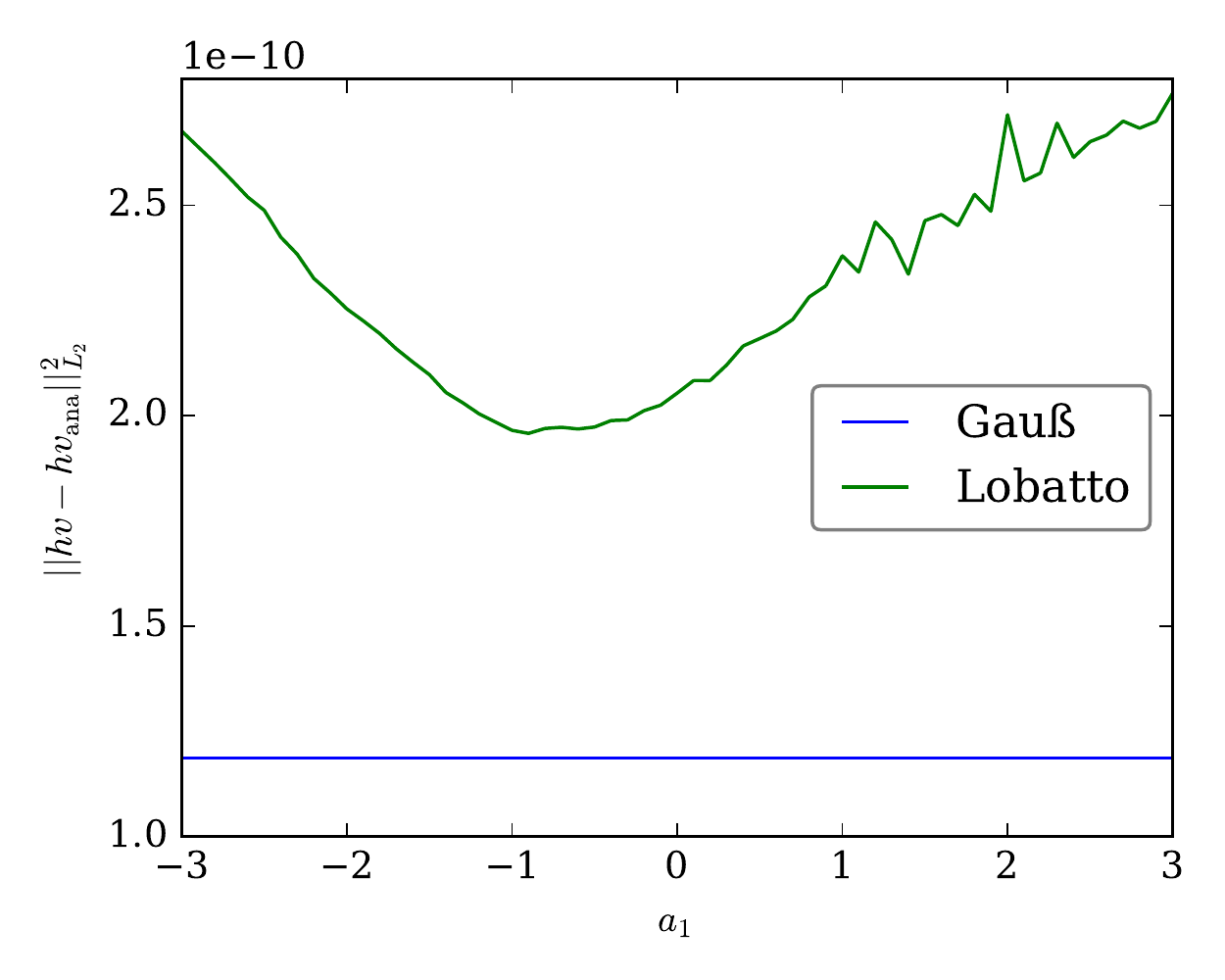}
    \caption{Error $\norm{hv - hv_\mathrm{ana}}_{L_2}^2$.}
  \end{subfigure}%
  ~
  \begin{subfigure}[b]{0.495\textwidth}
    \includegraphics[width=\textwidth]{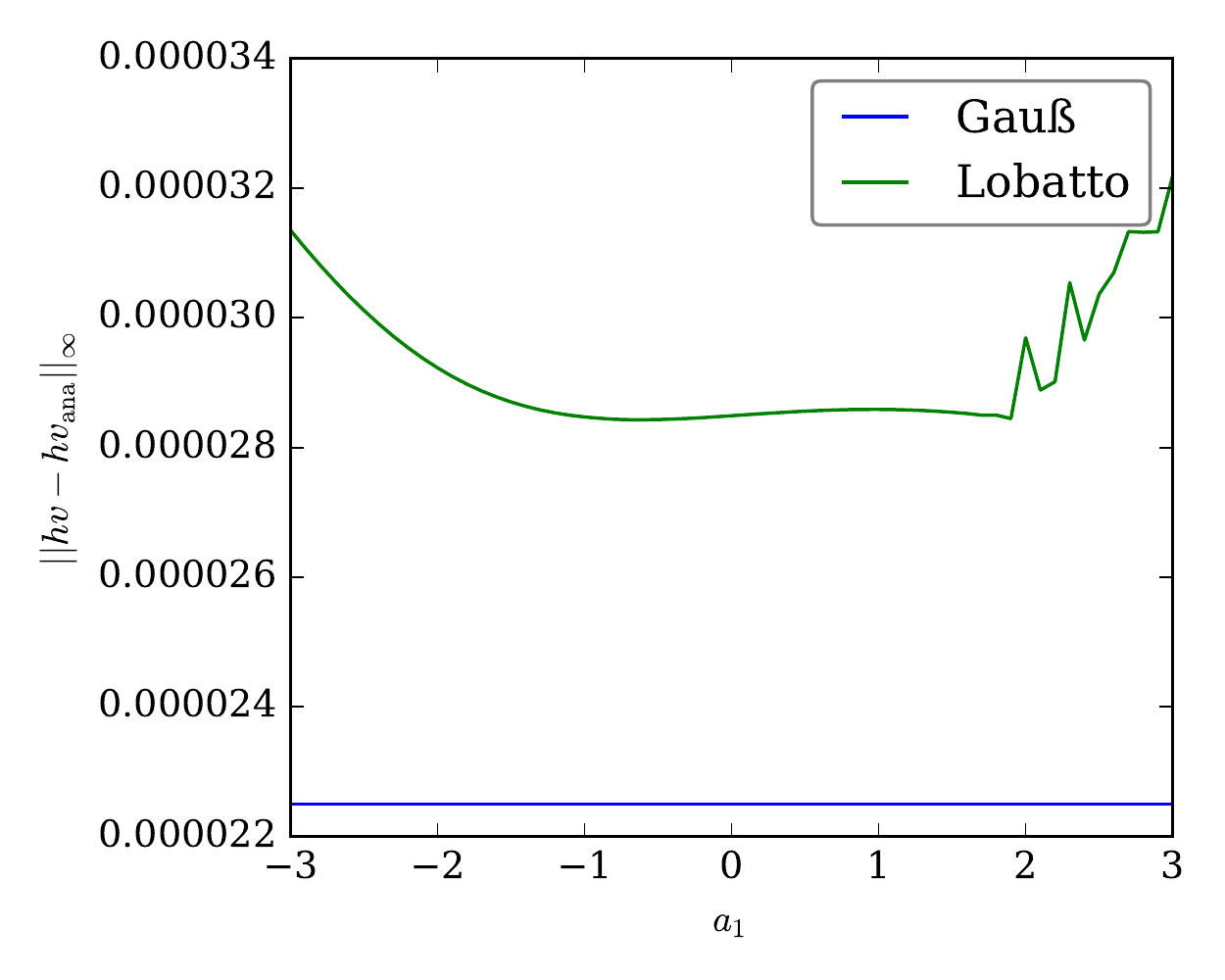}
    \caption{Error $\norm{hv - hv_\mathrm{ana}}_\infty$.}
  \end{subfigure}%
  \caption{Errors of the numerical solution for varying parameter $a_1$ and
           $N = 100$ elements of polynomials of degree $\leq p = 2$ using Gauß and
           Lobatto nodes.}          
  \label{fig:NUM_TEST_dam_breaks_p_2_N_100}
\end{figure}

Motivated by the result of section \ref{sec:NUM_TEST_moving_water}, only the
parameter $a_1$ has been varied for this problem, while the parameter $a_2$ is
fixed at $a_2 = \frac{2-a_1}{3}$. the results for $a_1 \in \set{-3 + \frac{k}{10},\;
0 \leq k \leq 60}$ are shown in Figures \ref{fig:NUM_TEST_dam_breaks_p_2_N_50}
and \ref{fig:NUM_TEST_dam_breaks_p_2_N_100}.
There, the $L_2$ errors have been computed exactly for the polynomials using Gauß
nodes and the $\norm{\cdot}_\infty$ errors are computed at the same nodes.

In these experiments, Gauß nodes yield a lower error in the solutions, both in
$\norm{\cdot}_{L_2}$ and $\norm{\cdot}_\infty$ and this error is nearly independent
of the parameter $a_1$ (it varies at most three orders of magnitude lower).
However, the error using Lobatto nodes are influenced by the choice of $a_1$
with variations up to \si{50}{\%}.

The corresponding errors in both norms $\norm{\cdot}_{L_2}$ and $\norm{\cdot}_\infty$
follow approximately the same trend for $h$ and $hv$, respectively, but there are
differences between the error curves of the height $h$ and the discharge $hv$.

These results, especially the ones for the discharge $hv$, suggest, that choosing
the parameter $a_1$ between $-1$ and $0$ might be optimal, but this has to be
investigated thoroughly.

\begin{remark}
  In the research code used for the simulations, the implementation has not been
  optimised for maximal performance or adapted to the special nodes. Instead,
  simple matrix vector products have been evaluated. However, if all possibilities
  for optimisation are used, the increased accuracy of Gauß nodes has to be
  compared to the increased computational performance of Lobatto nodes. Thus,
  especially for two-dimensional and unstructured grids, the choice of Lobatto
  nodes as used by \citet{wintermeyer2016entropy} seems to be superior.
  However, the influence of the parameter $a_1$ in this case would have to be
  investigated.
\end{remark}

\section{Summary and conclusions}
\label{sec:summary}

A new two-parameter family of entropy stable and well-balanced numerical fluxes and
corresponding split forms with adapted surface terms for general SBP bases
including Lobatto and Gauß nodes has been developed. The positivity preserving
framework of \citet{zhang2011maximum} can be used in this setting, but has to
be accompanied by some additional dissipation / stabilisation mechanism near
wet-dry fronts. Here, the subcell finite volume framework has been used and
extended naturally to diagonal-norm nodal SBP bases.

Numerical tests confirm the properties of the derived schemes. As suggested by
a first physicists intuition, the second parameter of the two-parameter family
should be chosen as $a_1 = \frac{2-a_1}{3}$ in order not to use some higher order
terms in the velocity $v$. This choice has been advantageous for the considered
moving water equilibrium in section \ref{sec:NUM_TEST_moving_water}.

However, the choice of the first parameter $a_1$ does not seem to be similarly
simple. There is no clear physical intuition at first and the dam break experiments
in section \ref{sec:NUM_TEST_dam_break} are not unambiguous. Thus, further
analytical and numerical studies have to be performed in order to understand
the influence of this parameter and possible optimal choices.

However, since the additional correction terms allowing the use of Gauß nodes
become more and more complicated, the gain in accuracy does not seem to justify
the use of these. Therefore, the computationally efficient flux differencing form
using Lobatto nodes seems to be advantageous. However, this does not confine
the results about positivity preservation and both entropy stability and
well-balancedness of the new two-parameter family of numerical fluxes. These
can be expected to be extendable to two-dimensional unstructured and curvilinear
grids using tensor product bases on quadrilaterals  similarly to
\citet{wintermeyer2016entropy}.

Additional topics of further research include the investigation of interactions
of curved elements with the parameter $a_1$, of other means performing finite
volume subcell projection, and other stabilisation techniques.

\appendix

\section*{Acknowledgements}

The author would like to thank the anonymous reviewer for some helpful comments,
resulting in an improved presentation of this material.

\bibliographystyle{spbasic}      % basic style, author-year citations
\bibliography{references_abbr}

\end{document}